\documentclass{birkbook}

\usepackage{amsfonts}
\usepackage{amsmath,  color, multicol, lscape, rotating, multirow}
\usepackage{amstext,amscd}
\usepackage{amssymb,txfonts}
\usepackage{epsfig, psfrag}

\def\D{\mathbb{D}}
\def\Z{\mathbb{Z}}
\def\N{\mathbb{N}}
\def\Naturals{\mathbb{N}}

\def\ni{\noindent}

\def\Area{\hbox{\rm Area}}
\def\Cay{\hbox{\it Cay}}

\def\Diam{\hbox{\rm Diam}}
\def\EDiam{\hbox{\rm EDiam}}
\def\IDiam{\hbox{\rm IDiam}}

\def\Rad{\hbox{\rm Rad}}

\def\Vol{\hbox{\rm Vol}}
\def\FVol{\hbox{\rm FVol}}
\def\FDiam{\hbox{\rm FDiam}}

\def\Star{\hbox{\rm Star}}

\def\CAT{\hbox{\rm CAT}}

\def\R{\hbox{\rm R}}
\def\FL{\hbox{\rm FL}}
\def\FFL{\hbox{\rm FFL}}
\def\F+L{\hbox{$\textup{F}\!_+\textup{L}$}}
\def\GL{\hbox{\rm GL}}
\def\DGL{\hbox{\rm DGL}}

\def\ssm{\smallsetminus}
\def\SL{\hbox{\rm SL}}

\def\ms{\medskip}

\def\onto{{\kern3pt\to\kern-8pt\to\kern3pt}}

\def\<{\langle}
\def\>{\rangle}
\def\|{{\ |\ }}

\def\s{\sigma}

\newcommand{\Cone}{\textup{Cone}_{\omega}\!\left(\mathcal{X},\textbf{e},\textbf{s}\right)}

\def\XX{\mathcal X}

\def\TT{\mathcal T}

 \def\CC{\mathcal C}
 
 \def\SS{\mathcal S}

 \def\PP{\mathcal{P}}
 
\def\QQ{\mathcal Q}

\def\s{\sigma}
\def\e{\mathbf{e}}
\def\x{\mathbf{x}}
\def\y{\mathbf{y}}

\newcommand{\set}[1]{\left\{#1\right\}}
\newcommand{\restricted}[1]{\left|_{#1} \right.}
\newcommand{\abs}[1]{\left|#1\right|}
\renewcommand{\ni}{\noindent}

\renewcommand{\ms}{\medskip}
\newcommand{\bs}{\bigskip}
\def\*{^{\ast}}

\def\newqed{{\ifhmode\unskip\nobreak\hfil\penalty50 \hskip1em \else\nobreak\fi
   \mbox{}\nobreak\hfil $\Box$%
    \parfillskip=0pt \finalhyphendemerits=0 \par}}

\newcommand{\tc}[2]{\textcolor{#1}{#2}}
\definecolor{dmagenta}{rgb}{.5,0,.5} \newcommand{\dmagenta}[1]{\tc{dmagenta}{#1}}
\definecolor{dred}{rgb}{.5,0,0} 
\definecolor{dgreen}{rgb}{0,.5,0} 
\definecolor{dblue}{rgb}{0,0,0.5} 
\definecolor{black}{rgb}{0,0,0} 
\definecolor{vdgreen}{rgb}{0,.3,0} 
\definecolor{vdred}{rgb}{.3,0,0}

 \newtheorem{thm}{Theorem}[section]
 \newtheorem{cor}[thm]{Corollary}
 \newtheorem{lem}[thm]{Lemma}
 \newtheorem{prop}[thm]{Proposition}
 \theoremstyle{definition}
 \newtheorem{defn}[thm]{Definition}
  \newtheorem{Conjecture}[thm]{Conjecture}
 \theoremstyle{remark}
 \newtheorem{rem}[thm]{Remark}
  \newtheorem{exer}[thm]{Exercise}
 \newtheorem{Open Question}[thm]{Open Question}
 \newtheorem{Question}[thm]{Question}
  \newtheorem{Open Problem}[thm]{Open Problem}
  \newtheorem{Open Questions}[thm]{Open Questions}
 \newtheorem{ex}[thm]{Example}
 \newtheorem{examples}[thm]{Examples}
 \newtheorem{example}[thm]{Example}
\begin{document}
%
\pagenumbering{roman}
\title{\textbf{Filling Functions} \\  \bs \bs Notes for an advanced course on \\  \emph{The Geometry of the Word Problem \\ for Finitely Generated Groups} \\ \bs Centre de Recerca Mathem\`atica \\ Barcelona }
\author{\Large{T.R.Riley}}
\date{\Large{July 2005 \\ Revised February 2006}} 
\maketitle

\hspace*{2cm} \thispagestyle{empty} \newpage

\setcounter{page}{5}
\tableofcontents



\newpage
\section*{Notation}
\addcontentsline{toc}{chapter}{Notation}
\markboth{Notation}{Notation}\thispagestyle{empty}

\begin{tabbing}
aaaaaaaaaaaaaaaaaaaaaaaaaaa  \= aaaaaaaaaaaaaaaaaaaaaaaaaaaaaaaa   \kill \\ 
$\preceq,  \ \succeq,  \ \simeq$ \> $\parbox{8.5cm}{$f,g : [0,\infty) \to [0,\infty)$ satisfy $ f \preceq g$ when there exists $C>0$ such that $f(n) \leq Cg(Cn+C) + Cn +C$ for all $n$, satisfy $ f \succeq g$ when $ g \preceq f$, and satisfy $ f \simeq g$ when $ f \preceq g$ and $g \preceq f$.  These relations are extended to functions $f: \N \to \N$ by considering such $f$ to be constant on the intervals $[n,n+1)$.}  $ \\
$a^b$, \ $a^{-b}$, \ $[a,b]$ \> \rule{0cm}{.5cm}$b^{-1}ab$, \ $b^{-1}a^{-1}b$, \ $a^{-1}b^{-1}ab$  \\ 
$\Cay^1(G, X)$ \> \parbox{8cm}{the Cayley graph of $G$ with respect to a generating set $X$} \\
$\Cay^2(\PP)$ \> \parbox{8cm}{the Cayley 2-complex of a presentation $\PP$} \\
$\mathbb{D}^n$ \> \parbox{8.5cm}{the $n$-disc $\set{(x_1, \ldots, x_n) \in \mathbb{R}^n \mid \sum_{i=1}^n {x_i}^2 \leq 1}$} \\   
$\Diam(\Gamma)$ \> \parbox{8cm}{$\max\set{ \, \rho(a,b) \, \mid \, \textup{vertices } a,b \textup{ in } \Gamma}$, where $\rho$ is the combinatorial metric on a finite connected graph $\Gamma$} \\  
$d_{X}$ \> the word metric with respect to a generating set $X$ \\
$\ell(w)$ \> word length; i.e.\ the number of letters in the word $w$ \\
$\ell(\partial \Delta)$ \> the length of the boundary circuit of $\Delta$ \\
$\mathbb{N}, \, \mathbb{R}, \, \mathbb{Z}$ \> the natural numbers, real numbers, and  integers \\
$R^{-1}$ \> $\set{r^{-1} \mid r \in R}$, the inverses of the words in $R$ \\
$ \mathbb{S}^n$ \> \parbox{8.5cm}{the $n$-sphere, $\set{(x_1, \ldots, x^{n+1}) \in \mathbb{R}^{n+1} \mid \sum_{i=1}^{n+1} = 1}$} \\   
$\textup{Star}(\Delta_0)$ \> \parbox{8.5cm}{ for a subcomplex $\Delta_0 \subseteq \Delta$, the union of all closed cells in $\Delta$ that have non-empty intersection with $\Delta_0$} \\
$(T, T^{\ast})$ \> a dual pair of spanning trees -- see Section~\ref{via vK diagrams}  \\
$w^{-1}$ \> the inverse ${x_n}^{-\varepsilon_n}\ldots{x_2}^{-\varepsilon_2}{x_1}^{-\varepsilon_1}$ of a word $w =   {x_1}^{\varepsilon_1}{x_2}^{\varepsilon_2}\ldots{x_n}^{\varepsilon_n}$. \\
$\langle X \mid R \rangle$ \> \parbox{8.5cm}{the presentation with generators $X$ and defining relations (or \emph{relators}) $R$} \\
$X^{-1}$ \> \parbox{8cm}{$\set{ x^{-1} \mid x \in X}$, the formal inverses $x^{-1}$ of letters $x$ in an alphabet $X$} \\
$(X \cup X^{-1})^{\ast}$ \> \parbox{8cm}{the free monoid (i.e.\ the words) on $X \cup X^{-1}$}  \\
$\Delta$ \> a van~Kampen diagram \\ 
$\varepsilon$ \> the empty word  \\

\ni \emph{Diagram measurements} \rule{0cm}{1cm} \> \\ 

\ms

$\Area(\Delta)$  \rule{0mm}{7mm} \> the number of 2-cells in $\Delta$ \\
$\DGL(\Delta)$ \>  $\min\set{ \ \Diam(T) + \Diam(T^{\ast}) \  \left| \  T \textup{ a spanning tree in } \Delta^{(1)} \right. \ }$  \\
$\EDiam(\Delta)$ \> the diameter of $\Delta$ as measured in the Cayley 2-complex \\
$\FL(\Delta)$ \> the filling length of $\Delta$ -- see Section~\ref{via vK diagrams} \\ 
$\GL(\Delta)$ \> the diameter of the 1-skeleton of the dual of $\Delta$ \\ 
$\IDiam(\Delta)$ \> the diameter of the 1-skeleton of $\Delta$ \\
$\Rad(\Delta)$ \> \parbox{8.5cm}{$\max \set{ \ \rho(a, \partial \Delta) \  \mid \  \textup{ vertices } a \textup{ of } \Delta \  }$ as measured in $\Delta^{(1)}$} \\

\emph{Filling functions}  \rule{0cm}{1cm} \> \\ 

$\Area : \N \to \N$  \rule{0mm}{7mm} \> the Dehn function \\
$\DGL : \N \to \N$ \> the simultaneous diameter and gallery length function \\
$\EDiam : \N \to \N$ \> the extrinsic diameter function \\ 
$\FL : \N \to \N$ \> the filling length function \\ 
$\GL : \N \to \N$ \> the gallery length function \\ 
$\IDiam : \N \to \N$ \> the intrinsic diameter function  \\
$\overline{\IDiam} : \N \to \N$ \> the upper intrinsic diameter function -- see Section~\ref{det}  \\
$\Rad : \N \to \N$ \> the radius function  \\
$\overline{\Rad} : \N \to \N$ \> the upper radius function -- see Section~\ref{hyperbolic section}  

\end{tabbing}


\chapter{Introduction}\thispagestyle{empty}
\pagenumbering{arabic}\setcounter{page}{1}

The Word Problem was posed by Dehn \cite{Dehn2} in 1912.  He asked, given a group, for a systematic method (in modern terms, an \emph{algorithm}) which, given a finite list (a \emph{word}) of basic group elements (\emph{generators} and their formal inverses), declares whether or not their product is the identity.  One of the great achievements of $20^{\textup{th}}$ century mathematics was the construction by Boone~\cite{Boone} and Novikov~\cite{Novikov} of finitely presentable groups for which no such algorithm can exist.  However, the Word Problem transcends its origins in group theory and rises from defeat at the hands of decidability and complexity theory, to form a bridge to geometry -- to the world of isoperimetry and curvature, local and large-scale invariants, as brought to light most vividly by Gromov \cite{Gromov}.

So where does geometry enter?  Groups act: given a  group, one seeks a space on which it acts in as favourable a manner as possible, so that the group can be regarded as a discrete approximation to the space.  Then a dialogue between geometry and algebra begins.  And where can we find a reliable source of such spaces?  Well, assume we have a finitely generating set $X$ for our group $G$. (All the groups in this study will be finitely generated.)  For $x,y \in G$, define the distance $d_{X}(x,y)$ in the \emph{word metric} $d_{X}$ to be the length of the shortest word in the generators and their formal inverses that represents $x^{-1}y$ in $G$.    
Then $$d_{X}(zx,zy) = d_{X}(x,y)$$ for all $x,y,z \in G$, and so left multiplication is action of $G$ on $(G, d_{X})$ by isometries.  

However $(G, d_{X})$ is discrete and so appears geometrically emaciated (``boring and uneventful to a geometer's eye'' in the words of Gromov~\cite{Gromov}).   Inserting a directed edge labelled by $a$ from $x$ to $y$ whenever $a \in X$ and $xa =y$ gives a skeletal structure known as the \emph{Cayley graph} $\Cay^1(G, X)$.  If $G$ is given by a  finite presentation $\PP = \langle X \mid R \rangle$ we can go further: attach flesh to the dry bones of $\Cay^1(G, X)$ in the form of 2-cells, their boundary circuits glued along edge-loops around which read words in $R$.  The result is a simply connected space  $\Cay^2(\PP)$ on which $G$ acts \emph{geometrically} (that is, properly, discontinuously and cocompactly) that is known as the \emph{Cayley 2-complex} -- see Section~\ref{vK_diagrams} for a precise definition.  
Further enhancements may be possible.  For example, one could seek to attach cells of dimension $3$ or above to kill off higher homotopy groups or so that the complex enjoys curvature conditions such as the $\CAT(0)$ property (see \cite{BradyNotes} or \cite{BrH}).  

Not content with the combinatorial world of complexes, one might seek continuous or smooth models for $G$.  For example, one could realise $G$ as the fundamental group of a closed manifold $M$, and then $G$ would act geometrically on the universal cover $\widetilde{M}$.  (If $G$ is finitely presentable then $M$ can be taken to be four dimensional  -- see \cite[A.3]{Bridson6}.)  Wilder non-discrete spaces, \emph{asymptotic cones}, arise from viewing $(G, d_{X})$ from increasingly distant vantage points (\emph{i.e.}\ scaling the metric to $d_{X} /s_n$ for some sequence of reals with $s_n \to \infty$) and recording recurring patterns using the magic of a non-principal ultrafilter.  Asymptotic cones see only some large-scale features of $(G, d_{X})$; they are the subject of Chapter~\ref{cones}.   

\bs

\ni \emph{Filling functions}, the subject of this study, capture features of discs spanning loops in spaces.  The best known is the classical \emph{isoperimetric} function for Euclidean space $\mathbb{E}^m$ -- any loop of length $\ell$ can be filled with a disc of area at most a constant times $\ell^2$.   To hint at how filling functions enter the world of discrete groups we mention a related algebraic result concerning the group $\Z^m$, the integer lattice in $m$-dimensional Euclidean space, generated by $x_1, \ldots, x_m$.  If $w$ is a word of length $n$ on $\set{ {x_1}^{\pm 1}, \ldots, {x_m}^{\pm 1}}$ and $w$ represents the identity in $\Z^m$ then, by cancelling off pairs ${x_i}{x_i}^{- 1}$ and  ${x_i}^{- 1}{x_i}$, and by interchanging adjacent letters at most  $n^2$ times, $w$ can be reduced to the empty word.   

This qualitative agreement between the number of times the commutator relations $x_ix_j=x_jx_i$ are applied and the area of fillings is no coincidence; such a relationship holds for all finitely presented groups, as will be spelt out in Theorem~\ref{Filling Thm} (\emph{The Filling Theorem}).  
The bridge between continuous maps of discs filling loops in spaces and this computational analysis of reducing words is provided by \emph{van~Kampen} (or \emph{Dehn}) \emph{diagrams}.  The Cayley 2-complex of the presentation $$\PP \ := \  \langle x_1, \ldots, x_m \mid [x_i,x_j], \forall i,j \in \set{1, \ldots, m} \, \rangle$$ of $\Z^m$ is the 2-skeleton of the standard decomposition of $\mathbb{E}^m$ into an infinite array of $m$-dimensional unit cubes.   A word $w$ that represents $1$ in $\PP$ (or, indeed, in any finite presentation $\PP$) corresponds to an edge-loop in $\Cay^2(\PP)$.  As Cayley 2-complexes are simply connected such edge-loops can be spanned by filling discs and, in this combinatorial setting, it is possible and appropriate to take these homotopy discs to be combinatorial maps of planar 2-complexes homeomorphic to (possibly \emph{singular}) 2-discs into $\Cay^2(\PP)$.  A \emph{van Kampen diagram} for $w$ is a graphical demonstration of how it is a consequence of the relations $R$ that $w$ represents $1$; Figure~\ref{vK diag in Z^3} is an example. So the Word Problem amounts to determining whether or not a word admits a van~Kampen diagram.  (See \emph{van~Kampen's Lemma: Lemma}~\ref{vK Lemma}.)

Filling functions for finite presentations of groups (defined in Chapter~\ref{filling functions}) record geometric features of van~Kampen diagrams.  The best known is the \emph{Dehn function} (or \emph{minimal isoperimetric function}) of Madlener \& Otto~\cite{MO} and Gersten~\cite{Gersten4};  it concerns \emph{area} -- that is, number of 2-cells.  In the example of $\Z^m$ this equates to the number of times commutator relations have to be applied to reduce $w$ to the empty word -- in this sense (that is, \emph{in the Dehn proof system} -- see Section~\ref{Dehn proof system}) the Dehn function can also be understood as a non-deterministic \textsc{Time} complexity measure  of the Word Problem for $\PP$ -- see Section~\ref{Dehn proof system}.  The corresponding \textsc{Space} complexity measure is called the \emph{filling length function} of Gromov~\cite{Gromov}.  It has a geometric interpretation -- the filling length of a loop $\gamma$ is the infimal length $L$ such that $\gamma$ can be contracted down to its base vertex through loops of length at most $L$.  Other filling functions we will encounter include the \emph{gallery length}, and \emph{intrinsic} and \emph{extrinsic diameter functions}.  All are group invariants in that whilst they are defined with respect to specific finite presentations, their qualitative growth depends only on the underlying group;  moreover, they are quasi-isometry invariants, that is, qualitatively they depend only on the large-scale geometry of the group -- see Section~\ref{qi section} for details.

In Chapter~\ref{interrelationships} we examine the interplay between different filling functions -- this topic bares some analogy with the relationships that exist between different algorithmic complexity measures and, as with that field, many open questions remain.  The example of nilpotent groups discussed in Chapter~\ref{nilpotent} testifies to the value of simultaneously studying multiple filling functions.  Finally, in Chapter~\ref{cones}, we discuss how the geometry and topology of the asymptotic cones of a group $G$ relates to the filling functions of $G$. 

\bs
\ni 
\emph{Acknowledgements.}  These notes build on and complement \cite{Bridson6}, \cite{Gersten} and \cite[Chapter 5]{Gromov} as well as the other two sets of notes in this volume \cite{BradyNotes, ShortNotes}.  For Chapter~\ref{cones} I am particularly indebted to the writings of Dru\c{t}u \cite{Drutu2, Drutu, DrutuThesis, DrutuSurvey}.  This is not intended to be a balanced or complete survey of the literature, but rather is a brief tour heavily biased towards areas in which the author has been involved.  

If I have done any justice to this topic then the influence of  Martin~Bridson and Steve~Gersten should shine through.  I am grateful to them both for stimulating collaborations, for their encouragement, and for communicating their deep vision for the subject.  I thank Emina~Alibegovic, Will~Dison, Cornelia Dru\c{t}u Badea, Steve~Gersten, Mark~Sapir and Hamish Short for comments on earlier drafts.  I also thank the Centre de Recerca Matem\`atica in Barcelona for their hospitality during the writing of this study, Jos\'e Burillo and Enric Ventura for organising the workshop, and the NSF for partial support via grant  DMS--0540830.  

\bs 

\begin{flushright}
TRR
\end{flushright}

\chapter{Filling functions}\thispagestyle{empty}\label{filling functions}

\section{Van Kampen diagrams} \label{vK_diagrams}

\begin{figure}[ht]
\psfrag{a_1}{{$x_1$}}
\psfrag{a_2}{{$x_2$}}
\psfrag{a_3}{{$x_3$}}
\psfrag{a_4}{{$x_4$}}
\psfrag{a_m}{{$x_m$}}
\psfrag{r_1}{{$r_1$}}
\psfrag{r_2}{{$r_2$}}
\psfrag{r_3}{{$r_3$}}
\psfrag{r_4}{{$r_4$}}
\psfrag{r_n}{{$r_n$}}
\centerline{\epsfig{file=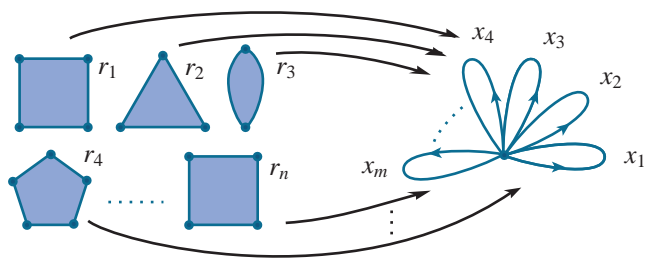}} 
\caption{The presentation 2-complex for $\langle x_1, \ldots, x_m \mid r_1, \ldots, r_n \rangle$}
\label{pres 2-complex}
\end{figure}

\begin{figure}[ht]
\psfrag{a}{{$a$}}
\psfrag{b}{{$b$}}
\centerline{\epsfig{file=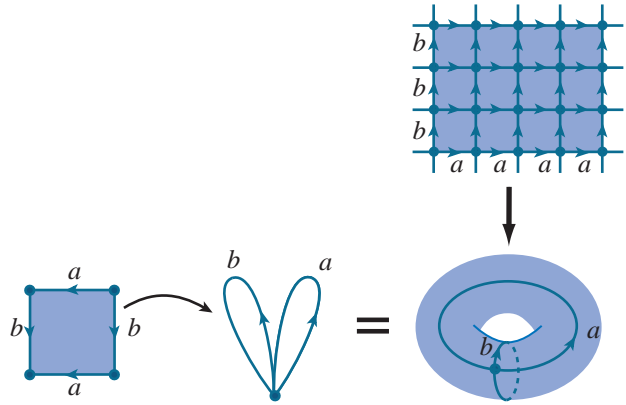}} 
\caption{The presentation and Cayley 2-complexes for $\langle a,b \mid [a,b] \rangle$} \label{Z^2 example}
\end{figure}

The \emph{presentation 2-complex} of $$\PP \ =  \ \langle X \mid R \rangle \  = \ \langle x_1, \ldots, x_m \mid r_1, \ldots, r_n \rangle$$ is constructed as shown in Figure~\ref{pres 2-complex}: take $m$ oriented edges, labelled by $x_1, \ldots, x_m$, identify all the vertices to form a  \emph{rose}, and then attach 2-cells $C_1, \ldots, C_n$, where $C_i$ has $\ell(r_i)$ edges, by identifying the boundary circuit of $C_i$ with the edge-path in the rose along which one reads $r_i$.  The \emph{Cayley 2-complex} $\Cay^2(\PP)$ is the universal cover of the presentation 2-complex of $\PP$.  The example of the free abelian group of rank 2, presented by $\langle a,b \mid [a,b] \rangle$ is shown in Figure~\ref{Z^2 example}.
 
\ni The edges of $\Cay^2(\PP)$ inherit labels and orientations from the presentation 2-complex and the 1-skeleton of $\Cay^2(\PP)$ is the Cayley graph  $\Cay^1(\PP)$ (cf.\ Definition~1.3 in \cite{ShortNotes}).  Identifying the group $G$ presented by $\PP$ with the 0-skeleton of $\Cay^1(\PP)$,  the path metric on $\Cay^1(\PP)$ in which each edge has length $1$  agrees with the word metric $d_{X}$ on $G$.

A word $w \in (X \cup X^{-1})^{\ast}$ is \emph{null-homotopic} when it represents the identity.  To such a $w$ one can associate an edge-circuit $\eta_w$ in $\Cay^2(\PP)$ based at some (and hence any) fixed vertex $v$, so that around $\eta_w$, starting from $v$, one reads $w$. A $\PP$-\emph{van~Kampen diagram} for $w$ is a combinatorial map $\pi : \Delta \to \Cay^2(\PP)$, where $\Delta$ is  $S^2 \ssm e_{\infty}$ for some combinatorial 2-complex $S^2$ homeomorphic to the 2-sphere and some open 2-cell $e_{\infty}$ of $S^2$, and the edge-circuit $\pi\restricted{\partial \Delta}$, based at vertex $\star$ on $\partial \Delta$, agrees with $\eta_w$.  (A map $\mathcal{X} \to \mathcal{Y}$ between complexes is \emph{combinatorial} if for all $n$, it maps the interior of each $n$-cell in $\mathcal{X}$ homeomorphically onto the interior of an $n$-cell in $\mathcal{Y}$.)  The base vertex $\star$ of $\Delta$ should not be ignored; we will see it plays a crucial role.  Van~Kampen's Lemma (see Section~\ref{algebraically}) says,  in particular, that  $w \in (X \cup X^{-1})^{\ast}$ is null-homotopic if and only if it admits a $\PP$-van~Kampen diagram. 
 
\begin{figure}[ht]
\psfrag{a}{\footnotesize{$a$}}
\psfrag{b}{\footnotesize{$b$}}
\psfrag{c}{\footnotesize{$c$}}
\centerline{\epsfig{file=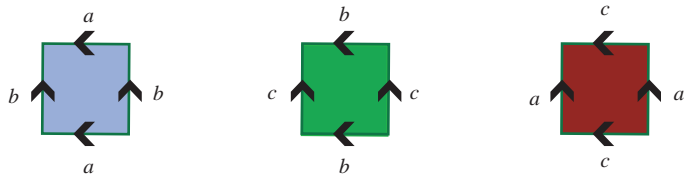}} 
\caption{The defining relations in   $\langle a, b, c \mid [a,b], [b,c], [c,a] \rangle$.}
\label{Z^3 2-cells}
\end{figure}

\begin{figure}[ht]
\psfrag{a}{\footnotesize{$a$}}
\psfrag{b}{\footnotesize{$b$}}
\psfrag{c}{\footnotesize{$c$}}
\psfrag{pi}{\footnotesize{$\pi$}}
\psfrag{s}{{$\star$}}
\centerline{\epsfig{file=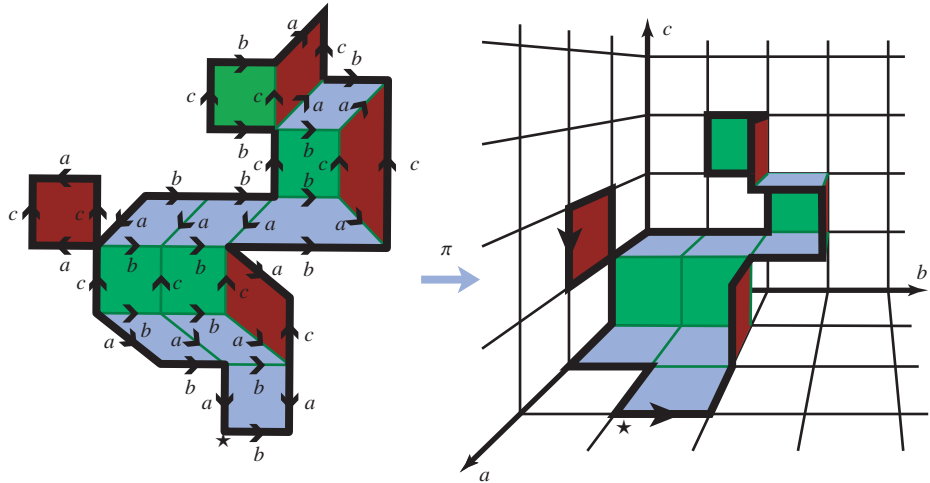}} 
\caption{A van~Kampen diagram for $ba^{-1}ca^{-1}bcb^{-1}ca^{-1}b^{-1}c^{-1}bc^{-1}b^{-2}acac^{-1}a^{-1}c^{-1}aba$ in $\langle a, b, c \mid [a,b], [b,c], [c,a] \rangle$.}
\label{vK diag in Z^3}
\end{figure}

It is convenient to have the following alternative definition  in which no explicit reference is made to $\pi$.  A $\PP$-van~Kampen diagram $\Delta$ for $w$ is a finite planar contractible combinatorial 2-complex with directed and labelled edges such that anti-clockwise around $\partial \Delta$ one reads $w$, and around the boundary of each 2-cell one reads a cyclic conjugate  (that is, a word obtained by cyclically permuting letters) of a word  in $R^{\pm 1}$.  From this point of view a $\PP$-van~Kampen diagram for $w$ is a filling of a planar edge-loop labelled by $w$ in the manner of a jigsaw-puzzle, where the pieces (such as those pictured in Figure~\ref{Z^3 2-cells}) correspond to defining relations and are required to be fitted together in such a way that orientations and labels match.
 The analogy breaks down in that the pieces may be distorted or flipped.  Figures~\ref{vK diag in Z^3}, \ref{BS fig}, \ref{rings} and \ref{Heisenberg fig} show examples of van~Kampen diagrams.    
 
The two definitions are, in effect, equivalent because, given the first, edges of $\Delta$ inherit
directions and labels from $\Cay^2(\PP)$, and given the second,  there is a combinatorial map $\pi : \Delta \to \Cay^2(\PP)$,  uniquely determined on the 1-skeleton of $\Delta$,  that sends $\star$ to $v$ and preserves the labels and directions of edges. 


\section{Filling functions via van~Kampen diagrams} \label{via vK diagrams}

Suppose $\PP = \langle X \mid R \rangle$ is a finite presentation for a group and $\pi : \Delta \to \Cay^2(\PP)$ is a van~Kampen diagram.    Being planar, $\Delta$ has a dual $\Delta^{\ast}$ -- the 2-complex with a vertex, edge and face dual to each face, edge and vertex in $\Delta$ (including a vertex dual to the face \emph{at infinity}, the complement of $\Delta$ in the plane).  

For a finite connected graph $\Gamma$ define $$\Diam(\Gamma) \  := \max\set{ \ \rho(a,b)  \ \mid \ \text{vertices } {a,b} \text{ of } \Gamma \ }$$ where $\rho$ is the combinatorial metric on $\Gamma$ -- the path metric in which each edge is given length $1$.     
Define the {\emph{area}},  
{\emph{intrinsic diameter}}, 
{\emph{extrinsic diameter}}, 
{\emph{gallery length}}, {\emph{filling length}} and $\DGL$ of $\pi: \Delta \to \Cay^2(\PP)$ by
\begin{eqnarray*}
\Area(\Delta) & = & \text{\# 2-cells in } \Delta \\
\IDiam(\Delta) & = & \Diam( \Delta^{(1)} ) \\  
\EDiam(\Delta) & = & \max \set{ \ d_{X}(\pi(a),\pi(b))  \ \mid \ \text{vertices } {a,b} \text{ of } \Delta \ } \\
\GL(\Delta) & = &  \Diam( {\Delta}^{\star(1)} ) \\ 
\FL(\Delta) & = & \min \set{ \ L \  \left| \  \exists \text{ a } \textit{shelling } (\Delta_i) \text{ of } \Delta  \text{ such that }  \max_i \ell(\partial \Delta_i) \leq L \right. \ } \\
\DGL(\Delta) & = & \min\set{ \ \Diam(T) + \Diam(T^{\ast}) \  \left| \  T \textup{ a spanning tree in } \Delta^{(1)} \right. \ }. 
\end{eqnarray*}

\ni 
These are collectively referred to as  \emph{diagram measurements}.
Note that $\IDiam$ measures diameter in the 1-skeleton of $\Delta$ and $\EDiam$ in the Cayley graph of $\PP$.  

\begin{figure}[ht]
\psfrag{1}{{{1-cell}}}%
\psfrag{2}{{{2-cell}}}%
\psfrag{f}{{{fragmentation}}}%
\psfrag{c}{{{collapse}}}%
\psfrag{e}{{{expansion}}}%
\psfrag{l(w_0)}{}
\psfrag{l(w_i)}{}
\psfrag{neq}{}
\psfrag{for}{}
\psfrag{s}{$\star$}
\psfrag{FL}{{{$\FL(\Delta)$} minimises {$\max_i\ell(w_i)$}:}}
\centerline{\epsfig{file=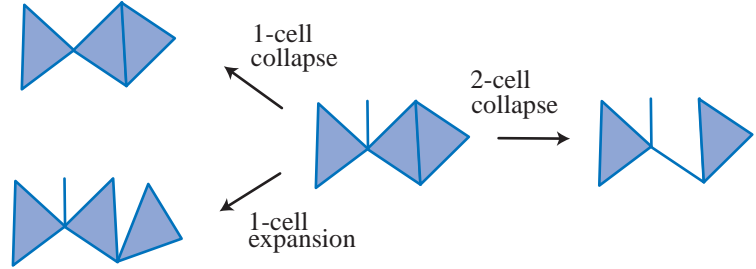}} 
\caption{Shelling moves} \label{shelling moves}
\end{figure}

The definitions of $\FL(\Delta)$ and $\DGL(\Delta)$ require further explanation.  
A \emph{shelling} of $\Delta$ is,  roughly speaking, a \emph{combinatorial null-homotopy} of $\Delta$ down to its base vertex $\star$.  More precisely, a shelling of $\Delta$ is sequence $\mathcal{S}= (\Delta_i)$ of diagrams $(\Delta_i)_{i=0}^m$ with $\Delta_0 = \Delta$ and  $\Delta_m = \star$ and such that  $\Delta_{i+1}$ is obtained from $\Delta_i$ by one of the \emph{shelling moves} defined below (illustrated in Figure~\ref{shelling moves}).  Moreover,  the base vertex $\star$ is \emph{preserved} throughout $(\Delta_i)$ -- that is, in every \emph{1-cell collapse} $e^0 \neq \star$, and in every \emph{1-cell expansion} on $\Delta_i$ where $e_0=\star$ a choice is made as to which of the two copies of $e_0$ is to be $\star$ in $\Delta_{i+1}$.

\begin{itemize}
\setlength{\itemsep}{0pt} \setlength{\parsep}{0pt}

\item \emph{1-cell collapse.} Remove a pair $(e^1,e^0)$ where
$e^1$ is a 1-cell with $e^0 \in \partial e^1$ and $e^1$ is
attached to the rest of $\Delta_i$ only by one of its end vertices
$\neq e^0$.  (Call such an $e^1$ a \emph{spike}.)

\item \emph{1-cell expansion.} Cut along some 1-cell $e^1$ in $\Delta_i$ that has a vertex $e^0$ in $\partial
\Delta_i$, in such a way that $e^0$ and $e^1$ are doubled.

\item \emph{2-cell collapse.} Remove a pair $(e^2,e^1)$ where
$e^2$ is a 2-cell which has some edge $e^1 \in (\partial e^2 \cap
\partial \Delta_i)$.  The effect on the boundary circuit is to replace
$e^1$ with $\partial e^2 \ssm e^1$.
\end{itemize}

\begin{figure}[ht]
\psfrag{D}{{\dmagenta{$\Delta$}}}
\psfrag{w_0}{{{$w_0$}}}
\psfrag{w_1}{{{$w_1$}}}
\psfrag{w_2}{{{$w_2$}}}
\psfrag{w_3}{{{$w_3$}}}
\psfrag{w_s}{{{$w_s$}}}
\psfrag{l(w_0)}{}
\psfrag{l(w_i)}{}
\psfrag{neq}{}
\psfrag{for}{}
\psfrag{s}{$\star$}
\psfrag{FL}{{\dmagenta{$\FL(\Delta)$} minimises \dmagenta{$\max_i\ell(w_i)$}:}}
\centerline{\epsfig{file=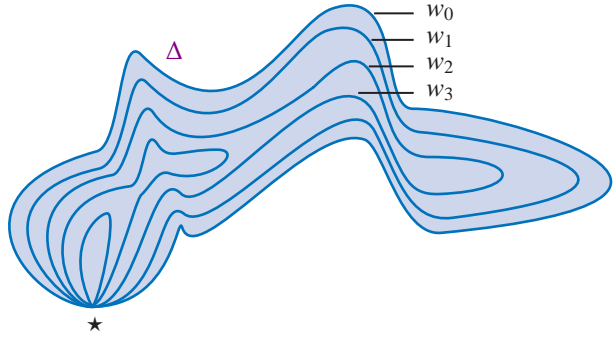}} 
\caption{The contracting loops in the course of a null-homotopy} \label{contracting loops}
\end{figure}

\ni There are natural combinatorial maps $\Delta_i \to \Delta$ whose restrictions to the interiors of $\Delta_i$ are injective, and which map the boundary circuits of $\Delta_i$, labelled by words $w_i$, to a sequence of edge-loops contracting to $\star$ as illustrated schematically in Figure~\ref{contracting loops}.

\ni Returning to the definition of  $\DGL(\Delta)$, if $T$ is a spanning tree in the 1-skeleton of $\Delta$ then define  $T^{\ast}$ to be the subgraph of the 1-skeleton of $\Delta^{\ast}$ consisting of all the edges dual to edges of $\Delta^{(1)} \ssm T$. The crucial property of $T^{\ast}$ is:  

\begin{exer}
$T^{\ast}$ is a spanning tree in the 1-skeleton of $\Delta^{\ast}$. 
\end{exer}

\ms

\ni We now define an assortment of \emph{filling functions} for  $\PP$. 
\begin{itemize}\setlength{\itemsep}{0pt} \setlength{\parsep}{0pt}
\item The  \textit{Dehn function} $\Area : \N \to \N$, 
\item the \textit{intrinsic diametric function} $\IDiam : \N \to
\N$, 
\item the \textit{extrinsic diametric function} $\EDiam : \N \to
\N$, 
\item  the \textit{gallery length function} $\GL : \N \to \N$, 
\item  the function $\DGL:\N \to\N$,
\item and the \textit{filling length function} $\FL : \N \to \N$
\end{itemize}
of $\PP$ are defined by 
\begin{eqnarray*}
\textup{M}(w) & := & \min \set{ \ \textup{M}(\Delta) \  \mid \ \Delta \textup{ a van Kampen diagram for } w \ }, \\ 
\textup{M}(n) & := & \max \set{ \ \textup{M}(w) \ \mid \
\text{words } w \text{ with } \ell(w) \leq n \text{ and } w=1 \text{ in } G \  }, 
\end{eqnarray*}
where
$\textup{M}$ is $\Area$, $\IDiam$, $\EDiam$, $\GL$, $\DGL$ and $\FL$
respectively. (The meanings of $\textup{M}(\Delta)$, $\textup{M}(w)$ and $\textup{M}(n)$ depend on their arguments: diagram, null-homotopic word, or natural number; the potential ambiguity is tolerated as it spares us from an overload of terminology.)

An \emph{isoperimetric function} (respectively, \emph{isodiametric function}) for $\PP$ is any $f : \N \to \N$ such that $\Area(n) \leq f(n)$ (respectively, $\IDiam(n) \leq f(n)$) for all $n$.  Chapter~5 of \cite{Gromov} is a foundational reference on filling functions.  
Many other references to isoperimetric functions, Dehn functions and isodiametric functions appear in the literature; \cite{Bridson6} and \cite{Gersten} are surveys.  Dehn functions were introduced by Madlener \& Otto~\cite{MO} and, independently, by Gersten~\cite{Gersten4}.  The filling length function of Gromov~\cite{Gromov} is discussed extensively in \cite{GR1};  a closely related notion $\textup{LNCH}$ was introduced by Gersten~\cite{Gersten5}. Gallery length and $\DGL$ were introduced in \cite{GR3}, and $\EDiam$ appears in \cite{BR1}.

\begin{defn}
We say that functions  $f_i:\N \to \N$ are \emph{simultaneously realisable} (upper or lower) bounds on a collection of filling functions $M_i: \N \to \N$  of $\PP$  (that is, each $M_i$ is one of $\Area$, $\FL$, $\IDiam$,\ldots) if for every null-homotopic word $w$, there exists a van~Kampen diagram $\Delta$ for $w$ such that $f_i(\ell(w))$ is at most  or at least (as appropriate) $M_i(\Delta)$, for all $i$.
\end{defn}

\ms

\ni The following two exercises are essentially elementary observations.  The second serves to describe
a combinatorial group theoretic adaptation from \cite{GR4} of a variant $\F+L$ of filling length defined by Gromov in \cite[page 101]{Gromov}.   $\F+L$ was used in \cite{GR4} to show that groups that enjoy Cannon's \emph{almost convexity condition} $\textup{AC}(2)$ have filling length functions (in the standard sense) growing $\preceq n$.  

\ms

\begin{exer}
Show that for a finite presentation, $d(u,v) := \Area(u^{-1}v)$ defines a metric on any set of reduced words, all representing the same group element.  
\end{exer}

\ms

\begin{exer}
Suppose $u,v\in (X^{\pm 1})^{\ast}$ represent the same group element and $\Delta$ is a  van~Kampen diagram for $uv^{-1}$ with two distinguished boundary vertices $\star_1, \star_2$ separating the $u$- and $v$-portions of the boundary circuit.  Define $\F+L(u,v,\Delta)$ to be the least $L$ such that there is a \emph{combinatorial homotopy of $u$ to $v$ across $\Delta$} through paths of length at most $L$ from  $\star_1$ to $\star_2$.  (More formally, such a combinatorial homotopy is a sequence of van~Kampen diagrams $(\Delta_i)_{i=0}^m$ with $\Delta_0=\Delta$ and $\Delta_m$ a simple edge-path along which one reads $w_2$, and such that $\Delta_{i+1}$ is obtained from $\Delta_i$ by \emph{1-cell collapse},  \emph{1-cell expansion}, or \emph{2-cell collapse}  in such a way that the $v$-portion of the boundary words $\partial \Delta_i$ is not broken.)  
Define $\F+L(u,v)$ to be the minimum of $\F+L(u,v,\Delta)$  over all van~Kampen diagrams $\Delta$ for $u{v}^{-1}$ 
and define $\F+L: \Naturals \to \Naturals$ by letting $\F+L(n)$ be the maximum of $\F+L(u,v)$ over all $u, v$ of length at most $n$ that represent the same element of the group.  Show that 
\begin{enumerate}
\item $\FL(2n) \leq \F+L(n)$ and $\FL(2n+1) \leq \F+L(n+1)$ for all $n$.
\item  $\F+L(u,v)$ defines a metric on any set of words that all represent the same group element.  
\end{enumerate}
\end{exer}

\section{Example: combable groups} \label{combable section}

\emph{Combable} groups form a large and well studied class that includes all automatic groups \cite{Epstein} and (hence) all hyperbolic groups.  A \emph{normal form} for a group $G$ with finite generating set $X$ is  a section $\s : G \to (X \cup
X^{-1})^{\ast}$ of the natural surjection $(X \cup X^{-1})^{\ast} \twoheadrightarrow G$. In other words, a normal form  is  a choice of representative $\s_g = \s(g)$ for each group element $g$.  View $\s_g$ as a continuous path $[0, \infty) \to \Cay^1(G, X)$ 
from the identity to $g$ (the ``\emph{combing line} of $g$''), travelling at unit speed from the identity until time $\ell(\s_g)$ when it halts for evermore at $g$. 

Following \cite{Bridson5, Gersten}, we say $\s$ \emph{synchronously
$k$-fellow-travels} when $$\forall g,h \in G, \ \
\left( d_{X}(g,h)=1 \implies \forall t \in \N, \ d_{X}(\, \s_g(t),\s_h(t)
\,  ) \leq k \right).$$
Define a \emph{reparametrisation} $\rho$ to be an unbounded
function $\Naturals \to \Naturals$ such that $\rho(0)=0$ and
$\rho(n+1) \in \set {\rho(n), \rho(n)+1 }$ for all $n$.  We say $\s$ \emph{asynchronously $k$-fellow travels} when for all $g,h \in G$ with $d_{X}(g,h)=1$, there
exist reparametrisations $\rho$ and $\rho'$ such that $$\forall
t \in \N, \ \ d_{X}( \, \s_g(\rho(t)),\s_h(\rho'(t)) \, ) \leq k.$$
(Note $\rho$ and $\rho'$ both depend on both $g$ and $h$.)

We say $(G,X)$ is (\emph{a})\emph{synchronously combable} when, for some $k \geq 0$, there is an (a)synchronous $k$-fellow-travelling normal  form $\s$ for $G$. Define the \emph{length function}  $\textup{L}: \Naturals \to \Naturals$ of  $\s$ by:
$$\textup{L}(n) \ := \  \max \set{\ \ell(\s_g) \  \mid  \  d_{X}(1, g) \leq n \ }.$$

\begin{examples} \label{combing examples} 
 \rule{0cm}{0cm} \\ \vspace*{-.5cm}
 \begin{enumerate}
\item \emph{Finite groups.} 
Suppose $G$ is a finite group with finite generating set $X$.  Let $k$ be the diameter of $\Cay^1(G, X)$.  If for all $g \in G$ we take $\sigma_g$ to be any word representing $g$ then $\sigma$ is synchronously $k$-fellow travelling.  Moreover, if $\sigma_g$ is a choice of geodesic (i.e.\ minimal length)  word representing $g$ then $\sigma$ has length function satisfying $\textup{L}(n) \leq \min \set{n,k}$ for all $n$.  

\item The (unique) geodesic words representing group elements in $F_m = \langle a_1, \ldots, a_m \mid \  \rangle$ form a synchronous 1-fellow-travelling combing.   The words $$\set{ \  {a_1}^{r_1} {a_2}^{r_2} \ldots {a_m}^{r_m} \mid r_1,\ldots, r_m \in \Z \ }$$ comprise a synchronous 2-fellow travelling combing of 
$$\Z^m = \langle a_1, \ldots , a_m \mid [a_i,a_j], \forall 1\leq i < j \leq m \rangle.$$  In both cases $L(n) = n$ for all $n$.    

\item  $\textup{BS}(1,2) = \langle a,b \mid b^{-1}ab = a^2 \rangle$  is one of a family of groups $\textup{BS}(m,n)$ with presentations $\langle a,b \mid b^{-1}a^mb = a^n \rangle$ that are named in honour of Baumslag and Solitar who studied them in \cite{BS}.   $\textup{BS}(1,2)$ is often referred to as  ``the Baumslag--Solitar group'' despite repeated public insistence from Baumslag that this is an inappropriate attribution.   

The words $$\set{ \ b^rua^s  \  \mid \ r,s \in \Z, \ u \in \set{ab^{-1},b^{-1}}^{\ast}, \ \textup{the first letter of } u \textup{ is not } b^{-1} \ } $$ define an asynchronous combing of $\textup{BS}(1,2)$.   In fact, the normal form element for $g \in \textup{BS}(1,2)$ results from applying the \emph{rewriting rules} 
\begin{tabbing}
 \hspace*{2.5cm} \= $a b$ \ \ \ \ \ \  \ \=   $\mapsto$ \ \  \= $ba^{2}$  \hspace*{2.5cm}  \= $aa^{-1} \ \mapsto \ \varepsilon$ \\
\> $a^{-1} b$ \>  $\mapsto$ \>  $ba^{-2}$ \> $a^{-1}a \  \mapsto \ \varepsilon$ \\
\> $a^2 b^{-1}$  \> $\mapsto$ \> $b^{-1}a$ \> $bb^{-1} \ \mapsto \ \varepsilon$ \\
\> $a^{-1} b^{-1}$ \>  $\mapsto$ \>  $ab^{-1}a^{-1}$ \> $b^{-1}b \  \mapsto \ \varepsilon$ 
\end{tabbing} 
to any word representing $g$.  See \cite[Chapter 7]{Epstein} for more details.

\begin{figure}[ht]
\psfrag{a}{$a$}
\psfrag{b}{$b$}
\psfrag{w}{$b^6ab^{-6}$}
\centerline{\epsfig{file=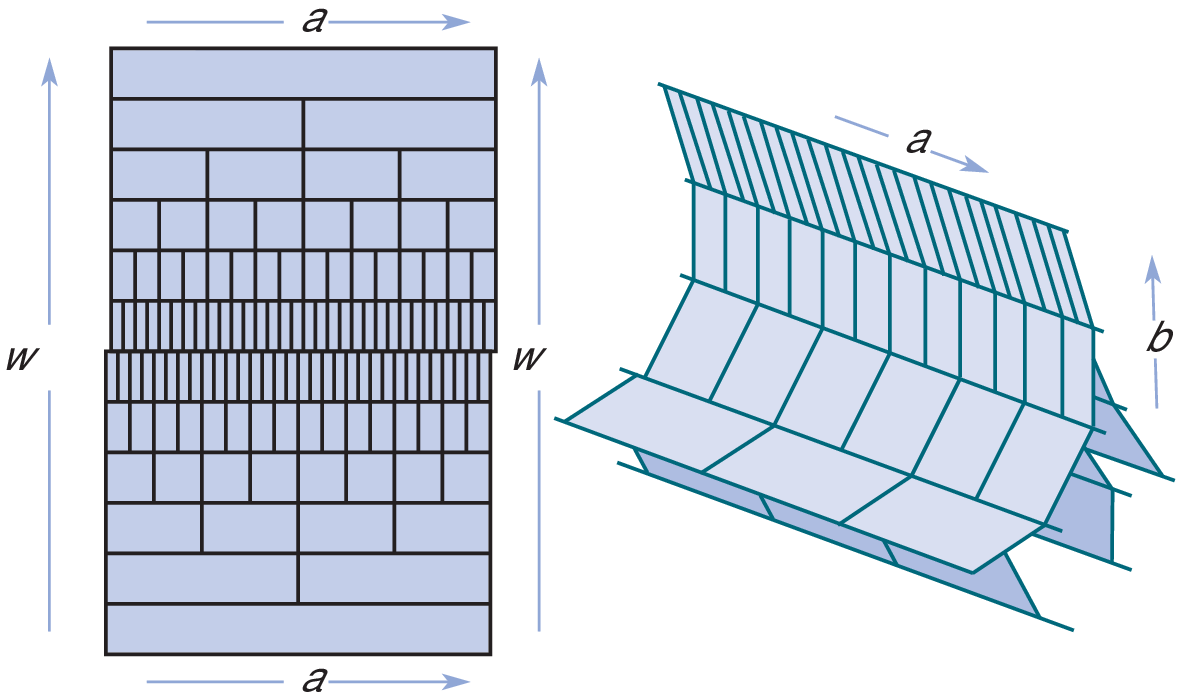}} 
\caption{A van~Kampen diagram for $\left[a , a^{b^6}\right]$ in $\langle a,b \mid b^{-1}ab = a^2 \rangle$ and a portion of $\Cay^2(\PP)$.} \label{BS fig}
\end{figure}

\item \emph{Automatic groups.} A group $G$ with finite generating set $X$ is \emph{synchronously automatic} when it admits a synchronous combing $\sigma$ such that the set of normal forms $\set{\sigma_g \mid g \in G}$ comprise a \emph{regular language}.   The foundational reference is \cite{Epstein} -- there, a group is defined to be synchronous automatic when its multiplication and inverse operations are governed by finite state automata, and agreement of the two definitions is Theorem 2.3.5.

The $\sigma_g$ are quasi-geodesics and so there is a constant $C$ such that $\textup{L}(n) \leq Cn +C$ for all $n$.  Synchronously automatic groups include all finitely generated virtually abelian groups, all hyperbolic groups, all the braid groups $B_n$, many 3-manifold groups, and the mapping class groups of closed surfaces with finitely many punctures -- see \cite{Epstein, Mosher}.

The larger class of \emph{asynchronously automatic groups}, defined in \cite{Epstein}, can be characterised (see Theorems~7.2.8 and 7.3.6 in \cite{Epstein}) as the asynchronously combable groups $G$ such that $\set{\sigma_g \mid g \in G}$ is regular and admits a function $D:\mathbb{N} \to \mathbb{N}$ such that for all $g \in G$, $r,s \geq 0$, $t \geq D(r)$, if $s+t \leq \ell(\sigma_g)$ then $d_X( \sigma_g(s),   \sigma_g(s+t)) >r$. Their length functions satisfy $\textup{L}(n) \leq C^n$ for some $C>0$  (see e.g.\ \cite[Lemma~2.3]{Bridson2}).  Thurston \cite[Section~7.4]{Epstein} showed that $\textup{BS}(1,2)$  is asynchronously automatic: the language of the asynchronous combing given above is regular but, on account of of the exponential Dehn function  we will establish in Proposition~\ref{BS Dehn fn}, $\textup{BS}(1,2)$ is not synchronously automatic.

For references of an introductory nature and for open questions see, e.g.\ \cite{BMSqus, Farb2, Gersten5, Gersten9, Ohshika}.

\item  \emph{$\CAT(0)$ groups.}  Paths in a Cayley graph of a $\CAT(0)$ group that run close to $\CAT(0)$ geodesics can be used to define a synchronous combing with length function $\textup{L}(n) \simeq n$.   See the proof of Proposition~1.6 in \cite[III.$\Gamma$]{BrH} for more details.  
\end{enumerate}
\end{examples}

\ni Further classes of groups can be defined by specifying other grammatical, length function or geometric constraints.  See, for example, \cite{Bridson8, Farb2, Gilman} for more details.

\begin{exer}
Check that the words in Example~\ref{combing examples}(3) define an asynchronous combing of $\textup{BS}(1,2)$ and show that its length function satisfies $\textup{L}(n) \simeq 2^n$.  Estimate the fellow-travelling constant.  
\end{exer}

\begin{exer}
Show that if $X$ and $X'$ are two finite generating sets for a group $G$ then $(G, X)$ is (a)synchronously combable if and only if  $(G, X')$ is.  Moreover, show that if $\textup{L}(n)$ is a length function for an (a)synchronous combing of  $(G, X)$ then there is  an (a)synchronous combing of  $(G, X')$ with length function $\textup{L}'(n)$ satisfying $\textup{L}(n) \simeq \textup{L}'(n)$.   
\end{exer}

\begin{thm}[cf.\ \textup{\cite{ Epstein, Gersten, GR3}}] \label{Combable groups}
Suppose a group $G$ with finite generating set $X$ admits a
combing $\s$ that  \textup{(}a\textup{)}synchronously $k$-fellow-travels. Let $\textup{L}: \Naturals \to \Naturals$
be the length function of $\s$.  Then there exists $C>0$ and a finite
presentation $\PP=\langle X \mid R \rangle$ for $G$ for which  
\begin{eqnarray*}
\Area(n) & \leq &  Cn \, \textup{L}(n) \ \  \text{ and  } \\ 
\EDiam(n) \ \leq \ \IDiam (n) \ \leq \  \FL (n) & \leq &  C\,n  
\end{eqnarray*}
for all $n$.  Moreover, these bounds are realisable simultaneously.
\end{thm}

\begin{proof} 
Suppose $w$ is a word representing $1$ in $G$.  We will construct the  \emph{cockleshell} van~Kampen diagram $\Delta$ for $w$, illustrated in Figure~\ref{cockleshell}.  We start with a circle $C$ in the plane subdivided into $\ell(w)$ edges, directed and labelled so that  anticlockwise from the base vertex $\star$ one reads $w$. Then we join each vertex $v$ to $\star$ by an edge-path labelled by the word $\sigma_{g_v}$, where $g_v$ is the group element represented by the prefix of the word $w$ read around $C$ between $\star$ and $v$.  For each pair $u, v$ of adjacent vertices on $C$, we use the fact that $\s_{g_u}$ and $\s_{g_v}$ asynchronously $k$-fellow travel (with respect to reparametrisations $\rho$ and $\rho'$)   to construct a \emph{ladder} whose \emph{rungs} are paths of length at most $k$; these would be straight if our combing was synchronous, but may be skew in the asynchronous case.  A rung connects $\s_{g_u}(\rho(t))$ to $\s_{g_v}(\rho'(t))$
for  $t =1,2,\ldots$.  Duplicate rungs may occur on account of $\rho$  and $\rho'$ simultaneously pausing, that is,  $\rho(m) = \rho(m+1)$ and $\rho'(m) = \rho'(m+1)$ for some $m$; counting each of these only once, there are a total of at most $\ell_{g_u} + \ell_{g_v}$ rungs and hence at most the same number of 2-cells in the ladder.   

Subdivide the rung from $\s_{g_u}(\rho(t))$ to $\s_{g_v}(\rho'(t))$ into $d_{X}(\s_{g_u}(\rho(t)), \s_{g_v}(\rho'(t))) \leq k$ edges and label it by a shortest word representing the same element  in $G$ as $[\s_{g_u}(\rho(t))]^{-1}\s_{g_v}(\rho'(t))$. 

\begin{figure}[ht]
\psfrag{D}{{{$\Delta$}}}
\psfrag{w}{{{$w$}}}
\psfrag{u}{{{$u$}}}
\psfrag{v}{{{$v$}}}
\psfrag{su}{{{$\sigma_{g_u}$}}}
\psfrag{sv}{{{$\sigma_{g_v}$}}}
\psfrag{s}{$\star$}
\centerline{\epsfig{file=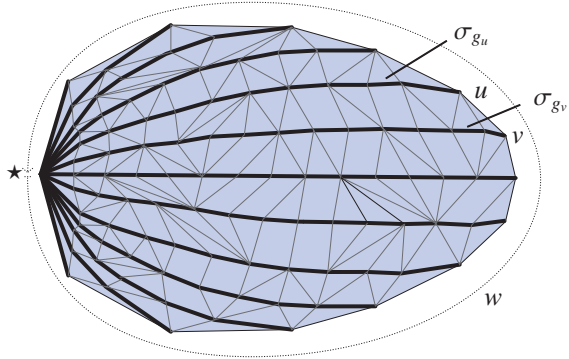}} 
\caption{The cockleshell diagram filling an edge-loop in an asynchronously combable group.} \label{cockleshell}
\end{figure}

Obtain $\Delta$ by regarding every simple edge-circuit in a ladder to be the boundary circuit of a 2-cell;
words around such edge-circuits have length at most $2k+2$ and represent the identity in $G$.  So, defining $R$ to be the set of null-homotopic words of length at most
$2k+2$, we find that  $\PP:=\langle X \mid R \rangle$ is a finite
presentation for $G$ and $\Delta$ is a $\PP$-van~Kampen diagram for
$w$.  We get the asserted bound on $\Area(n)$ from the observation 
$$\Area(\Delta) \ \leq \ 2 \, \ell(w) \, \textup{L}( \lfloor \ell(w)/2 \rfloor).$$

We will only sketch a proof of the linear upper bound on $\FL(n)$.  Consider paths $p$ connecting pairs of points on $\partial \Delta$, that run along successive rungs of adjacent ladders, no two rungs coming from the same ladder.  Observe that there is such a $p$ for which the region enclosed by $p$ and a portion of the the boundary circuit of $\Delta$ not passing through $\star$ is a topological 2-disc that includes no rungs in its interior.  Shell away this region using 1-cell and 2-cell collapse moves.  Repeating this process and performing 1-cells collapse moves will shell $\Delta$ down to $\star$  and the boundary loops of the intermediate diagrams will run (roughly speaking) perpendicular to the combing lines.

The inequalities $\EDiam(n) \leq \IDiam (n) \leq  \FL (n)$ are straightforward and are not specific to combable groups -- see Chapter~\ref{interrelationships}.
\end{proof}

\begin{rem}
Figure~\ref{cockleshell} suggests that $\DGL(\Delta)$ and (hence) $\GL(\Delta)$ are both at most a constant $C$ times $\textup{L}(\lfloor \ell(w)/2 \rfloor)$, taking $T$ (a spanning tree in the 1-skeleton of $\Delta$)  to be the union of all the combing lines, and hence that there is a further simultaneously realisable bound of 
\begin{eqnarray} \label{GL by L}
\GL(n) \ \leq \ \DGL(n) \ \leq \ C \, \textup{L}(n),
\end{eqnarray}
for some constant $C$.  But the picture is misleading -- some of the rungs could have zero length and hence $T$ could fail to be a tree.  
A way to circumvent this problem is to \emph{fatten} $\PP$, that is, introduce a spurious extra generator $z$, that will be trivial in the group and add, in a sense redundant, defining relations $z, z^2, z^3, zz^{-1}, z^2 z^{-1}$ and $[z,x]$ for all $x \in X$.  The combing lines can be kept apart in this new presentation and (\ref{GL by L}) then holds.  Details are in \cite{GR3}.
\end{rem}

\begin{rem} \label{Bridson combing}
It is tempting to think that an induction on word length would show that the length function of a $k$-fellow-travelling synchronous combing satisfies $L(n) \leq Cn$ for some constant $C=C(k)$.  But this is not so: Bridson \cite{Bridson8} gave the first example of a synchronously combable group with Dehn function growing faster than quadratic -- his example has cubic Dehn function.  In particular, on account of this Dehn function, his group is not automatic.
\end{rem}

\section{Filling functions interpreted algebraically} \label{algebraically}

By definition, when a group $G$ is presented by $\langle X \mid R \rangle$ there is a short exact sequence,  $\langle \! \langle R \rangle \! \rangle \hookrightarrow F(X) \onto G$, where $\langle \! \langle R \rangle \! \rangle$ is the normal closure of $R$ in $F(X)$.   So a word $w \in (X \cup X^{-1})^{\ast}$ represents $1$ in $G$ if and only if there is an equality
\begin{equation} \label{prod}
w \  = \  \prod_{i=1}^N {u_i}^{-1} {r_i}^{\epsilon_i} u_i
\end{equation} 
in $F(X)$ for some integer $N$, some $\epsilon_i = \pm 1$, some words $u_i$, and some $r_i \in R$.  Recent surveys including proofs of the following key lemma are \cite{Bridson6} (Theorem~4.2.2) and  \cite{ShortNotes} (Theorem~2.2).

\begin{lem}[van~Kampen's Lemma \cite{vK}] \label{vK Lemma}
A word $w$ is null-homotopic in a finite presentation $\PP$ if and only if it admits a $\PP$-van~Kampen diagram.  Moreover, $\Area(w)$ is the minimal $N$ such that there is an equality in $F(X)$ of the form \textup{(\ref{prod})}.  
\end{lem}

\ms

\begin{exer}
Adapt a proof of van~Kampen's Lemma to show that $$L \ \leq \ \IDiam(w)  \ \leq \ 2L +\max\set{\ell(r) \mid r \in \R} + \ell(w),$$ where $L$ is the minimal value of $\max_i\ell(u_i)$ amongst all equalities (\ref{prod}) in $F(X)$.  (The additive term $\ell(w)$ can be discarded when $w$ is freely reduced.)
\end{exer}

\ms

\begin{exer}
Similarly, relate $\EDiam(w)$ to the minimal value of $$\max \set{ \ d_X(1,w_0)  \ \left| \ \textup{ prefixes } w_0 \textup{ of the unreduced word } \prod_{i=1}^N {u_i}^{-1} {r_i}^{\epsilon_i} u_i
\right. \  }$$ amongst all equalities (\ref{prod}) in $F(X)$. 
\end{exer}



\section{Filling functions interpreted computationally} \label{Dehn proof system}

Given a finite presentation $\PP = \langle X \mid R \rangle$ and a word $w_0 \in (X \cup X^{-1})^{\ast}$ imagine writing $w_0$ on the tape of a Turing machine, one letter in each of $\ell(w_0)$ adjacent squares.  The moves described in the following lemma are ways of altering the word on the tape.  Using these moves with the aim of making all squares blank is the \emph{Dehn proof system} for trying to show $w$ represents the identity in $\PP$.  (The author recalls hearing this terminology in a talk by Razborov, but is unaware of its origins.)  We call $(w_i)_{i=0}^m$ a \emph{null-$\PP$-sequence} for $w_0$ (or a \emph{$\PP$-sequence} if we remove the requirement that the final word $w_m$ be the empty word).

\begin{lem}\label{comp lemma}
Suppose $w_0 \in (X \cup X^{-1})^{\ast}$.   Then $w_0$ is null-homotopic in $\PP = \langle X \mid R \rangle$ if and only if it can be reduced to the empty word $w_m$ via a sequence $(w_i)_{i=0}^m$  in which  each $w_{i+1}$ is obtained from $w_i$ by one of three moves:
\begin{enumerate}\setlength{\itemsep}{0pt} \setlength{\parsep}{0pt}
\item \emph{Free reduction:} $w_i = \alpha x x^{-1} \beta \mapsto \alpha \beta =w_{i+1}$, where $x \in X^{\pm 1}$.
\item \emph{Free expansion:} the inverse of  a free reduction.
\item \emph{Application of a relator:}  $w_i=\alpha u\beta \mapsto \alpha v\beta = w_{i+1}$, where a cyclic conjugate of $uv^{-1}$  is in $R^{\pm 1}$.
\end{enumerate}
\end{lem}

\begin{proof}
If $w_0$ is null-homotopic then, by Lemma~\ref{vK Lemma}, it admits a van~Kampen diagram $\Delta_0$.  Let $(\Delta_i)_{i=0}^m$ be any shelling of $\Delta_0$ down to its base vertex $\star$ by 1-cell- and 2-cell-collapse moves.  Let $w_i$ be the boundary word of $\Delta_i$ read from $\star$. Then $(w_i)_{i=0}^m$ is a null-sequence for $w_0$.
The converse is straight-forward.
\end{proof}

\ms

\begin{prop} \label{dehn proof prop}
Suppose $w_0$ is null-homotopic in $\PP$.  Then,  amongst all 
null-sequences $(w_i)$ for $w_0$,  
$\Area(w_0)$ is the minimum $A$ of the number of $i$ such that $w_{i+1}$ is obtained from $w_i$ by an \emph{application-of-a-relator} move, and 
$\FL(w_0)$ is the minimum of $\max_i \ell(w_i)$.
\end{prop}

\ms

\ni \emph{Sketch proof.}  
The proof of Lemma~\ref{comp lemma} shows that $A \leq \Area(w_0)$ and $F \leq \FL(w_0)$.  

For the reverse inequalities suppose $(w_i)$ is a null-sequence for $w_0$ and that $A$ and $F$ are the number of $i$ such that $w_{i+1}$ is obtained from $w_i$ by an \emph{application-of-a-relator} move, and of $\max_i \ell(w_i)$, respectively.  We seek a van~Kampen diagram $\Delta_0$ for $w_0$ such that $\Area(w_0) \leq A$ and 
$\FL(w_0) \leq F$.  We will describe a sequence of planar diagrams $A_i$ that topologically are singular annuli and have the property that the outer boundary of $A_i$ is labelled by $w_0$ and the inner boundary by $w_i$.  

Start with a planar, simple edge-circuit $A_0$ labelled by $w_0$.  Obtain $A_{i+1}$ from $A_{i}$ as follows.  If $w_i \mapsto w_{i+1}$ is a free reduction $w_i = \alpha x x^{-1} \beta \mapsto \alpha \beta =w_{i+1}$ then identify the two appropriate adjacent edges labelled $x$ and $x^{-1}$ in $A_i$.    If it is a free expansion $w_i = \alpha \beta \mapsto \alpha x x^{-1} \beta   =w_{i+1}$, then attach an edge labelled by $x$ by its initial vertex to the inner boundary circuit of  $A_i$.  If its is an application of a relator,  $w_i=\alpha u\beta \mapsto \alpha v\beta = w_{i+1}$, then attach a 2-cell with boundary labelled by $u^{-1}v$ along the $u$ portion of the inner boundary circuit of $A_i$.

However, the problem with this method of obtaining $A_{i+1}$ from $A_i$ is that there is no guarantee that $A_{i+1}$ will  be planar.  For example, if $w_i \mapsto w_{i+1}$ is a free expansion immediately reversed by a free reduction $w_{i+1} \mapsto w_{i+2}$  then an edge will be pinched off.  The resolution is that moves that give rise to such problems are redundant and $(w_i)$ can be amended to avoid the difficulties, and this can be done without increasing $A$ or $F$.   A careful treatment is in \cite{GR1}.   
\newqed

\bs

\ni The proof outlined above actually shows more: if $(w_i)$ is a null-sequence for $w$ then there is a van~Kampen diagram $\Delta$ for $w$ that simultaneously has area at most the number of application-of-a-relator moves in $(w_i)$ and has $\FL(\Delta) \leq \max_i \ell(w_i)$ .   

\ms

\begin{examples} \label{examples comp}  \rule{0cm}{0cm} \\ \vspace*{-.5cm}
\begin{enumerate} 
\item  \emph{Hyperbolic groups.} A  group $G$ with finite generating set $X$ is (\emph{Gromov}-) \emph{hyperbolic} if and only if $\Cay^1(G ,X)$ is hyperbolic in the sense we will define in Section~\ref{hyperbolic section}.  This condition turns out to be independent of the finite generating set $X$ (as follows from Theorem~\ref{R-trees}, for example).  Equivalently \cite[page~450]{BrH}, a group is hyperbolic if and only if it admits a \emph{Dehn presentation} -- a finite presentation $\PP$ such that if $w$ is null-homotopic then one can perform an application-of-a-relator move $w=\alpha u\beta \mapsto \alpha v\beta$ for which $\ell(v) < \ell(u)$.  It follows that $\Area(n), \FL(n) \leq n$, and these bounds are simultaneously realisable.     

So the Dehn and filling length functions of hyperbolic groups grow at most linearly.  

\item \emph{Finitely generated abelian groups.}  Suppose $G$ is a finitely generated abelian group $$\langle \ x_1, \ldots, x_m \ \mid \ {x_i}^{\alpha_i}, \ [x_i,x_j], \ \forall 1\leq i < j \leq m \ \rangle$$ where each $\alpha_i \in \set{0,1,2,\ldots}$.  Any word representing the identity can be reduced to the empty word by shuffling letters past each other to collect up powers of $x_i$, for all $i$,  and then freely reducing and applying the ${x_i}^{\alpha_i}$ relators.  So $\Area(n) \preceq n^2$ and $\FL(n) \leq n$ and these bounds are simultaneously realisable. 
\end{enumerate}
\end{examples}

\ms 

\ni Implementing the Dehn Proof System on a Turing machine shows the following .  
[A Turing machine is \emph{symmetric} when the transition relation is symmetric: the reverse of any transition is again a transition.] 

\begin{thm}[\textup{\cite[Theorem 1.1]{SBR}, \cite[Proposition 5.2]{Birget}}] \label{Turing machine for P} 
For every finite presentation $\PP$, there is a symmetric Turing machine accepting the language of words $w$ representing $1$ in $\PP$ in non-deterministic time equivalent to the Dehn function of $\PP$ and in non-deterministic space equivalent to the filling length function of $\PP$.   
\end{thm}

\ni \emph{Sketch proof.}   The construction of a non-deterministic Turing machine we give here is based on \cite{SBR}, however we shirk the formal descriptions of machine states and transitions (which are set out carefully in \cite{SBR}).  Let $\langle X \mid R \rangle$ be a finite presentation of $\PP$ where for all  $r \in R$ we find $R$ also contains all cyclic permutations of $r$ and $r^{-1}$. Our machine has two tapes and two heads; its alphabet is $X^{\pm 1}$ and it will start with an input word $w$ written on the first tape and with the second tape blank.   At any time during the running of the machine there will be a word (possibly empty, but containing no blank spaces) on each tape and the remainder of the tapes will be blank; the head on each tape will always be located between two squares, that to its right being blank and that to its left containing the final letter of the word on the tape (assuming that word is not empty).  We allow a head to see only  the one square to its right and the $\rho:=\max_{r \in R} \ell(r)$ squares to its left (a convenience to make the \emph{apply-a-relator} transition reversible). The transitions are:

\ms \ni \emph{Shift.}  Erase the final letter $x$ from the word on the end of one tape and insert $x^{-1}$ on the end of the word on the other tape. 

\smallskip \ni \emph{Apply-a-relator.}  Append/remove a word from $R$ to/from the word on the first tape as a suffix. 

\smallskip \ni \emph{Expand/Reduce.}  Append/remove a letter in $X^{\pm1}$ to/from both words. 

\smallskip \ni \emph{Accept.}  Accept $w$ if the squares seen by the two heads are empty. 

\ms \ni Notice that acceptance occurs precisely when both tapes are empty.  Consider the words one gets by concatenating the word on the first tape with the inverse of the word on the second at any given time -- as the machine progresses this gives a null-sequence;  indeed, (modulo the fact that it may take multiple transitions to realise each move) any null-sequence can be realised as a run of the machine.  So, as the \textsc{space} of a run of the machine is the maximum number of non-blank squares in the course of the run, the space complexity of this machine is equivalent to the filling length function of $\PP$.

Checking that the time complexity of our machine is equivalent to the Dehn function of $\PP$ requires more care.  Suppose $w$ is a null-homotopic word.  It is clear that any run of the machine on input $w$ takes time at least $\Area(w)$ because it must involve at least $\Area(w)$ \emph{apply-a-relator} transitions.  The reverse bound takes more care -- we wish to show that there is a run of the machine on input $w$ in which not only the \emph{apply-a-relator} transitions, but also the number \emph{shift} and \emph{expand/reduce} transitions, can bounded above by a constant times $\Area(w) + \ell(w)$. 
Sapir \cite[Section~3.2]{SapirSurvey} suggests using \cite[Lemma~1]{OS} to simplify the proof of \cite{SBR}.  We outline a similar approach based on the proof of Theorem~\ref{FLbyDGL}.  

Let $\Delta$ be a minimal area van~Kampen diagram for $w$ and let $T$ be a spanning tree in the 1-skeleton of $\Delta$.  Note that the total number of edges in $T$ is at most $\rho \, \Area(w) + \ell(w)$. 
A shelling of $\Delta$ is described in the proof of Theorem~\ref{FLbyDGL} in which one collapses 1-cells and 2-cells as one encounters them when following the boundary circuit of a small neighbourhood of $T$.  Realise the null-sequence of this shelling (that is the sequence of boundary words of its diagrams), as the run of our Turing machine.  The total number of moves is controlled as required because the shift moves correspond to the journey around $T$.  
\qed


\bs

\ni The message of these results is that $\Area: \N \to \N$ and $\FL: \N \to \N$ are complexity measures of the word problem for $\PP$ --  they are the \textsc{Non-Deterministic Time} (up to $\simeq$-equivalence)  and the \textsc{Non-Deterministic Space} of the na\"ive approach to solving the word problem in the Dehn proof system, by haphazardly applying relators, freely reducing and freely expanding in the hope of achieving an empty tape.  
This theme is taken up by Birget \cite{Birget} who describes it as a ``connection between \emph{static} fillings (like length and area) and \emph{dynamic} fillings (e.g., space and time complexity of calculations).'' 
It opens up the tantalising possibility of, given a computational problem, translating it, suitably effectively, to the word problem in the Dehn proof system of some finite presentation.  Then geometric considerations can be brought to bear on questions concerning algorithmic complexity.  

Striking results in this direction have been obtained by Sapir and his collaborators using \emph{$\mathcal{S}$-machines}.  See \cite{SapirSurvey} for a recent survey.  Their results include the following.  

\begin{thm} [\textup{Sapir-Birget-Rips \cite{SBR}}]  
 \label{realising a TM}
If $f(n)$ is the time function of a non-deterministic Turing machine such that $f^4(n)$ is super-additive \textup{(}that is, $f^4(m+n) \geq f^4(m) + f^4(n)$ for all $m,n$\textup{)}  then there is a finite presentation with $\Area(n) \simeq f^4(n)$ and $\IDiam(n) \simeq f^3(n)$.  
\end{thm}

\ni Building on \cite{Birget2, Ol3, SBR}, a far-reaching \emph{controlled} embedding result (Theorem~1.1)  is proved in \cite{BORS}, of which the following theorem is a remarkable instance.  Theorem~\ref{Turing machine for P} gives the \emph{if} implication; the hard work is in the converse.  
[Recall that a problem is NP when it is decidable in non-deterministic polynomial time.] 

\begin{thm}[\textup{Birget-Ol'shanskii-Rips-Sapir \cite{BORS}}] \label{NP}
The word problem of a finitely generated group $G$ is an NP-problem if and only if $G \hookrightarrow \hat{G}$ for some finitely presentable group $\hat{G}$ such that the Dehn function of  $\hat{G}$ is bounded above by a polynomial.  Indeed, $G \hookrightarrow \hat{G}$ can be taken to be a quasi-isometric embedding. 
\end{thm}

\ni Birget conjectures the following analogue of Theorem~\ref{NP} for filling length, relating it to \textsc{non-deterministic symmetric space} complexity.   The \emph{if} part of the conjecture is covered by Theorem~\ref{Turing machine for P}.  

\begin{Conjecture}[\textup{Birget \cite{Birget}}]
For a finitely generated group $G$ and function $f: \N \to \N$, there is a non-deterministic symmetric Turing machine accepting the language of words $w$ representing $1$ in $G$ within space  $\preceq f(\ell(w))$ if and only if $G$ embeds in a finitely presentable group with filling length function $\preceq f$.  
\end{Conjecture}

\smallskip

\ni Another approach is pursued by Carbone in \cite{CarboneBook, Carbone} where she relates van~Kampen diagrams to \emph{resolution proofs} in logic, in such a way that the Dehn function corresponds to proof length.    

\ms

\ni Our discussions so far in this section get us some way towards a proof of the following theorem.  We encourage the reader to complete the proof as an exercise (or to refer to \cite{Gersten} or \cite{ShortNotes}).  

\begin{thm}[\textup{Gersten \cite{Gersten}}] \label{word problem thm}
For a finite presentation the following are equivalent.  
\begin{itemize}
\item The word problem for $\PP$ is solvable.  
\item  $\Area: \N \to \N$ is bounded above by a recursive function.
\item  $\Area: \N \to \N$ \emph{is} a recursive function.
\end{itemize}

\end{thm}

\ms

\ni So, strikingly, finite presentations of groups with unsolvable word problem (which do exist \cite{Boone, Novikov}) have Dehn functions out-growing  all recursive functions.   Moreover, with the results we will explain in Chapter~\ref{interrelationships}, it can be shown that for a finite presentation,  \emph{one} of the filling functions 
$$\Area, \FL, \IDiam, \GL, \DGL: \N \to \N$$ is a recursive function if and only if \emph{all} are recursive functions. 

  However, in general, Dehn function is a poor measure of the time-complexity of the word problem for a group $G$.  Cohen \cite{CMO} and Madlener \& Otto \cite{MO} showed that amongst all finite presentations of groups, there is no upper bound on the size of the gap between the Dehn function and the time-complexity of the word problem (in the sense of the Grzegorczyk hierarchy).  

Indeed \cite{SapirSurvey}, the constructions of \cite{BORS, SBR} can be used to produce finite presentations with word problems solvable in quadratic time but arbitrarily large (recursive) Dehn functions.  The seed idea is that an embedding of a finitely generated group $G$ into another group $\hat{G}$ that has an efficient  algorithm to solve its word problem leads to an efficient algorithm for the word problem in $G$.  
For example, $\langle \, a,b \, \mid \, a^b =a^2 \, \rangle$ has an exponential Dehn function (Proposition~\ref{BS Dehn fn}), but a polynomial time word problem:

\begin{exer} 
Find a deterministic polynomial time algorithm for $\langle \, a,b \, \mid \, a^b =a^2 \, \rangle$.  \emph{Hint}:  $\langle \, a,b \, \mid \, a^b =a^2 \, \rangle$ is subgroup of $\GL_2(\mathbb{Q})$ via $a = \left(\begin{array}{cc}1 & 1 \\0 & 1\end{array}\right)$ and $b= \left(\begin{array}{cc}1/2 & 0 \\0 & 1\end{array}\right)$.
\end{exer}

\ni In fact, this approach leads to an $n (\log n)^2 (\log \log n)$ time solution -- see \cite{BORS}.

\ms 
\ni The group $$\langle \  a,b  \ \mid \  a^{a^b} a^{-2} \  \rangle,$$ introduced by Baumslag in \cite{Baumslag2}, is an even more striking example. Its Dehn function was identified by Platonov \cite{Platonov} as $$\Area(n) \  \simeq  \  \stackrel{\log_2 n}{\overbrace{f(\, f \ldots (f}(1))} \ldots),$$ where $f(n):=2^n$.  
Earlier Gersten \cite{Gersten6, Gersten} had shown it's Dehn function to grow faster than every iterated exponential and Bernasconi \cite{Bernasconi} had found a weaker upper bound.
However I.Kapovich \& Schupp  \cite{SchuppPersonal} and Miasnikov,  Ushakov \& Wong \cite{MUW} claim its word problem is solvable in polynomial time. In the light of this and Magnus' result that every 1-relator group has solvable word problem \cite{Lyndon}, Schupp sets the challenge \cite{SchuppPersonal}:

\ms

\begin{Open Problem}
Find a one-relator group for which there is, provably, no algorithm to solve the word problem within linear time. 
\end{Open Problem}

\ms

\ni In this context we mention Bernasconi's result \cite{Bernasconi} that all Dehn functions of 1-relator groups are at most (that is, $\preceq$) Ackermann's function (which is recursive but not primitive recursive), and a question of Gersten:  

\ms

\begin{Open Problem}
Is the Dehn function of every one-relator group bounded above by (that is, $\preceq$) the Dehn function of $\langle \  a,b  \ \mid \  a^{a^b} a^{-2} \  \rangle$?
\end{Open Problem}


\section{Filling functions for Riemannian manifolds}

For a closed, connected Riemannian manifold $M$, define $\Area_M:[0,\infty) \to[0,\infty)$ by 
\begin{equation*}
\Area_{M}(l) \  := \  \sup_c\inf_{D}\{ \ \Area(D) \  \mid \ D:\mathbb D^2\to \widetilde{M},\, D\restricted{\partial \D^2}=c,\, \ell(c)\le l \ \}, 
\end{equation*}
where $\Area(D)$ is 2-dimensional Hausdorff measure and the infimum is over all Lipschitz maps $D$.  Some remarks, following Bridson~\cite{Bridson6}:  by Morrey's solution to Plateau's problem \cite{Morrey}, for a fixed $c$, the infimum is realised; and, due to the regularity of the situation, other standard notions of area would be equivalent here.   

One can construct a $\Gamma$-equivariant map $\Phi$ of $\Cay^2(\PP)$ into the universal cover $\widetilde{M}$ of $M$ by first mapping $1$ to some base point $v$, then extending to $\Gamma$ (which we identify with the 0-skeleton of $\Cay^2(\PP)$) by mapping $\gamma \in \Gamma$ to the translate of $v$ under the corresponding deck transformation, then extending to the 1-skeleton by joining the images of adjacent vertices by an equivariantly chosen geodesic for each edge, and then to the 2-skeleton by equivariantly chosen finite area fillings for each 2-cell.  Loops in $M$ can be approximated by \emph{group-like paths} -- the images under $\Phi$ of edge-paths in $\Cay^2(\PP)$ -- and maps of discs $D:\mathbb D^2 \to \widetilde{M}$ can be related to images under  $\Phi$ of  van~Kampen diagrams   into $\widetilde{M}$ (or, more strictly, compositions with $\Phi$).  This is the intuition for the following result, however a full proof requires considerable technical care.   

\begin{thm}[The Filling Theorem \textup{\cite{Bridson6, Burillo-Taback, Gromov3}}] \label{Filling Thm}
For a closed, connected Riemannian manifold $M$, the Dehn function of any finite presentation of $\Gamma := \pi_1 M$ is $\simeq$-equivalent to $\Area_{M}:  [0,\infty) \to[0,\infty)$.
\end{thm}

Similar results hold for intrinsic diameter, extrinsic diameter and filling length.
The \emph{intrinsic diameter} of a continuous map $D: \mathbb{D}^2 \to \widetilde{M}$ is $$\IDiam_M(D) \ := \ \sup \set{\ \rho(a,b) \mid a,b \in \mathbb{D}^2 \ }$$ where $\rho$ is the  pull back of the Riemannian metric:
 $$\rho(a,b) \ = \  \inf \set{ \ \ell(D \circ p) \mid p \textup{ a path in } \mathbb{D}^2 \textup{ from } a \textup{ to } b \  }.$$  [Note that there may be distinct points $a,b \in \mathbb{D}^2$ for which $\rho(a,b)=0$ -- that is, $\rho$ may be only a \emph{pseudo}-metric $\mathbb{D}^2$.]
The extrinsic diameter $\EDiam_M(D)$ is the diameter of $D(\mathbb{D}^2)$, as measured with the distance function on $\widetilde{M}$.  These two notions of diameter give functionals on the space of rectifiable loops in 
$\widetilde{M}$:  the \emph{intrinsic} and \emph{extrinsic} diameter functionals $\IDiam_M, \EDiam_M:[0,\infty) \to[0,\infty)$ -- at $l$ these functionals take the values
\begin{eqnarray*}
& & \sup_c\inf_{D \in \mathcal{D}}\{ \ \IDiam_M(D) \  \mid \ D:\mathbb D^2\to \widetilde{M},\, D\restricted{\partial \D^2}=c,\, \ell(c)\le l \ \} \ \  \text{ and} \\ 
& & \sup_c\inf_{D \in \mathcal{D}}\{ \ \EDiam_M(D) \ \mid \  D:\mathbb D^2\to \widetilde{M},\, D\restricted{\partial \D^2}=c,\, \ell(c)\le l \ \},
\end{eqnarray*}
respectively. 
 
The filling length $\FL(c)$ of a rectifiable loop $c: [0,1] \to \widetilde{M}$ is defined in \cite{Gromov} to be the infimal length $L$ such that there is a null-homotopy $H: [0,1]^2 \to \widetilde{M}$ where $H(s,0) = c(s)$ and $H(0,t) = H(1,t) = H(s,1)$ for all $s,t \in [0,1]$, such that for all $t  \in [0,1]$ we find $s \mapsto H(s,t)$ is a loop of length at most $L$.  Then the \emph{filling length functional} $\FL_M$ is defined to be the supremum of $\FL(c)$ over all rectifiable loops $c: [0,1] \to \widetilde{M}$ of length at most $c$.       
 
The analogue of the Filling Theorem is:
 
\begin{thm} \textup{\cite{BR1, BR2}} \label{Riem-Comb}
If $\PP$ is a finite presentation of the fundamental group $\Gamma$ of a closed, connected Riemannian manifold $M$ then $$\IDiam_{\PP} \simeq \IDiam_{M}, \ \  \EDiam_{\PP} \simeq \EDiam_{M}, \ \  \textup{ and } \ \  \FL_{\PP} \simeq \FL_{M} .$$
\end{thm}

\section{Quasi-isometry invariance}\thispagestyle{empty} \label{qi section}

The Dehn function is a \emph{presentation invariant} and, more generally, a \emph{quasi-isometry} invariant, in the sense of the following theorem. 
For a proof and background see \cite{Alonso} or \cite[Theorem~4.7]{Short}. 

\begin{thm} \label{q.i.}
If $G, H$ are groups with finite generating sets $X,Y$, respectively, if $(G, d_X)$ and $(H, d_Y)$ are quasi-isometric,  and if  $G$ is finitely presentable, then $H$ is finitely presentable.  Moreover, if $G$ and $H$ have finite presentations  $\PP$ and $\QQ$, respectively, then the associated Dehn functions $\Area_{\PP}, \Area_{\QQ} : \N \to \N$ satisfy $\Area_{\PP} \simeq \Area_{\QQ}$.  In particular, the Dehn functions of two finite presentations of the same group are $\simeq$ equivalent. 
\end{thm}

\ni Thus it makes sense to say that a finitely presentable group $G$ has a linear, quadratic, $n^{\alpha}$ (with $\alpha \geq 1$), exponential etc.\ Dehn function, meaning $$\Area(n) \ \simeq \  n, n^2, n^{\alpha}, \exp(n), \ldots$$ with respect to some, and hence any, finite presentation.  Note also, that for $\alpha, \beta \geq 1$ and  $f,g: [0, \infty) \to [0, \infty)$ defined by $f(n) = n^{\alpha}$ and $g(n) = n^{\beta}$  we have $f \simeq g$ if and only if $\alpha = \beta$.  

\ms

\ni Results analogous to Theorem~\ref{q.i.} also apply to $\EDiam$, $\IDiam$ and $\FL$ -- see \cite{BR1, BR2}.  
In each case the proof involves monitoring diagram measurements is the course of the standard proof that finite presentability is a quasi-isometry invariant -- the first such quantitative version was \cite{Alonso}.  Proofs of the independence, up to $\simeq$ equivalence, of $\IDiam : \N \to \N$ and $\FL : \N \to \N$ on the finite presentation predating  \cite{BR1, BR2} use Tietze transformations and are in \cite{Gersten-Short} and \cite{GR1}.  In \cite{GR3} it is shown that the gallery length functions of two \emph{fat} finite presentations of the same group are $\simeq$ equivalent.  One expects that  \emph{fat} finite presentations of quasi-isometric groups would have $\simeq$ equivalent gallery length functions, and $\DGL$ functions, likewise.

\begin{center}
\begin{table}
\caption{A summary of the multiple interpretations of filling functions}

\vspace*{5mm}
\begin{tabular}{| l || cccc | }
\hline

\multirow{5}{*}{}  &     &   &  & \multirow{5}{*}{$G = \pi_1(M)$} \\
& van~Kampen  &$P \  = \  \prod_{i=1}^N {u_i}^{-1} {r_i}^{\epsilon_i} u_i $\hspace*{-.7cm}  & Dehn proof & \\
& diagrams & equalling & system & \\
& $\pi: \Delta \to \Cay^2(\PP)$ &  $w$ in $F(X)$ & & \\
& for $w$ & & & \\

\hline \hline

\multirow{5}{*}{$\Area(w)$} 	&  \multirow{5}{*}{$\displaystyle{\min_{\Delta}\set{\# \text{ of 2-cells in } \Delta}}$}  & \multirow{5}{*}{$\displaystyle{\min_P N}$} & \multirow{5}{*}{} & \multirow{5}{*}{} \\
& & & \textsc{Non-deterministic} & Riemannian  \\
&  & & \textsc{time} & $\Area$ \\
& & & & \\
& & & & \\

\hline

\multirow{5}{*}{$\FL(w)$} 	&    & \multirow{5}{*}{?} & \multirow{5}{*}{\textsc{Space}} & \multirow{5}{*}{} \\
& $\min\{ \, \FL( \SS) \ |  \textup{ shellings }   $  & & & Riemannian  \\
& $ \ \  \  \ \ \ \ \SS = (\Delta_i) \textup{ of } \Delta \, \},$ & & & $\FL$ \\
& where $\displaystyle{\FL(\SS) = \max_i \ell(\partial \Delta_i)}$\hspace*{-.7cm}  & & & \\
& & & & \\

\hline

\multirow{5}{*}{$\IDiam(w)$} 	&  \multirow{5}{*}{$\displaystyle{\min_{\Delta} \Diam \ \Delta^{(1)}}$}  & \multirow{5}{*}{$\displaystyle{\min_P \max_i \ell(u_i)}$} & \multirow{5}{*}{?${}^{\ast}$} & \multirow{5}{*}{} \\
& & & & Riemannian  \\
&  & & & $\IDiam$ \\
& & & & \\
& & & & \\

\hline

\multirow{5}{*}{$\EDiam(w)$} 	&  \multirow{5}{*}{$\displaystyle{\min_{\Delta} \Diam \ \pi(\Delta^{(1)})}$}  & & \multirow{5}{*}{?} & \multirow{5}{*}{} \\
&  & $\displaystyle{ \min_P\max \, \{ \, d_X(\star,u) \, | }$\hspace*{-.7cm}   & & Riemannian  \\
&  & $ \ \ \ u \textup{ a prefix of } w \}$ & & $\EDiam$ \\
& & & & \\
& & & & \\

\hline
\multirow{5}{*}{$\GL(w)$} 	&  \multirow{5}{*}{$\displaystyle{\min_{\Delta}\Diam \ \left(\Delta^{\ast}\right)^{(1)}}$}  & \multirow{5}{*}{?} & \multirow{5}{*}{?} & \multirow{5}{*}{?} \\
& & & & \\
&  & & & \\
& & & & \\
& & & & \\

\hline
\multirow{5}{*}{$\DGL(w)$} 	&    & \multirow{5}{*}{?} & \multirow{5}{*}{?$^{\ast}{}^{\ast}$} & \multirow{5}{*}{?$^{\ast}{}^{\ast}$} \\
& $\displaystyle{\min_{\Delta} \, \{ \, \Diam T + \Diam T^{\ast} \,  |  }$\hspace*{-.7cm}  & & & \\
&  $ \ \ \ T \textup{ a spanning tree in } \Delta^{(1)} \, \} $\hspace*{-.7cm}  & & & \\
& & & & \\
& & & & \\

\hline
  \end{tabular} \\
\rule{0mm}{8mm}\hspace*{5mm} \parbox{12.5cm}{It remains a challenge to supply appropriate interpretations (if they exist) in place of each question mark (modulo ${}^{\ast}$ and ${}^{\ast}{}^{\ast}$).} 

\smallskip
 \ni \hspace*{5mm} \parbox{12.5cm}{${}^{\ast}$ Birget \cite[Proposition~5.1]{Birget} relates $\IDiam(w)$ to an exponential of the deterministic time of an approach to solving the word problem by constructing a non-deterministic finite automata based on part of the Cayley graph.} \\ 

\smallskip
\ni \hspace*{5mm}  \parbox{12.5cm}{${}^{\ast}{}^{\ast}$ By Theorem~\ref{FL and DGL} we have $\DGL \simeq FL$ for \emph{fat} presentations.}

\end{table}
\end{center}


\chapter{Relationships between filling functions}\thispagestyle{empty}\label{interrelationships}

This chapter concerns relationships between filling functions that apply irrespective of the group being presented.  A more comprehensive account of known relationships is in \cite{GR4}, along with a description of how there are many additional coincidences between filling functions if one only uses van~Kampen diagrams whose vertices have valence at most some constant.  

We fix a finite presentation $\PP = \langle X \mid R \rangle$ of  a group $G$ for the whole chapter.  Some relationships are easy to find.  Suppose $\pi : \Delta \to \Cay^2(\PP)$ is a van~Kampen diagram with base vertex $\star$.   Then $$\IDiam(\Delta) \ \leq \ \FL(\Delta)$$ because the images in $\Delta$ of the boundaries $\partial \Delta_i$ of the diagrams in a shelling form a family of contracting loops that at some stage pass through any given vertex $v$  and so provide paths of length at most $\FL(\Delta)/2$ to $\star$.  
Defining $K := \max \set{ \ell(r) \mid r \in R}$, the total number of edges in $\Delta$ is at most  $K \Area(\Delta) + \ell(\partial \Delta)$ and so 
$$\FL(\Delta)  \ \leq \   K \Area(\Delta) + \ell(\partial \Delta).$$ 
It follows that for all $n$, 
$$  \IDiam(n) \ \leq \ \FL(n)  \ \leq \  K \Area(n) + n.$$ 
Similarly, one can show that 
\begin{equation*}
\GL(n) \  \leq \  \DGL(n)  \ \leq \ K \Area(n) +  n \ \ \text{ and } \ \ 
\IDiam(n) \ \leq \ K \GL(n).
\end{equation*}

Next we give the ``\emph{space-time}''  bound \cite[Corollary 2]{GR1}, \cite[5.C]{Gromov} -- an analogue of a result in complexity theory -- 

\begin{prop}\label{space-time bound}
Define $K:= 2 \abs{X} +1$.  Then for all $n$, $$\Area(n) \leq K^{\mbox{$\FL(n)$}}.$$ 
\end{prop}

\begin{proof} 
This result is most transparent in the Dehn proof system of Section~\ref{Dehn proof system}.  The number of words of length at most $\FL(n)$, and hence the number of different words occurring in a null-sequence for a null-homotopic word $w$ of length at most $n$,  is at most $K^{\mbox{$\FL(n)$}}$.  Every such $w$ has a null-sequence in which no word occurs twice and so there are at most $K^{\mbox{$\FL(n)$}}$ application-of-a-relator moves in this null sequence.  
\end{proof}

\smallskip
\ni Proposition~\ref{space-time bound} and Theorem~\ref{Combable groups} combine to give the following result which may be surprising in that it makes no reference to the length function of the combing.

\begin{cor} \textup{\cite{Gersten5, GR1}}
The Dehn function of an \textup{(}a\textup{)}synchronously combable group grows at most exponentially fast.   
\end{cor}

\section{The Double Exponential Theorem} \label{det}

Gallery length can be used to prove a theorem of D.E.Cohen \cite{Cohen}, known as \emph{the Double Exponential Theorem}.  Cohen's proof involves an analysis of the Nielsen reduction process.  A proof, using \emph{Stallings folds}, was given by Gersten \cite{Gersten2} and was generalised by Papasoglu \cite{Papasoglu2} to the more general setting of filling edge-loops in a simply connected complex of uniformly bounded local geometry. Birget found a proof based on context free languages \cite{Birget}.
The proof below is from \cite{GR3}.  It combines two propositions each of which establish (at most) exponential leaps, the first from $\IDiam$ up to $\GL$, and the second from $\GL$ up to $\Area$. 

A key concept used is that of a \emph{geodesic spanning tree} based at a vertex $\star$ in a graph $\Gamma$ -- that is, a spanning tree such that for every vertex $v$ in $\Gamma$, the combinatorial distances from $v$ to $\star$ in $\Gamma$ and in $T$ agree.  

\begin{exer}
Prove that for every finite connected graph $\Gamma$ and every vertex $\star$ in $\Gamma$ there is a geodesic spanning tree in $\Gamma$ based at $\star$.  
\end{exer} 

\begin{thm} \label{Double exp thm} There exists $C>0$, depending only on $\PP= \langle X \mid R \rangle$, such that 
$$\Area (n) \ \leq \  n \, C^{\mbox{$C^{\mbox{$\IDiam(n)$}}$}}$$ for all
$n \in \N$.
\end{thm}

\ni Define the \emph{based} intrinsic diameter of a diagram $\Delta$ with base vertex $\star$ by  
$$ \IDiam_{\star}(\Delta) \ := \  \max \set{ \ \rho ( \star, v) \ \mid \ v \in \Delta^{(0)} \ },$$ where $\rho$ denotes the combinatorial metric on $\Delta^{(1)}$.  Note that, by definition, $\IDiam_{\star}(\Delta) \leq \IDiam(\Delta)$.  

\begin{prop}\label{GL by Diam}
Suppose $w$ is a null-homotopic word and $\Delta$ is a van~Kampen diagram for $w$ such that $\IDiam_{\star}(\Delta)$ is minimal.  Moreover, assume $\Delta$ is of minimal area amongst all such diagrams.  For  $ A:= 2\abs{X}+1$,
\begin{eqnarray*}
\GL(\Delta) & \leq & 2 A^{\mbox{$1+ 2 \, \IDiam(\Delta)$}}.
\end{eqnarray*}
So $\GL(n) \leq 2 A^{\mbox{$1+ 2 \, \IDiam(n)$}}$ for
all $n$.
\end{prop}

\begin{proof}  
Let $T$ be a geodesic spanning tree in $\Delta^{(1)}$ based at
$\star$.  (See Figure~\ref{double_exp_figs}.)  To every edge $e$ in $\Delta^{(1)} \ssm T$ associate an
anticlockwise edge-circuit $\gamma_e$ in $\Delta^{(1)}$, based at $\star$, by
connecting the vertices of $e$ to $\Delta^{(1)}$ by geodesic paths in $T$.
Let $w_e$ be the word one reads along $\gamma_e$. Then $\ell(w_e)
\leq 1+ 2 \, \IDiam(\Delta^{(1)})$.

\begin{figure}[ht]
\psfrag{f}{{$e_{\infty}\*$
}}
\psfrag{e1}{{$e_1$}}
\psfrag{e2}{{$e_2$}}
\psfrag{eta}{{$\eta$}}
\psfrag{T}{{$T$}}
\psfrag{s}{{$\star$}}
\centerline{\epsfig{file=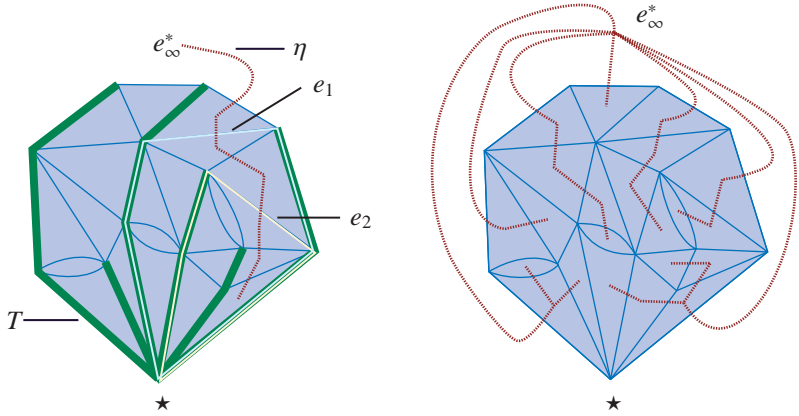}}
\caption{Illustrations of the proofs of Propositions~\ref{GL by Diam} and \ref{Area by GL}} \label{double_exp_figs}
\end{figure}

Let $T\*$ be the dual tree to $T$ (as defined in
Section~\ref{via vK diagrams}) and $e_{\infty}\*$ be the vertex of
$T\*$ dual to the face $e_{\infty}$.
Suppose that $e_1$ and $e_2$ are two distinct edges of $\Delta^{(1)} \ssm T$
dual to edges $e_1\*$ and $e_2\*$ of $T\*$ on some geodesic $\eta$ in
$T\*$ from $e_{\infty}\*$ to a leaf and that $e_2\*$ is further from $e_{\infty}\*$ than
$e_1\*$ along $\eta$.  Then $\gamma_{e_1}$ and $\gamma_{e_2}$ bound
subcomplexes  $C(\gamma_{e_1})$ and $C(\gamma_{e_2})$ of $\Delta$ in
such a way that $C(\gamma_{e_2})$ is a subcomplex of
$C(\gamma_{e_1})$.  Suppose that  $w_{e_1}$ and $w_{e_2}$
are the same words.  Then we could cut $C(\gamma_{e_1})$ out of $\Delta$
and glue $C(\gamma_{e_2})$ in its place, producing a new
van~Kampen diagram for $w$ with strictly smaller area and no
increase in $\IDiam_{\star}$ --  a contradiction.

Therefore $\ell(\eta) \leq A^{\mbox{$1+ 2\,\IDiam(\Delta)$}}$, as the right-hand side is an upper bound on the number of distinct words of length at most $1+ 2 \, \IDiam(\Delta)$, and the result follows.
\end{proof}

\begin{prop}\label{Area by GL}  If $\Delta$ is a van Kampen diagram
for a null-homotopic word $w$ then
\begin{eqnarray*}
\Area(\Delta) \leq \ell(w) (B+1)^{\mbox{$\GL(\Delta)$}}, 
\end{eqnarray*} where
$B:= \max \set{\ell(r) \mid r \in R}$. It follows that $\Area(n) \leq n
(B+1)^{\mbox{$\GL(n)$}}$ for all $n$.
\end{prop}

\begin{proof}
Sum over $1\leq k \leq \GL(\Delta)$, the geometric series whose $k$-th term $n(B+1)^{k-1}$ dominates the number of vertices at distance $k$ from $e^{\ast}_{\infty}$  in a geodesic spanning tree in the 1-skeleton of $\Delta^{\ast}$, based at $e^{\ast}_{\infty}$, as illustrated in Figure~\ref{double_exp_figs}.  
\end{proof}

\smallskip

\begin{rem} \label{GL in middle}
In fact, the two propositions above establish more than claimed in Theorem~\ref{Double exp thm}.  They show that there exists $C>0$, depending only on $\PP$, such the bounds  $$ \Area(n) \ \leq \  n \,  C^{\mbox{$C^{\mbox{$ \Diam(n)$}}$}}  \ \ \text{ and } \ \  \GL(n) \  \leq \  C^{\mbox{$ \Diam(n)$}}$$ are simultaneously realisable on van~Kampen diagrams that are of minimal $\IDiam$ amongst all diagrams for their boundary words.  This together with an inequality of \cite{GR1} show (as observed by Birget~\cite{Birget}) that $\IDiam(n)$ is always no more than a single exponential of $\FL(n)$.   
\end{rem}

\smallskip

\begin{prop} \label{BS Dehn fn}
The presentation $\langle \, a,b \, \mid  \, a^b =a^2 \, \rangle$ has filling length and Dehn functions satisfying $$\FL(n) \simeq n  \ \ \text{ and } \ \  \Area(n) \simeq \exp(n).$$  
\end{prop}

\begin{proof} 
That $\FL \preceq n$ follows from asynchronous combability and Theorem~\ref{Combable groups}, and so $\Area(n) \preceq \exp(n)$ by Proposition~\ref{space-time bound};    $\FL \succeq n$  by definition of $\succeq$, and $\Area(n) \succeq \exp(n)$ by the following lemma (see e.g.\ \cite[Section~2.3]{GR2},  \cite[Section~4C${}_2$]{Gromov}) applied to fillings for words $\left[a , a^{b^n}\right]$ as illustrated in Figure~\ref{BS fig} in the case $n=6$.  
\end{proof}

A presentation  $\PP$ is \emph{aspherical} when $\Cay^2(\PP)$ is contractible.

\begin{lem}[Gersten's Asphericity Lemma \cite{Bridson6, Gersten4}] \label{asphericity}
If $\PP$ is a finite aspherical presentation and $\pi: \Delta \to \Cay^2(\PP)$ is a van~Kampen diagram such that $\pi$ is injective on the complement of $\Delta^{(1)}$ then $\Delta$ is of minimal area amongst all van~Kampen diagrams with the same boundary circuit.    
\end{lem}  

\begin{exer}
Give another proof that  $\Area(n) \succeq \exp(n)$ for  $\langle \, a,b \, \mid  \, a^b =a^2 \, \rangle$ by studying the geometry of   \emph{$b$-corridors} and the words in $\set{a^{\pm 1}}^{\ast}$ along their sides in van~Kampen diagrams for $[a,a^{b^n}]$.  (A \emph{$b$-corridor} is a concatenation of 2-cells, each joined to the next along an edge labelled by $b$, and running between two oppositely oriented edges on $\partial \Delta$ both labelled $b$.  For example, every horizontal strip of 2-cells in the van~Kampen diagram in Figure~\ref{BS fig} is a $b$-corridor.)   
\end{exer}

\bs

\ni Remarkably, no presentation is known for which the leap from $\IDiam(n)$ to $\Area(n)$ is more than the single exponential of  $\langle \, a,b \, \mid  \, a^b =a^2 \, \rangle$.  So many (e.g.\ in print \cite{Cohen, Gersten, Gromov}, but it is often attributed to Stallings and indeed said that he conjectured a positive answer) have raised the natural question: 

\ms

\begin{Open Question} \label{sep} 
Given a finite presentation, does there always exist $C>0$ such that $\Area(n) \leq {C}^{\mbox{$\IDiam(n) +n$}}$ for all $n$?
\end{Open Question}

\ni Indeed, one can ask more \cite{GR3}: is a bound $\Area(n) \leq {C}^{\mbox{$\IDiam(n) +n$}}$ always realisable on van~Kampen diagrams of minimal $\IDiam$ for their boundary words.

\ms 

\ni Gromov \cite[page 100]{Gromov} asks:

\begin{Open Question}  \label{flbyd} 
Given a finite presentation, does there always exist $C>0$ such that $\FL(n) \leq C (\IDiam(n)+n)$ for all $n$?
\end{Open Question}

\ni This question has also been attributed to Casson~\cite[page~101]{Gromov}.  Gromov notes that an affirmative answer to \ref{flbyd} would imply an affirmative answer to \ref{sep} by Proposition~\ref{space-time bound}.  At the level of combinatorial or metric discs the analogue of \ref{flbyd} is false -- Frankel \& Katz \cite{Frankel-Katz} produced a family of metric discs $D_n$ with diameter and perimeter $1$ but such that however one null-homotopes $\partial D_n$ across $D_n$, one will (at some time) encounter a loop of length at least $n$; we will describe combinatorial analogues of these discs in Section~\ref{duality etc}.   

In \cite{GR3} Gersten and the author ask \ref{flbyd} but with $\FL$ is replaced by $\GL$: 
\begin{Open Question}  \label{glbyd} 
Given a finite presentation, does there always exist $C>0$ such that $\GL(n) \leq C (\IDiam(n)+n)$ for all $n$?
\end{Open Question}

\ni As mentioned in Remark~\ref{GL in middle}, like filling length, gallery length sits  at most an exponential below the Dehn function and at most an exponential above the intrinsic diameter function, in general, and a positive answer to \ref{glbyd} would resolve \ref{sep} positively.  However, \ref{glbyd} may represent a different challenge to \ref{flbyd} because separating the $\GL$ and $\IDiam$ functions has to involve a family of van~Kampen diagrams for which there is no uniform upper bound on the valence of vertices.  The features of combinatorial discs $\Delta_n$ discussed in  Section~\ref{duality etc} may  be reproducible in van~Kampen diagrams in some group in such a way as to answer \ref{flbyd}  negatively, but those discs have uniformly bounded vertex valences and so will not be of use for  \ref{glbyd}. 

\smallskip

\ms

\ni We now turn to what can be said in the direction of an affirmative answer to \ref{sep}.
For a null-homotopic word $w$, define
$$\overline{\IDiam}(w) \ := \  \max \set{ \ \IDiam(\Delta) \mid \text{minimal area van~Kampen diagrams } \Delta \text{ for } w \ }.$$  
Then, as usual, define a function $\N \to \N$, the \emph{upper intrinsic diameter function}, by 
$$\overline{\IDiam}(n) \ := \  \max \set{ \ \overline{\IDiam}(w) \mid \text{words } w \text{ with } \ell(w) \leq n \text{ and } w=1 \text{ in } G \ }.$$ A reason \ref{sep} is a hard problem is that, in general, minimal area  and minimal diameter fillings may fail to be realisable simultaneously as illustrated by the example in the following exercise.  

\ms 

\begin{exer}
By finding embedded van Kampen diagrams $\Delta_n$ for the words $w_n:=\left[x^{t^n}, y^{s^n}\right]$ and applying Lemma~\ref{asphericity}, show that the Dehn and $\overline{\IDiam}$ functions of the aspherical presentation $$\langle \  x,y,s,t \ \mid \ [x,y],  \ x^tx^{-2},  \ y^sy^{-2} \ \rangle$$ of Bridson~\cite{BridsonNotes, GR2}  satisfy $$\Area(n) \  \succeq \ 2^n \ \ \text{ and } \ \  \overline{\IDiam}(n) \ \succeq  \ 2^n.$$  But show that   $\IDiam(w_n) \preceq n$ for all $n$.  (\emph{Hint}.  Insert many copies of suitable van~Kampen diagrams for $x^{t^n}x^{-2^n}$ into $\Delta_n$.  In fact, it is proved in \cite{GR2} that $\IDiam(n) \preceq n$.)       
\end{exer}

\bs

\ni Using $\overline{\IDiam}$ side-steps the problem of minimal area and diameter bounds not being realised on the same van~Kampen diagram and the following  single exponential bound on $\Area$ can be obtained.  

\begin{thm} \label{set} \textup{\cite{GR2}} Given a finite presentation $\PP$, there exists $K >0$ such that $$\Area(n) \ \leq \ nK^{\mbox{$ \overline{\IDiam}(n)$}}.$$  
\end{thm}

\begin{figure}[ht]
\psfrag{0}{$A_0$}
\psfrag{1}{$A_1$}
\psfrag{2}{$A_2$}
\centerline{\epsfig{file=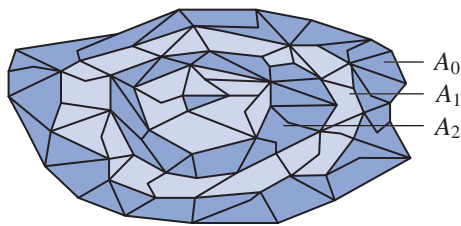}} 
\caption{Successive star-neighbourhoods of the boundary of a diagram}  
\label{star nbhds} 
\end{figure}

\ms \ni  \emph{Sketch proof.} 
Suppose $\Delta$ is a minimal area van~Kampen diagram for $w$.  Recall that topologically $\Delta$ is a singular 2-disc and so is a tree-like arrangement of arcs and topological 2-discs.  
Consider successive \emph{star-neighbourhoods} $\Star^{i}(\partial \Delta)$ of $\partial \Delta$ defined by setting  $\Star^{0}(\partial \Delta)$ to be $\partial \Delta$ and  $\Star^{i+1}(\partial \Delta)$ to be the union closed cells that share at least a vertex with $\Star^i(\partial \Delta)$.   Define $A_i:= \Star^{i+1}(\partial \Delta) \ssm \Star^{i}(\partial \Delta)$, a sequence of \emph{annuli}  that 
partition $\Delta$ into $(A_i)_{i=0}^m$ where $m \leq K_0 \IDiam(\Delta)$ for some constant $K_0$ depending only on $\PP$.  This is illustrated in Figure~\ref{star nbhds}.  (The $A_i$ will not be topological annuli, in general.)  

\ni Suppose that for every vertex $v$ on $\partial \Delta$ and every topological 2-disc component $\mathcal{C}$ of $\Delta$, there are at most three edges in the interior of $\mathcal{C}$  that are all incident with $v$, all oriented away from or all towards $v$, and all have the same label.  Then by considering  2-cells in $\mathcal{C}$ that are incident with $v$, we find there is constant $K_1$, depending only on $\PP$, such that $\Area(A_0)$ and the total length of the inner boundary circuit of $A_0$ are both at most $K_1 \ell(\partial \Delta)$.  

We then seek to argue similarly for $\Delta \ssm A_0$, to get an upper bound on the length of the inner boundary and the area of $A_1$.  And then continuing inductively at most $K_0 \IDiam(\Delta)$ times and summing an appropriate geometric series we will have our result.  

\begin{figure}[ht]
\psfrag{1}{$a$}
\psfrag{2}{$a$}
\psfrag{v}{$v$}
\psfrag{star}{$\star$}
\centerline{\epsfig{file=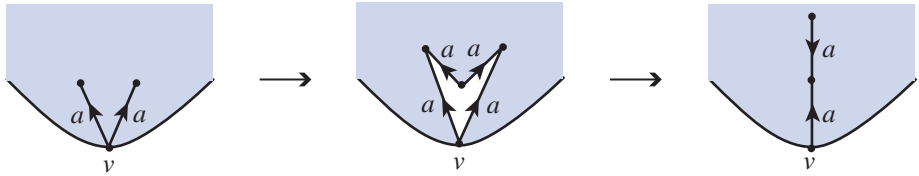}} 
\caption{A diamond move}
\label{diamond move}
\end{figure}

\ni The way we get this control on edges incident with vertices is to use
\emph{diamond moves}. These are performed on pairs of edges that are in the same 2-disc component of a van~Kampen diagram, that are oriented towards or away from a common end-vertex $v \in \partial \Delta$, that have the same label, and whose other vertices, $v_1$ and $v_2$, are either both in the interior of $\Delta$ or both on $\partial \Delta$.  The case where $v_1$ and $ v_2$ are in the interior is illustrated in Figure~\ref{diamond move}.  Performing a diamond moves when $v_1,  v_2 \in \partial \Delta$ creates a new 2-disc component.  

Do such diamond moves until none remain available -- the process does terminate because diamond moves either increase the number of 2-disc components (which can only happen finitely many times) or decreases the valence of some boundary vertex without altering the number of 2-disc components.  Then peal off $A_0$ and repeat the whole process on 
the boundary vertices of $\Delta \ssm A_0$.  

Continue likewise until the diagram has been exhausted.  This takes at most $\sim \overline{\IDiam}(w)$ iterations because, whilst the diamond moves may change $\IDiam$, they preserve $\Area$.   
\newqed

\bs

\begin{exer}
Use $\langle \, a,b \, \mid \, a^b =a^2 \, \rangle$ to show that the result of \ref{set} is sharp, in general.  
\end{exer}

\bs

\begin{Open Problem} \cite{GR2}
Give a general upper bound for $\overline{\IDiam} : \N \to \N$ in terms of $\IDiam : \N \to \N$.  Theorem~\ref{det} implies a double exponential bound holds; one might hope for a single exponential.
\end{Open Problem}

\section{Filling length and duality of spanning trees in planar graphs} \label{duality etc}

For a while, the following question and variants set out in \cite{GR4}, concerning planar graphs and dual pairs of trees in the sense of Section~\ref{via vK diagrams}, constituted an impasse in the study of relationships between filling functions.  [A \emph{multigraph} is a graph in which edges can form loops and pairs of vertices can be joined by the multiple edges.]

\ms

\begin{Question} \label{H}  Does there exist $K>0$ such that if $\Gamma$ is a finite connected planar
graph (or multigraph) then there is a spanning tree $T$ in $\Gamma$ such
that 
\begin{eqnarray*}
\Diam(T) & \leq & K \, \Diam(\Gamma) \ \  \text{ and} \\
\Diam(T^{\ast}) & \leq & K \, \Diam(\Gamma^{\ast})\,?
\end{eqnarray*}
\end{Question}

\ms

\ni It is easy to establish either one of the equalities with $K=2$: take a geodesic spanning tree based at some vertex (see Section~\ref{det}).  The trouble is that $T$ and $T^{\ast}$ determine each other and they fight -- altering one tree to reduce its diameter might increase the diameter of the other.  The following exercise sets out a family of examples in which taking one of $T$ and $T^{\ast}$ to be a geodesic tree based at some vertex, does not lead to $T$ and $T^{\ast}$ having the properties of Question~\ref{H}.  

\ms 

\begin{figure}[ht]
\psfrag{T}{$T_4$}
\psfrag{TS}{${T_4}^{\ast}$}
\centerline{\epsfig{file=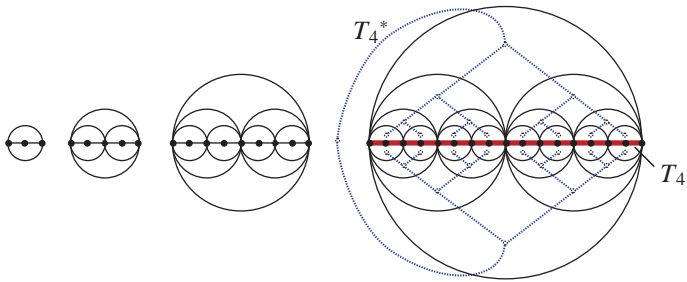}} 
\caption{The graphs $\Gamma_1, \ldots, \Gamma_4$ of Exercise~\ref{geodesic trees fail}.}
\label{dual trees ex fig}
\end{figure}

\begin{exer}  \label{geodesic trees fail}
Figure~\ref{dual trees ex fig} shows the first four of a family of connected (multi)graphs $\Gamma_n$. The horizontal
path through $\Gamma_n$ of length $2^n$ is a spanning tree $T_n$ in $\Gamma_n$ and  ${T_n}^{\ast}$ is a geodesic spanning tree in ${\Gamma_n}^{\ast}$ based at the vertex dual to the face \emph{at infinity} .  
Calculate $\Diam({T_n}^{\ast})$, $\Diam({\Gamma_n})$ and $\Diam({\Gamma_n}^{\ast})$ .
Find a spanning tree $S_n$ in $\Gamma_n$ such that $\Diam(S_n)$ and
$\Diam({S_n}^{\ast})$ are both at most a constant times $n$.
\end{exer}

\ms

\ni A \emph{diagram} is a finite planar contractible 2-complex -- in effect, a van~Kampen diagram bereft of its group theoretic content.  Question~\ref{H}  can be regarded as a question about diagram measurements because a multigraph is finite, planar and connected if and only it it is the 1-skeleton of a diagram.  Question~\ref{H} was recently resolved negatively in \cite{RT} as  we will explain.  Filling length plays a key role via the following result (compare \cite[5.C]{Gromov}).   

\begin{thm} \label{FLbyDGL}
For all diagrams $\Delta$ $$\FL(\Delta) \ \leq \ C (\DGL(\Delta) + \ell(\partial \Delta)),$$ where $C$ is a constant depending only on the maximum length of the boundary cycles of the 2-cells in $\Delta$.     
\end{thm}

\ms
\ni \emph{Proof  \textup{(}adapted from \textup{\cite{GR3}).}}
It suffices to show that given any spanning tree $T$ in $\Delta^{(1)}$, %
\begin{eqnarray}\label{corrected bound}
\FL(\Delta) \  \leq \ \Diam(T) \, +  \, 2\, \lambda \, \Diam(T^{\ast}) \, +
\,  \ell(\partial \Delta), 
\end{eqnarray}
where $\lambda$ is the maximum length of the boundary circuits of 2-cells in $\Delta$.

Let the base vertex $\star$ of $\Delta$ be the root of $T$ and the vertex $\star_{\infty}$ dual to the \emph{2-cell at infinity} be the root of $T^{\ast}$. 
Let $m$ be the number of edges in $T$.  Let $\gamma: [0,2m]
\to T$ be the edge-circuit in $T$ that starts from $\star$ and traverses every edge of $T$ twice, once
in each direction, running close to the anticlockwise boundary
loop of a small neighbourhood of $T$.  For $i=1,2,\ldots,2m$ let
$\gamma_i$ be the edge traversed by $\gamma
\restricted{[i-1,i]}$ and consider it directed from $\gamma(i-1)$ to $\gamma(i)$. 
Let $\tau_i$ be the geodesic in $T^{\ast}$ from
${\star}_{\infty}$ to the vertex $v_i^{\ast}$ dual to the 2-cell to the right of $\gamma_i$ (possibly $v_i^{\ast} = \star_{\infty}$). Let
$\overline{\tau_i}$ be the union of all the 2-cells dual to vertices on $\tau_i$.

\begin{figure}[ht]
\psfrag{T}{$T_4$}
\psfrag{TS}{${T_4}^{\ast}$}
\centerline{\epsfig{file=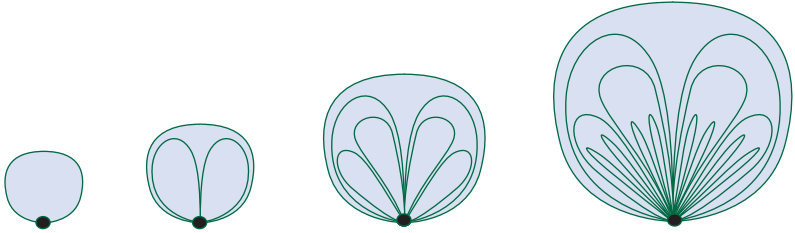}} 
\caption{Subdiagrams arising in the proof of Theorem~\ref{FLbyDGL}.}
\label{FL_examples fig}
\end{figure}

There may be some 2-cells $C$ in $\Delta \ssm \bigcup_i\overline{\tau_i}$.  For such a $C$, let $e$ be the edge of $(\partial  C) \ssm T$ whose dual edge is closest in $T^{\ast}$ to $\star_{\infty}$.   Call such an edge $e$ \emph{stray}. Then $e$ must be an edge-loop because otherwise the subdiagram with boundary circuit made up of $e$ and the geodesics in $T$ from the end vertices of $e$ to $\star$ would contain some $v_i^{\ast}$.  It follows that $e$ is the boundary of a subdiagram that is of the form of Figure~\ref{FL_examples fig} save that shelling away a 2-cell reveals at most $(\lambda -1)$ more 2-cells instead of two.  We leave it to the reader to check that such a diagram has filling length at most $(\lambda -1)\Diam(T^{\ast})$.  

To realise (\ref{corrected bound}), shell the 2-cells of $\Delta$ in the following order: for $i=1$, then $i=2$, and so on, shell all the (remaining) 2-cells along  $\overline{\tau_i}$, working from $\star_{\infty}$ to $v_i^{\ast}$.  Except, whenever a \emph{stray} edge $e$ appears in the boundary circuit, pause and entirely shell away the diagram it contains.  In the course of all the 2-cell collapses, do a \emph{1-cell-collapse} whenever one is available. 

One checks that in the course of this shelling, aside from detours into subdiagrams enclosed by stray edges, the anticlockwise boundary-circuit starting from $\star$ follows a geodesic path in $T$, 
then a path embedded in the 1-skeleton of some $\overline{\tau_i}$, before returning to $\star$
along $\partial \Delta$.   These three portions of the circuit have
lengths at most $\Diam(T)$, $\lambda \, \Diam(T^{\ast})$ and $\ell(\partial \Delta)$, respectively.  Adding $(\lambda -1)\Diam(T^{\ast})$ for the detour gives an estimate within the asserted bound.
\newqed

\ms




\begin{figure}[ht] 
\psfrag{p}{$p$}
\centerline{\epsfig{file=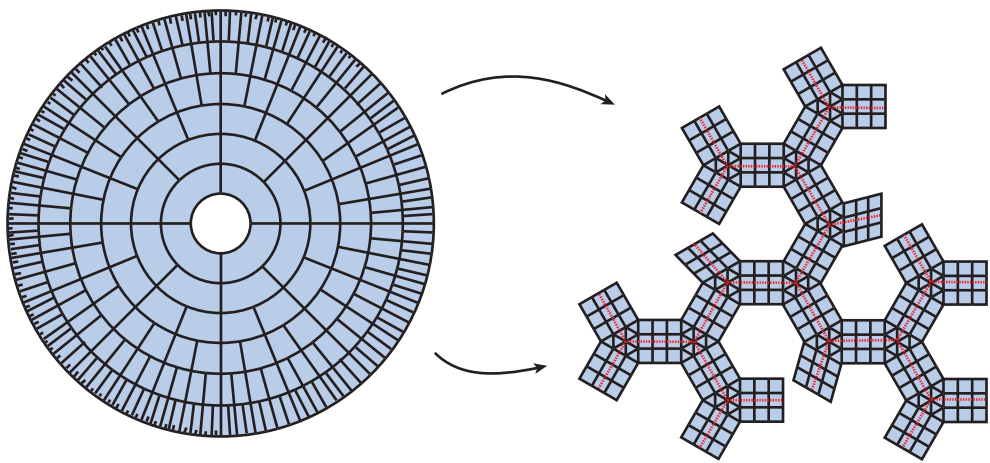}}
\caption{Attaching the outer boundary of the annulus to the boundary of the fattened tree gives the third of a family $(\Delta_n)$ of diagrams with filling length outgrowing diameter and dual diameter.} \label{fat trees} 
\end{figure}

\ni A family $(\Delta_n)$ of diagrams, essentially combinatorial analogues of metric 2-discs from \cite{Frankel-Katz}, is constructed as follows to resolve Question~\ref{H}.  One inductively defines trees $\mathcal{T}_n$ by taking $\mathcal{T}_1$ to be a single edge and obtaining $\mathcal{T}_n$ from three copies of $\mathcal{T}_{n-1}$ by identifying  a leaf vertex of each. (It is not important that this does not define  $\mathcal{T}_n$  uniquely.)  Fatten $\mathcal{T}_n$ to an $n$-thick diagram as shown in Figure~\ref{fat trees} in the case $n=3$.  Then attach a combinatorial annulus around the fattened tree by identifying the outer boundary of the annulus (in the sense illustrated) with the boundary of the fattened tree.  
The annulus is made up of concentric rings, each with twice as many faces as the next, and is constructed in such a way that the length of its outer boundary circuit is the same as the length of the boundary circuit of the fattened tree; it serves as \emph{hyperbolic skirt} bringing the diameters of ${\Delta_n}^{(1)}$ and its dual down to $\preceq n$.  
 
The filling length of $\Delta_n$ grows $\preceq n^2$ -- the reason being that a curve sweeping over $\mathcal{T}_n$ will, at some time, meet $n+1$ different edges of $\mathcal{T}_n$;  more precisely:

\begin{lem} \label{many intersections}
Suppose $\TT_n$ is embedded in the unit disc $\D^2$.  Let $\star$ be a basepoint on $\partial \D^2$.
Suppose $H : [0,1]^2 \to \D^2$ is a homotopy satisfying $H(0,t) = H(1,t) = \star$ for all $t$, and $H_0(s) = e^{2 \pi i s}$ and $H_1(s) = \star$ for all $s$,  where $H_t$ denotes the restriction of $H$ to $[0,1] \times \set{t}$.   Further, assume $H([0,1] \times [0,t]) \cap  H([0,1] \times [t,1]) = H([0,1] \times \set{t})$ for all $t$. Then $H_t$ meets at least $n+1$ edges in $\TT_n$ for some $t \in [0,1]$.
\end{lem}

\ni The final details of the proof that $(\Delta_n)$ resolves Question~\ref{H} negatively are left to the following exercise. 

\begin{exer}  \rule{0cm}{0cm} \\ \vspace*{-.5cm}
\begin{enumerate}
\item Show that $\Diam(\Delta_n), \GL(\Delta_n) \preceq n$.
\item Prove Lemma~\ref{many intersections} (\emph{hint}: induct on $n$) and deduce that $\FL(n) \succeq n^2$. 
\item Using Theorem~\ref{FLbyDGL}, deduce a negative answer to Question~\ref{H}.  
\end{enumerate}
\end{exer}

\ni The lengths of the boundary circuits in a $\PP$-van~Kampen diagram are at most the length of the longest defining relator in $\PP$.  So Theorem~\ref{FLbyDGL} shows that $\FL : \N \to \N$  and $\DGL : \N \to \N$ satisfy $\FL \preceq \DGL$.  The reverse, $\DGL \preceq \FL$, is also true for \emph{fat} presentations (defined in Section~\ref{combable section}).  The proof in \cite{GR3} is technical but, roughly speaking, the idea is that a shelling of $\Delta$ gives a family of concentric edge-circuits in $\Delta$ contracting down to $\star$ which we can \emph{fatten} and then have the arcs of both $T$ and $T^{\ast}$  follow these circuits.   Thus we get

\begin{thm}[\textup{\cite{GR3}}] \label{FL and DGL}  
The filling functions $\GL$, $\DGL$ and $\FL$ for any finite fat presentation
satisfy $\GL \leq \DGL \simeq \FL$. 
\end{thm}

\ni Despite the negative answer to Question~\ref{H}, whether the inequality in this theorem can be replaced by a $\simeq$, rendering all three functions equivalent, remains an open question.


\section{$\EDiam$ versus $\IDiam$}

It is clear that for a van~Kampen diagram $\pi: \Delta \to \Cay^2(\PP)$ we have $\EDiam(\Delta) \leq \IDiam(\Delta)$ because paths in $\Delta^{(1)}$ are sent by $\pi$ to paths in $\Cay^2(\PP)$.    But, in general, one expects a shortest path between $\pi(a)$ and $\pi(b)$ in the 1-skeleton of $\Cay^2(\PP)$ to take a shortcut through the space and so not lift to a path in $\Delta^{(1)}$, so  $\EDiam(\Delta)$ could be strictly less than  $\IDiam(\Delta)$.  In \cite{BR1} finite presentations are constructed that confirm this intuition; indeed, the examples show that (multiplicative) polynomial gaps of arbitrarily large degree occur:  

\begin{thm} \textup{\cite{BR1}} \label{diameters gap}
Given $\alpha >0$, there exists a finite presentation with
$$n^{\alpha}\,\EDiam(n)  \ \preceq \  \IDiam(n).$$
\end{thm}

\ni This result prompts the question of how far the gap can be extended:

\begin{Open Question} \cite{BR1}
What is the optimal upper bound for $\IDiam: \N \to \N$ in terms of $\EDiam: \N \to \N$, in general? 
\end{Open Question}

\section{Free filling length}

\begin{figure}[ht]
\psfrag{FL}{{$\FL(\Delta)$}}
\psfrag{f}{{$\FFL(\Delta)$}}
\psfrag{s}{$\star$}
\centerline{\epsfig{file=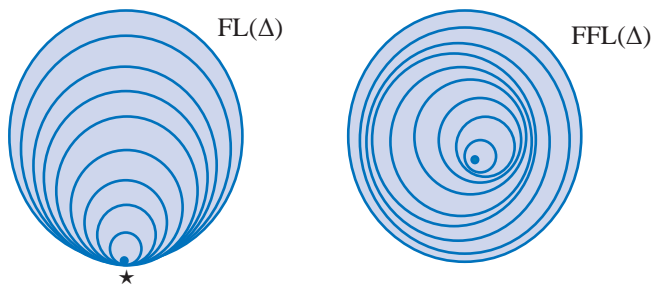}} 
\caption{Base-point-fixed vs.\ free filling length} \label{free fill}
\end{figure}

Recall from Section~\ref{via vK diagrams} that $\FL(\Delta)$ is defined with reference to a base vertex $\star$ on $\partial \Delta$.  If we remove the requirement that $\star$ be preserved throughout a shelling then we get $\FFL(\Delta)$: the minimal $L$ such that there is a \emph{free} combinatorial null-homotopy of $\partial \Delta$ across $\Delta$ -- see Figure~\ref{free fill}.  We then define $\FFL(w)$ for null-homotopic $w$, and the \emph{free filling length function} $\FFL: \N \to \N$ in the usual way.  

\ms

\begin{exer}
Re-interpret $\FFL(w)$ in terms of the Dehn Proof System.  (\emph{Hint.}  Allow one extra type of move to occur in null-sequences.) 
\end{exer}  

\ms

\ni It is questionable whether the base point plays a significant role in the definition of $\FL$: 

\begin{Open Problem} \cite{BR2}
Does there exist a finite presentation for which $\FL(n) \nsimeq \FFL(n)$? 
\end{Open Problem}

\ms

\ni However, the following example gives null-homotopic words $w_n$ for which $\FL(w_n)$ and $\FFL(w_n)$ exhibit dramatically different qualitative behaviour.   None-the-less  the filling functions agree: $\FL(n) \simeq \FFL(n) \simeq 2^n.$  (Exercise~\ref{ffl exer} provides some hints towards proof of these remarks; details are in \cite{BR2}.)   

\begin{ex}  \cite{BR2} \label{ffl ex}
Let $G$ be the group with presentation
$$\langle \ a,b, t, T, \tau \ \mid \ b^{-1}aba^{-2}, [t,a], [\tau, at], [T,t], [\tau,T]  \ \rangle$$ and define $w_n :=  \left[T, a^{-b^n} \tau  a^{b^n} \right].$ Then
$w_n$ is null-homotopic, has length $8n+8$, and $\FFL(w_n) \simeq n$, but $\FL(w_n) \simeq  2^n$. 
\end{ex}

\begin{figure}[ht]
\psfrag{a}{\footnotesize{$a$}}
\psfrag{b}{\footnotesize{$b$}}
\psfrag{T}{\footnotesize{$T$}}
\psfrag{t}{\footnotesize{$t$}}
\psfrag{tau}{\footnotesize{$\tau$}}
\psfrag{s}{\footnotesize{$\star$}}
\centerline{\epsfig{file=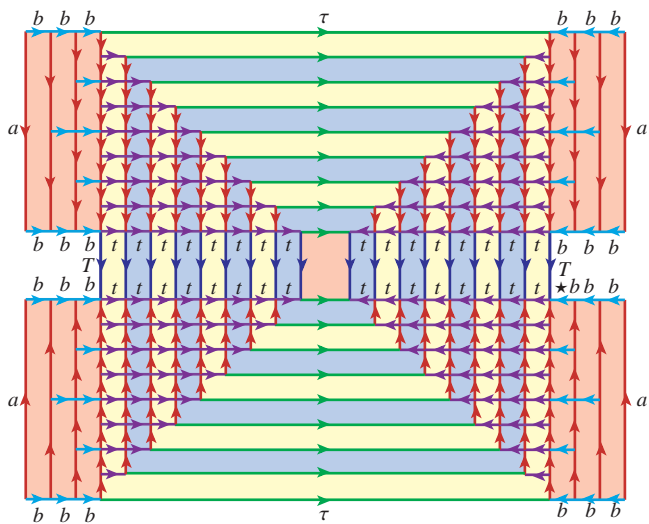}} 
\caption{A van Kampen diagram $\Delta_3$ for $w_3$ with concentric $t$-rings.}
\label{rings} 
\end{figure}

\begin{exer} \label{ffl exer} \rule{0cm}{0cm} \\ \vspace*{-.5cm}
\begin{enumerate}
\item  Show that $\FL(w_n) \simeq  2^n$.  (\emph{Hint.}  The van~Kampen diagram $\Delta_n$ for $w_n$ shown in Figure~\ref{rings} in the case $n=3$, has $2^n$ concentric \emph{$t$-rings} -- concatenations of 2-cells forming an annular subdiagram whose internal edges are all labelled by $t$.  Show that every van~Kampen diagram for $w_n$ has $t$-rings nested $2^n$ deep.)  

\item  Show that $\FFL(w_n) \simeq  n$. (\emph{Hint.} Produce a diagram $\Delta'_n$ from $\Delta_n$ such that $\FFL(\Delta'_n) \preceq n$ by  making cuts along paths labelled by powers of $a$ and gluing in suitable diagrams over the subpresentation $\langle a,b \ \mid \ b^{-1}aba^{-2} \rangle$.)  

\item Find a family of words $w'_n$ with $\FFL(w'_n) \simeq 2^n$ and $\ell(w'_n) \simeq n$.   
\end{enumerate}
\end{exer}

\chapter{Example: nilpotent groups}\thispagestyle{empty} \label{nilpotent}

Recall that a group $G$ is nilpotent of class $c$ if it has lower central series 
 $$G \ = \ \Gamma_1 \gneq \Gamma_2 \gneq  \cdots \gneq  \Gamma_{c+1} \ = \
\set{1},$$  defined inductively by $\Gamma_1 := G$ and
$\Gamma_{i+1} := [ G, \Gamma_i].$

\section{The Dehn and filling length functions}

The following theorem and its proof illustrate the value of two key ideas:  studying not only $\Area$ but how it interacts with other filling functions, and using the \emph{Dehn proof system} (see Section~\ref{Dehn proof system}) to manipulate words and monitor the Dehn and filling length function. 

\begin{thm} \label{c+1-theorem}
If $G$ is a finitely generated nilpotent group of class $c$ then the Dehn and filling length functions of $G$ satisfy $$\Area(n) \  \preceq \ n^{c+1} \ \  \text{ and } \ \ \FL(n) \  \preceq \ n.$$  Moreover, these bounds are simultaneously realisable. 
\end{thm}

\ni The combinatorial proof  in  \cite{GHR} is the culmination of a number of results:  \cite{Conner, Gersten6,  Hidber}.  A more general result was proved by Gromov \cite[$5.A'_5$]{Gromov}, \cite{Gromov6}  using masterful but formidable analytical techniques.  The linear upper bound on filling length was proved prior to Theorem~\ref{c+1-theorem} by the author in \cite{Riley} via asymptotic cones: finitely generated nilpotent groups have simply  connected  (indeed contractible) asymptotic cones by work of Pansu \cite{Pansu}, and so the result follows as we will see in Theorem~\ref{various bounds}.
Pittet \cite{Pittet}, following Gromov, proved that a lattice in a simply connected homogeneous nilpotent Lie group of class $c$ admits a polynomial isoperimetric function of degree
$c+1$. (A nilpotent Lie group is called \emph{homogeneous} if its Lie algebra is \emph{graded}.) 

No better general upper bound in terms of class is possible because, for example, the Dehn function of a free nilpotent group of class $c$ is $\simeq n^{c+1}$  -- see \cite{ShortNotes} or \cite{BMS, Gersten}. 
However, it is not best possible for individual nilpotent groups: the $2k+1$ dimensional integral Heisenberg groups are of class $2$ but have Dehn functions $\simeq n^2$ when $k>1$ \cite{Allcock}, \cite[$5.\textup{A}'_4$]{Gromov}, \cite{OS4}.  By way of contrast, the 3-dimensional integral Heisenberg group $\mathcal{H}_3$ is free nilpotent of class 2 and so has a cubic Dehn function.

\bs

\ni A full proof of Theorem~\ref{c+1-theorem} is beyond the scope of this text.  We will illustrate the ideas in the proof in \cite{GHR} in the case of $\mathcal{H}_3$, presented by $$\PP \  := \  \langle  \ x,y,z \ \mid  \ [x,y]z^{-1}, [x,z],  [y,z]  \ \rangle.$$  

\ni The proof is via the Dehn proof system.  We will show that there is a sequence $(w_i)$ that uses $\preceq n^{c+1}$ application-of-a-relator moves and has $\max_i \ell(w_i) \preceq n$.  This suffices -- see Section~\ref{Dehn proof system}.

 The following lemma gives a means of \emph{compressing} powers of $z$.  It is a special case of \cite[Corollary~3.2]{GHR}, which concerns non-identity elements in the final non-trivial group $\Gamma_c$ in the lower central series.  We use the terminology of Section~\ref{Dehn proof system}.  Fix $n \in \N$ (which in the proof will be the length of the word we seek to fill).  For integers $s \geq 0$, define \emph{compression words} $u(s)$ representing $z^s$ in $\mathcal{H}_3$:  if $0 \leq s \leq n^2-1$ then 
\begin{eqnarray*}
u(s)  & := &  z^{s_0}[x^n,y^{s_1}],
\end{eqnarray*}
where $s =s_0 + s_1n$ for some $s_0,s_1 \in \set{0, 1, \ldots, n-1}$;  define $$u \left( n^2 \right) \ :=  \  \left[x^n, y^n\right]$$ and $$u\left(A+Bn^2\right) \ := \ u(A) \, u\!\left(n^2\right)^B$$ for all integers $A,B$ with $0 \leq A \leq n^2
-1$ and $B > 0$.  Note that $\ell(u(A+Bn^2)) \leq K_0 n $ where $K_0$ depends only on $B$.

\begin{lem} \label{compression}
Fix $K_1>0$.  There exists $K_2>0$ such that for all integers $K_1 n^2 \geq s \geq 0$, there is a concatenation of  $\PP$-sequences $$z^{s} \  \to \ z^{s} \, u(0) \ \to \ z^{s-1} \, u(1) \ \to \ z^{s-2} \, u(2) \ \to \ \cdots \  \to \  z \, u(s-1)  \ \to \
u(s),$$ that converts $z^{s}$ to $u(s)$, and uses at most $K_2 n^{3}$ application-of-a-relator moves.
\end{lem}

\begin{figure}[ht]
\psfrag{s}{$\star$}
\psfrag{x}{$x^4$}
\psfrag{y}{$y^4$}
\psfrag{z}{$z$}
\centerline{\epsfig{file=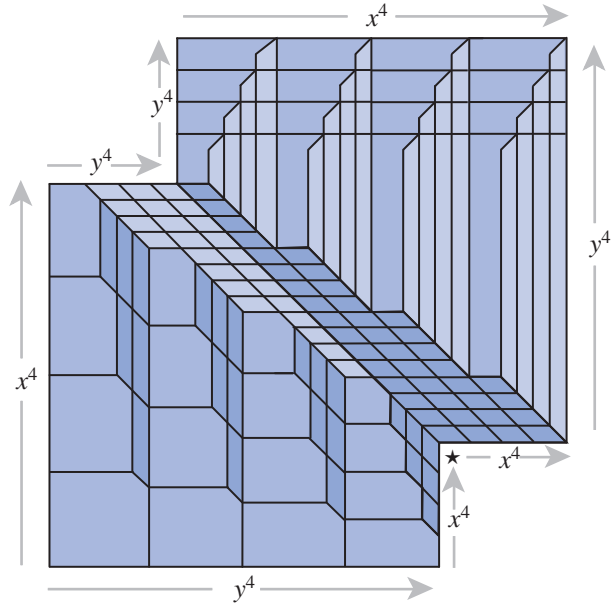}} 
\caption{A van~Kampen diagram for $[x^{-4}, y^{-4}][y^4, x^4] $ in the presentation $\langle x,y,z \mid [x,y]z^{-1}, [x,z],  [y,z] \rangle$ for the 3-dimensional integral Heisenberg group $\mathcal{H}_3$.  The diagonal line beginning at $\star$ is labelled by $z^{16}$.} \label{Heisenberg fig} 
\end{figure}

\ni  We leave the proof of the lemma as an exercise except to say that the key is absorbing $z^n$ subwords into commutators $\left[ x^n , y^k \right]$ in which $k<n$, via:  
$$z^n \left[ x^n , y^k \right]  \ = \  z^n x^{-n} y^{-k} x^n y^k  \ = \   z^n x^{-n} y^{-k} y^{-1} y x^n y^k \ = \   z^n x^{-n} y^{-(k+1)}   (xz^{-1})^n y^{k+1} \ = \ \left[ x^n , y^{k+1} \right].$$ 
A diagrammatic understanding of this calculation can be extracted from Figure~\ref{Heisenberg fig} -- the diagonal line running across the diagram from $\star$ is labelled by $z^{16}$ and the half of the diagram below this diagonal demonstrates the equality  $z^{16} = [x^4,y^4]$.

\ni The proof  of Theorem~\ref{c+1-theorem} is an induction on the class.  The base case of class $1$ concerns finitely generated abelian groups $G$ and is straightforward -- see Example~\ref{examples comp}~(2).  For the induction step one considers $G / \Gamma_c$, which is a finitely generated nilpotent group of class $c-1$.   Let $\overline{\PP} = \langle \overline{X} \mid \overline{R} \rangle$ be a finite presentation for $G / \Gamma_c$. 
Let $X = \overline{X} \cup X'$ where $X'$ is a finite set of  length $c$ commutators that span $\Gamma_c$.  The relators $\overline{R}$ represent elements of $\Gamma_c$ in $G$, so there is a finite presentation for $\langle X \mid R \rangle$ in which the elements of $R$ express the equality of elements of $\overline{R}$ with words on $X'$.  

Given a length $n$ null-homotopic word $w$ in $\PP$, the word $\overline{w}$ obtained from $w$ by deleting all letters in $X'$ is null-homotopic in $\overline{G}$ and so admits a null-$\overline{\PP}$-sequence $(\overline{w}_i)_{i=0}^{\overline{m}}$ which, by induction, involves words of length $\preceq n$ and uses $\preceq n^c$ application-of-a-relator moves.  

We seek to lift  $(\overline{w}_i)_{i=0}^{\overline{m}}$  to a null-$\PP$-sequence for $w$ in which the words have length $\preceq n$ and  $\preceq n^{c+1}$ application-of-a-relator moves are used.   
We illustrate how this works in the case of $G=\mathcal{H}_3$. We have  $$\overline{G} \ = \ \Gamma_1 / \Gamma_2  \ =  \  \mathcal{H}_3 / \langle z \rangle,$$ which is presented by $\overline{\PP} = \langle x,y  \mid [x,y] \rangle$.  Start with $w$ and use commutator relations to shuffle all $z$ and $z^{-1}$ to the right and left ends of the word, respectively, leaving $\overline{w}=\overline{w}_0$ in the middle of the word.  Each move $\overline{w}_i \mapsto \overline{w}_{i+1}$ that applies the relator $[x,y]$  lifts to a move that introduces a $z^{\pm 1}$ because $[x,y] \in \overline{R}$ lifts to $[x,y]z^{-1} \in R$.   Whenever a $z$ or $z^{-1}$ appears in this way, immediately shuffle it to the right or left end of the word, respectively.  As they arrive, compress the letters $z$ at the right end of the word as per Lemma~\ref{compression}, and compress the $z^{-1}$ at the left end by using the $\PP$-sequence  obtained by inverting every word in the  $\PP$-sequence of the lemma.  Once  $(\overline{w}_i)_{i=0}^{\overline{m}}$ has been exhausted we have a compression word at the right end of the word and its inverse at the left end, and the empty word  $\overline{w}_{\overline{m}}$ in between.  Freely reduce to obtain the empty word.
   
The total number of $z^{\pm 1}$ being sent to each end in this process is $\preceq n^2$; the lengths of the words stay $\preceq n$ on account of the compression process, Lemma~\ref{compression}, and the induction hypothesis; and, as required, the total number of application-of-a-relator moves used is $\preceq n^3$, the dominant term in the estimate coming from Lemma~\ref{compression}.  Note that the reason for collecting $z$s and $z^{-1}$s at different ends of the word is that the compression process of  Lemma~\ref{compression} has some steps that are more expensive that others, and having the power of $z$ (or $z^{-1}$) monotonically increasing avoids crossing these thresholds unduly often.

\section{Open questions}
 
By Theorem~\ref{Combable groups}, an affirmative answer to the following question would be a stronger result than Theorem~\ref{c+1-theorem}.

\begin{Open Question}
Do all class $c$ finitely generated nilpotent groups admit (a)synchronous combings with length functions $\textup{L}(n) \preceq n^c$?  (Cf.\ \cite{Gromov, Pittet}.) 
\end{Open Question}

\ni The 3-dimensional integral Heisenberg group cannot have a synchronous combing with length function growing at most linearly fast because that could contradict the $n^3$ lower bound on its Dehn function \cite{BMS, Gersten, ShortNotes}.  This makes it a candidate for a natural example of a group exhibiting the behaviour of the example of Bridson discussed in Remark~\ref{Bridson combing}.

\begin{Open Question} \label{Heis qn etc}
Is the 3-dimensional integral Heisenberg group synchronously combable?  More generally, which nilpotent groups are synchronously combable?
\end{Open Question}

\ni  In reference to the second part of \ref{Heis qn etc} we mention that Holt~\cite{Epstein} showed that the finitely generated nilpotent groups that are \emph{a}synchronously combable are those that are virtually abelian.    
Our next question shows how little we yet know about Dehn functions of nilpotent groups.

\begin{Open Question} \cite[Question~N2]{BMSqus}
Does there exist a finitely generated nilpotent group whose Dehn function is not $\simeq$-equivalent to $n^{\alpha}$ for some integer $\alpha$? 
\end{Open Question}

\ni Young \cite[Section~5]{Young} describes a nilpotent group put forward by Sapir as a candidate to have Dehn function $\simeq n^2 \log n$.  The upper bound has been established; the challenge of proving it to be sharp remains.

\begin{Open Question} \label{linear filling length qu}
Do linear upper bounds on  the filling length functions for nilpotent groups extend to 
finite presentations of polycyclic groups or, more generally, linear groups?
\end{Open Question}

\ni Progress (in the affirmative direction) towards \ref{linear filling length qu} would, by Proposition~\ref{space-time bound},  increase the evidence for the following.  

\begin{Conjecture} \textbf{\textup{(Gersten \textup{\cite{GR4}})}} 
There is a recursive function $f: \N \to \N$ (e.g.\ $f(n)=2^n$) such that finite presentations of linear groups all have $\Area(n) \preceq f(n)$. 
\end{Conjecture}

\chapter{Asymptotic cones}\thispagestyle{empty} \label{cones}

\section{The definition}

In \cite{Gromov} Gromov says of a finitely generated group $G$
with a word metric $d$:
\begin{quote}
``\emph{This space may at first appear boring and uneventful to a
geometer's eye since it is discrete and the traditional local
\emph{(}e.g.\ topological and infinitesimal\emph{)} machinery does not
run in $G$.}''
\end{quote}
Imagine viewing $G$ from increasingly distant vantage points, i.e.\ scaling the metric by a sequence $\mathbf{s}=(s_i)$ with $s_i \to \infty$.  An asymptotic cone is a limit of $\left(G,\frac{1}{s_n} d\right)$, a coalescing of $G$ to a more continuous object that is amenable to attack by ``\emph{topological and infinitesimal machinery}'', and which ``\emph{fills our geometer's heart with joy}'' (\cite{Gromov} again).  

The asymptotic cones $\Cone$ of a metric space $(\XX,d)$ are defined with reference to three ingredients: the first is a sequence $\e= (e_i)$ of basepoints in $\XX$ (if $\XX$ is a group then, by homogeneity, all the $e_i$ may as well be the identity element); the second is a sequence of strictly positive real numbers $\mathbf{s} = (s_i)$ such that $s_i \to \infty$ as $i \to \infty$;  and the third is a \emph{non-principal ultrafilter} on $\N$, the magic that forces convergence.  Van den Dries and Wilkie \cite{Wilkie} recognised the usefulness of non-principal ultrafilters for cutting through delicate arguments concerning extracting a convergent subsequence (with respect to the Gromov-Hausdorff distance)  from a sequence of metric spaces  in Gromov's proof  \cite{Gromov2} that groups of polynomial growth are virtually nilpotent.  

\begin{defn}
A \emph{non-principal ultrafilter} is a finitely additive probability measure on $\N$ that takes values in $\set{0,1}$ and gives all singleton sets measure $0$. 
\end{defn}

\ni The existence of non-principal ultrafilters is equivalent to the Axiom of Choice (see Exercises~\ref{Ultrafilters Remarks2}) and so we can be barred from describing asymptotic cones explicitly and most proofs involving asymptotic cones are non-constructive. Generally, asymptotic cones are wild beasts: the only finitely generated groups whose asymptotic cones are locally compact metric spaces are virtually nilpotent groups.  [En route to his Polynomial Growth Theorem \cite{Gromov2}, Gromov proves that the asymptotic cones of virtually nilpotent groups are proper.  Dru\c{t}u proves the converse in \cite{DrutuSurvey} and she adds \cite{DrutuPersonal} that by the Hopf-Rinow Theorem (see, for example, \cite[Proposition~3.7]{BrH}) proper can be replaced by \emph{locally compact} as asymptotic cones of groups are complete (see \cite{Wilkie}) geodesic spaces (see Exercise~\ref{Basic cones properties}).]  Indeed, some asymptotic cones contain $\pi_1$-injective copies of the Hawaiian earring and (so) have uncountable fundamental groups \cite{Burillo}. But rather than being put off, one should see these features as part of the spice of the subject.

The way a non-principal ultrafilter $\omega : \PP(\N) \to \set{0,1}$ forces convergence is by \emph{coherently} (see Exercise~\ref{Ultrafilters Remarks2}~(4), below) choosing limit points from sequences of real numbers: given a sequence $(a_{i})$ in $\mathbb{R}$ we say that  $a \in \mathbb{R}$ is an $\omega $-\emph{ultralimit} of $(a_{i})$ when $\forall \varepsilon
>0,~\omega \set{i \mid \, \abs{a-a_{i}} <\varepsilon}=1$, that $\infty$ is an $\omega $-ultralimit of $(a_{i})$ when $\omega \set{i \mid a_{i} >N}=1$ for all $N >0$, and that $-\infty$
is an $\omega$-ultralimit when $\omega
\set{i \mid a_{i} <-N}=1$ for all $N >0$.

\ms
 
\begin{exer} \label{Ultrafilters Remarks2}  \rule{0cm}{0cm} \\ \vspace*{-.5cm}
\begin{enumerate}
\item Show that if $\omega$ is a non-principal ultrafilter on $\N$ and $A, B \subseteq \N$ satisfy $\omega(A) = \omega(B) =1$ then $\omega(A \cap B) =1$.  

\item Establish the existence of non-principal ultrafilters $\omega$ on $\N$.  (\emph{Hint.} Consider the set of functions $\omega: \PP(\N) \to \set{0,1}$ such that $\omega^{-1}(1)$ is closed under taking intersections and taking supersets, and does not include the empty set, but does include all  $A \subseteq \N$ for which $\N \ssm A$ is finite.  Use Zorn's Lemma.)

\item  Let $\omega$ be a non-principal ultrafilter.  Prove that every $\omega$-ultralimit is also a limit point in the usual sense.

\item  Show that every sequence $(a_{n})$ of real numbers has a unique $\omega$-ultralimit in $\mathbb{R} \cup \set{ \pm \infty}$ denoted $\lim_{\omega
}a_{n}$. (\emph{Hint.} First show that if $+\!\infty$ or $-\!\infty$  is an
ultralimit then it is the unique ultralimit.  Next assume $\pm \infty$  is not an
ultralimit and prove the existence of an ultralimit by suitably adapting a proof that every bounded sequence of real numbers has a limit point. 
Finally prove uniqueness.)

\item Let $\lim_{\omega}a_i$ denote the $\omega$-ultralimit of $(a_i)$.  Prove that for all $\lambda, \mu \in \mathbb{R}$ and sequences $(a_i), (b_i)$ $${\lim_{\omega} ( \lambda a_i + \mu b_i) \ = \  \lambda \lim_{\omega}  a_i + \mu \lim_{\omega} b_i,}$$ and $\lim_{\omega} a_i < \lim_{\omega} b_i$ if and only if $a_i < b_i$ for all $n$ in a set of $\omega$-measure 1. 

\item Suppose $f: \mathbb{R} \to \mathbb{R}$ is continuous.  Show that for sequences $(a_i)$ of real numbers, $f(\lim_{\omega} a_i) = \lim_{\omega} f(a_i)$.      
\end{enumerate}
\end{exer}

\ni Points in $\Cone$ are equivalence classes of sequences $(x_i)$ in $\XX$ such that $\lim_{\omega} d(e_i,x_i)/s_i < \infty$,  where two sequences $(x_i)$, $(y_i)$ are equivalent if and only if $\lim_{\omega} d(x_i,y_i)/s_i =0$.  The metric on $\Cone$ is also denoted by $d$ and is given by $d(\x, \y) = \lim_{\omega} d(x_i, y_i)/s_i$ where $(x_i)$ and $(y_i)$ are sequences representing $\x$ and $\y$. 

\ms

\begin{exer}
Check $d$ is well-defined and is a metric on $\Cone$.   
\end{exer}

\ms

\ni The next exercise gives some important basic properties of asymptotic cones.  In particular, it shows that if a metric space is $(1,\mu)$-quasi-isometric to a geodesic metric space (as is the case for a group with a word metric, for example) then its asymptotic cones are all geodesic metric spaces.  [A \emph{geodesic} in a metric space $(\mathcal{X},d)$ is an isometric embedding $\gamma: I \to \mathcal{X}$, where $I \subseteq \mathbb{R}$ is a closed interval (possibly infinite or bi-infinite).  We say that $(\mathcal{X},d)$ is a \emph{geodesic space} if every two points in $\mathcal{X}$ are connected by a geodesic.]

\smallskip

\begin{exer}
 \label{Basic cones properties}  \rule{0cm}{0cm} \\ \vspace*{-.5cm}
\begin{enumerate}
\item Show that a $(\lambda, \mu)$-quasi-isometry $\Phi: \mathcal{X} \to \mathcal{Y}$  between metric spaces induces $\lambda$-bi-Lipschitz homeomorphisms $$\textup{Cone}_{\omega}(\mathcal{X}, (e_i),  \textbf{s}) \to \textup{Cone}_{\omega}(\mathcal{Y}, (\Phi (e_i)),  \textbf{s}),$$
for all $\omega, (e_i)$ and $\textbf{s}$. 

\item Show that if $\textbf{a} = (a_i)$ and $\textbf{b} = (b_i)$ represent points in $\textup{Cone}_{\omega}(\mathcal{X}, \textbf{e},  \textbf{s})$  and $\gamma_i : [0,1] \to \mathcal{X}$ are geodesics from $a_i$ to $b_i$, parametrised proportional to arc length, then $\pmb{\gamma} : [0,1] \to \textup{Cone}_{\omega}(\mathcal{X}, \textbf{e},  \textbf{s})$ defined by $\pmb{\gamma}(t) = (\gamma_i(t))$ is a geodesic from $\textbf{a}$ to $\textbf{b}$.
\end{enumerate}
\end{exer}

\ms 

\ni It is natural to ask how asymptotic cones depend on $\omega$ and $\mathbf{s}$.  Varying one is similar to varying the other (see \cite[Appendix~B]{Riley3}) but the precise relationship is not clear.  In this chapter we will discuss results about statements about all the asymptotic cones of a group for a fixed $\omega$ but  $\mathbf{s}$ varying.  (Other authors fix $\mathbf{s}$ as $(s_i) = (i)$ and vary $\omega$, and others allow both to vary.)  This will be crucial because another quirk of the subject is that the topological type of the asymptotic cones of a group, even, may depend on $\omega$ or $\mathbf{s}$.  
Thomas \& Velickovic \cite{Thomas} found the first example of a finitely generated group with two non-homeomorphic asymptotic cones -- they used an infinite sequence of defining relations satisfying small cancellation to create holes in the Cayley graph on an infinite sequence of scales, and then, depending on whether or not its ultrafilter caused the asymptotic cone to see a similar sequence of scales, it either has non-trivial fundamental group, or is an $\mathbb{R}$-tree.  Later Kramer, Shelah, Tent \& Thomas~\cite{KSTT} found an example of a finitely presented group with the  mind-boggling property of having $2^{2^{\aleph_0}}$ non-homeomorphic cones if the Continuum Hypothesis (CH) is false but only one if it is true.  Also they showed that under CH a finitely generated group has at most $2^{\aleph_0}$ non-homeomorphic asymptotic cones, a result Dru\c{t}u \& Sapir \cite{DS} proved to be sharp when they found a finitely generated group which (independent of CH) has $2^{\aleph_0}$ non-homeomorphic cones.   
Most recently,  Sapir \& Ol'shanskii \cite{OS1, Ol} constructed a finitely presented group for which the vanishing of the fundamental groups of their asymptotic cones depends on $\mathbf{s}$ (independent of CH). 


\section{Hyperbolic groups} \label{hyperbolic section}

A metric space $(\XX,d)$ is \emph{$\delta$-hyperbolic} in the sense of Gromov~\cite{Gromov4}  if for all $w,x,y,z \in \XX$, $$d(x,w) + d(y,z) \ \leq \ \max\set{d(x,y) + d(z,w), d(x,z) + d(y,w)} + \delta.$$  We say $(\XX,d)$ is \emph{hyperbolic} when it is $\delta$-hyperbolic  for some $\delta \geq 0$. 

This \emph{four-point condition} for hyperbolicity is one of a number of equivalent formulations.  A geodesic metric space $(\XX,d)$ is \emph{hyperbolic} if and only if it satisfies the \emph{thin-triangles condition}:  there exists $\delta>0$ such that, given any geodesic triangle in $\XX$,  any one side is contained in a $\delta$-neighbourhood of the other two.   

\begin{exer}
Relate the $\delta$ that occurs in the  thin-triangles condition to the $\delta$ in the four-point condition. 
\end{exer}

\ms

\ni We say that a group $G$ with finite generating set $X$ is hyperbolic when $\Cay^1(G,X)$ is a hyperbolic metric space or, equivalently, $(G, d_X)$ satisfies the four-point condition for some $\delta \geq 0$.

\ms

\ni For a van~Kampen diagram $\Delta$ over a  finite presentation $\PP$, define $$\Rad(\Delta)  \  := \max\set{ \ \rho(a, \partial \Delta)  \ \mid \ \text{vertices } a \text{ of } \Delta \ },$$ where $\rho$ is the combinatorial metric on $\Delta^{(1)}$.   For a null-homotopic word  define $\Rad(w)$ be the minimum of $\Rad(\Delta)$ amongst all van~Kampen diagrams for $w$ in the usual way, but also define $\overline{ \Rad}(w)$ to be the maximum of $\Rad(\Delta)$ amongst all minimal area diagrams for $w$.  Then, as usual, define the corresponding filling functions, the \emph{radius} and \emph{upper radius} functions, $\Rad, \overline{\Rad} : \N \to \N$ for  $\PP$ to be the maxima of $\Rad(w)$ and $\overline{\Rad}(w)$, respectively, over all null-homotopic words of length at most $n$.  If $\Rad(n) \succeq n$ then $\IDiam (n ) \simeq \Rad(n)$, but in hyperbolic groups the extra sensitivity of $\Rad$ for low radius diagrams is useful, as we will see in the following theorem, which gives formulations of hyperbolicity for groups in terms of Dehn functions, radius functions, and asymptotic cones.  A geodesic metric space $\mathcal{X}$ is an $\mathbb{R}$-tree when each pair of points $a,b \in \mathcal{X}$ is joined by a unique geodesic segment $[a,b]$, and if $[a,b]$ and $[b,c]$ are two geodesic segments with unique intersection point $b$, then $[a,c]$ is the concatenation of  $[a,b]$ and $[b,c]$.  

\begin{thm} \label{R-trees} \label{hyperbolicity via Dehn fns}
Let $\omega$ be a non-principal ultra-filter on $\N$.  For a group $G$ with finite  generating set $X$, the following are equivalent. 
\renewcommand{\theenumi}{\textup{(\arabic{enumi})}}
\begin{enumerate}
\item \label{hyp} $G$ is hyperbolic.
\item \label{trees} For all sequences of real numbers $\mathbf{s}=(s_i)$ with $s_i \to \infty$, the cones $\textup{Cone}_{\omega}(G, \pmb{1}, \mathbf{s})$ are $\mathbb{R}$-trees.
\item \label{Dehn pres cond} $G$ admits a Dehn presentation \textup{(}as defined in Section~\ref{Dehn proof system}\textup{)}.
\item  \label{linear isop} There exists $C>0$ and a finite presentation for $G$ with Dehn function satisfying $\Area(n) \leq C n$ for all $n$.
\item  \label{subquadratic} If \: \!\! $\PP$ is a finite presentation for $G$ then the Dehn function of $\PP$ satisfies $\Area(n)/n^2 \to 0$ as $n \to \infty$. 
\item  \label{log upper rad} There exists $C>0$ and a finite presentation for $G$ with  $\overline{\Rad}(n) \leq C \log n$ for all $n$.
\item \label{log rad}  There exists $C>0$ and a finite presentation for $G$ with  $\Rad(n) \leq C \log n$ for all $n$.
\item  \label{sublinear rad} If \: \!\! $\PP$ is a finite presentation for $G$ then the radius function of $\PP$ satisfies $\Rad(n)/n \to 0$ as $n \to \infty$. 
\end{enumerate}
\renewcommand{\theenumi}{\arabic{enumi}}
\end{thm}

\ni In fact, all the asymptotic cones of non-virtually-cyclic hyperbolic groups are the same: they are the universal $\mathbb{R}$-trees with $2^{\aleph_0}$ branching at every vertex \cite{DyubinaPolterovich}.     The example \cite{Thomas}, discussed earlier, shows that quantifying over all $\mathbf{s}$ is necessary in Theorem~\ref{R-trees}~\ref{trees}.  However,  M.Kapovich \& Kleiner \cite{KK} have shown that if a \emph{finitely presented} group $G$ has one asymptotic cone that is a $\mathbb{R}$-tree then $G$ is hyperbolic.  

\ms

\ni \emph{Proof of Theorem~\ref{R-trees}.}  We begin by proving the equivalence of \ref{hyp} and \ref{trees}, which is due to Gromov~\cite{Gromov4}; subsequent accounts are \cite{Drutu2}, \cite{DrutuThesis}, and \cite[Chapter~2, \S1]{GhysHarpe}.  The proof here is based closely on that of Dru\c{t}u and, in fact, amounts to not just a group theoretic result but to a characterisation of hyperbolic geodesic metric spaces.

\smallskip

\ni \ref{hyp} $\implies$ \ref{trees}. Assume $G$ is hyperbolic.  From Exercise~\ref{Basic cones properties} we know $\textup{Cone}_{\omega}(G, \pmb{1}, \mathbf{s})$ is a geodesic space.  One can show  that any four points $\textbf{w}, \textbf{x}, \textbf{y}, \textbf{z} \in  \textup{Cone}_{\omega}(G, \pmb{1},  \mathbf{s})$ satisfy the four point condition with $\delta =0$ by applying 
the four point condition to $w_i, x_i, y_i, z_i$,  where   $(w_i), (x_i), (y_i), (z_i)$ are representatives for $\textbf{w}, \textbf{x}, \textbf{y}, \textbf{z}$, and using properties of $\lim_w$ such as those in Exercise~\ref{Ultrafilters Remarks2}~(4).   It follows that all geodesic triangles are in $\textup{Cone}_{\omega}(G, \pmb{1}, s_i)$ are $0$-thin and from that one can deduce that  
$\textup{Cone}_{\omega}(G, \pmb{1}, s_i)$ is an $\mathbb{R}$-tree.  

\ms

\ni \ref{trees} $\implies$ \ref{hyp}.  Assume $\textup{Cone}_{\omega}(G, \pmb{1}, \mathbf{s})$ are $\mathbb{R}$-trees for all $\mathbf{s}$. Let $d$ denote the path metric on $\Cay^1(G,X)$.  Suppose, for a contradiction, there is no $\delta >0$ such that all the geodesic triangles in $\Cay^1(G,X)$ are $\delta$-thin.  Then there are geodesic triangles $[x_i, y_i, z_i]$  in $\Cay^1(G,X)$ such that, defining $s_i$ to be the infimal distance such that every point on any one side of $[x_i, y_i, z_i]$ is in an $s_i$-neighbourhood of the other two, we have $s_i \to \infty$.  

\begin{figure}[ht]
\psfrag{s}{$\star$}
\psfrag{x}{$\textbf{x}$}
\psfrag{y}{$\textbf{y}$}
\psfrag{z}{$\textbf{z}$}
\psfrag{xn}{$x_i$}
\psfrag{yn}{$y_i$}
\psfrag{zn}{$z_i$}
\psfrag{an}{$a_i$}
\psfrag{bn}{$b_i$}
\psfrag{cn}{$c_i$}
\psfrag{sn}{$s_i$}
\psfrag{tn}{$t_i$}
\psfrag{a}{$\textbf{a}$}
\psfrag{b}{$\textbf{b}$}
\psfrag{c}{$\textbf{c}$}
\psfrag{A}{$A$}
\psfrag{B}{$B$}
\psfrag{C}{$C$}
\psfrag{p}{$p$}
\psfrag{q}{$q$}
\psfrag{1}{$1$}
\centerline{\epsfig{file=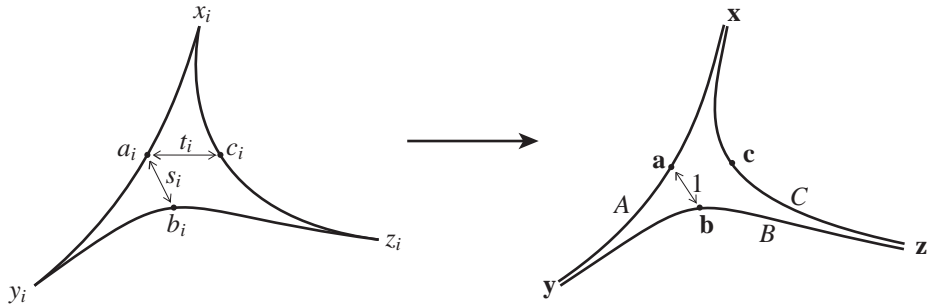}} 
\caption{The geodesic triangle $[x_i, y_i, z_i]$.} \label{labelled triangle} 
\end{figure}

By compactness and by relabelling if necessary, we can assume that $s_i =  d(a_i ,b_i)$ for some $a_i \in [x_i,y_i]$ and $b_i \in [y_i, z_i]$.  Then $$t_i \  := \ d(a_i, [x_i, z_i]) \  \geq \  s_i $$ and $t_i = d(a_i, c_i)$ for some $c_i \in [x_i , z_i]$ -- see Figure~\ref{labelled triangle}.  
Let $\textbf{a} = (a_i)$ and  $\textbf{s} = (s_i)$.  Define $\CC: = \textup{Cone}_{\omega}(\Cay^1(G,X) , \textbf{a}, \textbf{s})$, which is an $\mathbb{R}$-tree by hypothesis.  

\ni Let $\gamma_i$ be a geodesic running from $y_i$ to $x_i$ via $a_i = \gamma_i(0)$ at constant speed $s_i$ -- that is, $d(a_i , \gamma_i(r)) = \abs{r} s_i$.  Define $k_1 :=  \lim_{\omega} d(a_i ,  x_i)/s_i$ and $k_2 :=  \lim_{\omega} d(a_i ,  y_i)/s_i$, which will be $\infty$ if $\mathbf{x}$ or $\mathbf{y}$ (respectively) fail to define a point in $\mathcal{C}$.  If $k_1, k_2 < \infty$ then let $\pmb{\gamma} : [-k_2,k_1] \to \CC$ be the unit speed geodesic defined by 
\begin{eqnarray*}
\pmb{\gamma}(r) & := & (\gamma_i(r)) \ \ \textup{ for } -k_2 < r < k_1, \\ 
\pmb{\gamma}(k_1) & := & \mathbf{x},  \ \ \  \pmb{\gamma}(-k_2) \ := \ \mathbf{y}.
\end{eqnarray*} 
(For $-k_1 < r< k_2$, the expression $\gamma_i(r)$ is well defined for all $n$ in a set of $\omega$-measure $1$, and this is enough for  $(\gamma_i(r))$ to define a point in $\CC$.) 
If either or both of $k_1$ and $k_2$ is infinite then let $\pmb{\gamma}$ be an infinite or bi-infinite unit-speed geodesic $\pmb{\gamma}(r) = ( \gamma_i(r))$.   Let $A$ denote the image of $\pmb{\gamma}$.   In the same way, define $B$ to be (the image of) the geodesic  through $\mathbf{b} = (b_i)$ and between $\mathbf{y}$ and $\mathbf{z}$.  

Suppose $\mathbf{c} = (c_i)$ is a well defined point in $\CC$ (that is,  $\lim_{\omega} t_i/s_i < \infty$). Then, as before, define $C$ to be the geodesic through $\mathbf{c}$ and between $\mathbf{x}$ and $\mathbf{z}$.  
We claim that every point $\mathbf{p} = (p_i)$ on one of $A,B,C$ is in a $1$-neighbourhood of the other two.
Then, as $\CC$ is an $\mathbb{R}$-tree, the only way for $A,B$ and $C$ to remain close in this way, will be for $A \cup B \cup C$ to form a tripod (whose sides may be infinite), but that contradicts $d(\mathbf{a}, B \cup C) =1$.  

Suppose $\mathbf{p}$ is on $A$. In other words, $\mathbf{p} = \pmb{\gamma}(r) = (\gamma_i(r))$ for some $r$.  We will show $d(p, B \cup C) \leq 1$.   Similar arguments, which we omit, apply to $\mathbf{p}$ in other locations.  Let $q_i$ be a point on $ [x_i, z_i]\cup [y_i,z_i] $ closest to $p_i$.  Then $d(p_i, q_i) \leq s_i$ and so $d(\mathbf{p},  \mathbf{q}) \leq 1$ where $\mathbf{q} = (q_i)$.  We will show that $\mathbf{q} \in B \cup C$.

Let $\mathcal{S} = \set{i \mid q_i \in [z_i, x_i]}$ and $\mathcal{S}' = \set{i \mid q_i \in [y_i, z_i]}$. Suppose $\mathcal{S}$ has $\omega$-measure $1$.  Let $\pmb{\gamma}'$ be the (possibly infinite or bi-infinite) unit speed geodesic in $\CC$ with $\pmb{\gamma}'(0) = \mathbf{c}$, running between $\mathbf{z}$ and $\mathbf{x}$, and defined in the same way as $\pmb{\gamma}$.   
Now $d(p_i, q_i) \leq s_i$ and $$d(c_i, q_i) \ \leq \ d(c_i, a_i) + d(a_i, p_i) + d(p_i, q_i) \ \leq \ t_i + \abs{r} s_i + s_i.$$    So for all $i \in \mathcal{S}$ we find $q_i = \gamma'(r_i)$ for some $r_i$ with $\abs{r_i} \leq t_i/s_i +\abs{r}+1$ and therefore $\abs{\lim_{\omega} r_i} < \infty$.  It follows that $(\gamma'(r_i)) = (\gamma'(\lim_{\omega} r_i))$ in $\mathcal{C}$ and hence that $\mathbf{q}  \in C$, as required.  
If $\mathcal{S}$ has $\omega$-measure $0$ then $\mathcal{S}'$ has $\omega$-measure $1$ and a similar argument shows $\mathbf{q} \in B$.  

On the other hand, suppose $\mathbf{c}$ is not a well-defined point in $\CC$ or, equivalently, $\lim_{\omega} t_i/s_i = \infty$.  Then $A$ and $B$ are either both infinite or both bi-infinite as $\mathbf{x}$ and $\mathbf{z}$ cannot define points in $\mathcal{C}$.  We will show that every point on $A$ is a distance at most $1$ from $B$.  In an $\mathbb{R}$-tree this is not compatible with $d(\mathbf{a}, B) =1$.

Suppose $\mathbf{p} = (p_i)$ is a point on $A$.  Then  $\mathbf{p} = \pmb{\gamma}(r) = (\gamma_i(r))$ for some $r \in \mathbb{R}$.
As before, let $q_i$ be a point on $ [x_i, z_i] \cup [y_i,z_i] $ closest to $p_i$. For all $i$ such that $q_i \in [c_i, x_i] \cup [c_i, z_i]$ we have $d(a_i,q_i) \geq t_i$.  So as $d(\mathbf{p}, \mathbf{q}) =1$ and  $\lim_{\omega} t_i/s_i = \infty$, there is a set $\mathcal{S} \subseteq \N$ of $\omega$-measure $1$ such that   $q_i \in [b_i, y_i] \cup [b_i, z_i]$ for all $i \in \mathcal{S}$.  A similar argument to that used earlier shows that $\mathbf{q}$ is a point on $B$, and completes the proof.  

\ms

\ni \ref{hyp} $\implies$ \ref{Dehn pres cond},  \ref{log upper rad}, \ref{log rad}.  There is an account of an elegant proof   due to N.Brady  that \ref{hyp} $\implies$ \ref{Dehn pres cond} in H.Short's notes \cite[Proposition~4.11]{ShortNotes}.   That  \ref{hyp} $\implies$ \ref{log rad} is proved in Dru\c{t}u  \cite{Drutu2} in the course of establishing equivalence of \ref{hyp}, \ref{log rad} and \ref{sublinear rad}.  A proof of \ref{hyp} $\implies$ \ref{log rad} can also be found in \cite[Lemma~6.1]{ShortNotes} -- the idea is to use the thin-triangles condition to find a finite presentation for $G$ to construct suitable van~Kampen diagrams.   The stronger result \ref{hyp} $\implies$ \ref{log upper rad} that the logarithmic radius upper bound is realised  on \emph{all} minimal area van~Kampen diagrams is stated in \cite[\S5.C]{Gromov} and proved in \cite{GR2} using an estimate of the radius of a minimal area diagram in terms of the areas and of the \emph{annular} regions and the lengths of the curves separating them in the decomposition of the proof of Theorem~\ref{set}.  

\ms

\ni \ref{Dehn pres cond} $\implies$ \ref{linear isop}, \ref{linear isop} $\implies$ \ref{subquadratic}, \ \ref{log upper rad} $\implies$ \ref{log rad}, and \ref{log rad} $\implies$ \ref{sublinear rad} are all immediate (modulo standard arguments about changing the finite presentation). 

\ms

\ni \ref{subquadratic} $\implies$ \ref{trees}.  The following argument, combined with \ref{trees} $\implies$ \ref{hyp}, provides an alternative to the better known proofs of \ref{subquadratic} $\implies$ \ref{hyp}, given by Bowditch~\cite{Bowditch},  Ol'shanskii~\cite{Olshanskii2} and Papasoglu~\cite{Papasoglu5}.  The following is a sketch of the proof of Dru\c{t}u~\cite{Drutu} .   We  begin with two lemmas:

\begin{lem}[\textup{Bowditch~\cite{Bowditch}}] \label{quadrangle}
Suppose we express an edge-circuit $\gamma$ in the Cayley 2-complex $\Cay^2(\PP)$ of a finite presentation as a union of four consecutive arcs: $\gamma = \alpha^1 \cup    \alpha^2 \cup  \alpha^3 \cup  \alpha^4$.  Define $d_1:= d(\alpha^1, \alpha^3)$ and  $d_2:= d(\alpha^2, \alpha^4)$ where $d$ is the combinatorial metric on the 1-skeleton of $\Cay^2(\PP)$.   Then $\Area(\gamma) \geq K d_1d_2$, where $K$ is a constant depending only on $\PP$.  
\end{lem} 

\ni Bowditch's lemma can be proved by taking successive star-neighbourhoods $S^i$ of $\alpha_1$ in a $\PP$-van~Kampen diagram $\Delta$ filling $\gamma$.  That is, $S^0= \alpha_1$ and $S^{i+1} = \textup{Star}(S^i)$.  [For a subcomplex $A$ of a complex $B$, the star neighbourhood of $A$ is the union of all the (closed) cells that intersect $A$.]  For some constant $K'$, we have $S^i \cap \alpha_3 = \emptyset$ for all $i \leq K'd_1$, and for all such $i$ the portion of $\partial S^i$ in the interior of $\Delta$ has length at least $d_2$.  The result then follows by summing estimates for the areas of $S^{i+1} \ssm S^i$.   

\ms

\ni Our second lemma is straight-forward.  

\begin{lem} \label{triangle lemma}
If a geodesic space is not a $\mathbb{R}$-tree then it contains a geodesic triangle whose sides only meet at their vertices.  
\end{lem} 

\ni Now assume \ref{subquadratic} and suppose $\textup{Cone}_{\omega}(G, \pmb{1}, \mathbf{s})$ is not an $\mathbb{R}$-tree for some $\mathbf{s}$.  So there is a geodesic triangle $[\mathbf{x}, \mathbf{y}, \mathbf{z}]$  satisfying the condition of Lemma~\ref{triangle lemma}.  Let $l$ be the sum of the lengths of the sides of $[\mathbf{x}, \mathbf{y}, \mathbf{z}]$.  Decompose $[\mathbf{x}, \mathbf{y}, \mathbf{z}]$ into four consecutive arcs $\alpha^1, \alpha^2, \alpha^3, \alpha^4$ such that $d(\alpha^1, \alpha^3), d(\alpha^2, \alpha^4) >  \epsilon$ for some $\epsilon >0$.  

Write $\mathbf{x} =(x_i)$, $\mathbf{y} =(y_i)$ and $\mathbf{z} =(z_i)$ where $x_i,y_i,z_i$ are vertices in $\Cay^2(\PP)$.  The four points $\textbf{a}^j = \alpha^j \cap \alpha^{j+1}$ (indices $j=1,2,3,4$ modulo $4$) are represented by sequences $\textbf{a}^j = (a^j_i)$ of vertices.  Connect up $x_i,y_i,z_i, a^1_i, a^2_i,a^3_i,a^4_i$ (in the cyclic order of the corresponding points $\mathbf{x}, \mathbf{y}, \mathbf{z}, \mathbf{a}^1,  \mathbf{a}^2,  \mathbf{a}^3,  \mathbf{a}^4$ around $[\mathbf{x}, \mathbf{y}, \mathbf{z}]$) by geodesics to make an edge-circuit 
$\gamma_i$ in $\Cay^2(\PP)$.  
Let $w_i$ be the word one reads around $\gamma_i$ and let $l_i := \ell(w_i)$.  Then  $$ N_1 \ := \  \set{ \ i  \  \mid \ \abs{ l - l_i/s_i} < l/2 \ }$$ satisfies $\omega(N_1) =1$.  In particular, $s_i > 2l_i / l$ for all $i \in N_1$.  Also 
\begin{eqnarray*}
 N_2 & := &  \set{ \ i  \  \left| \ d( \alpha_i^1, \alpha_i^3 ) > \epsilon s_i \ \right. } \\ 
 N_3 & := &  \set{ \ i  \  \left| \ d( \alpha_i^2, \alpha_i^4 ) > \epsilon s_i \  \right. }  
\end{eqnarray*}
have $\omega(N_2) = \omega(N_3)  =1$.  So $N := N_1 \cap N_2 \cap N_3$ has measure 1 and so is infinite, and by Lemma~\ref{quadrangle} and the inequalities given above,
$$\Area(w_i) \ \succeq \ (\epsilon s_i)^2 \ \geq \ \epsilon^2 \left( \frac{2l_i}{l} \right)^2  \ \succeq \ {l_i}^2 \ = \  \ell(w_i) ^2$$ for all $i \in N$. This contradicts \ref{subquadratic}.

\ms

\ni \ref{sublinear rad} $\implies$ \ref{trees}.  This is  the final step in our proof of the equivalence of the eight statements in the theorem.  Assume \ref{sublinear rad}.  Let $\mathbf{s} = (s_i)$ be a sequence of real numbers with $s_i \to \infty$.  It is straightforward to check that every geodesic triangle $[\mathbf{x}, \mathbf{y}, \mathbf{z}]$ in  $\textup{Cone}_{\omega}(G, \pmb{1}, \mathbf{s})$ is a tripod: lift to a sequence of geodesic triangles $[x_i, y_i, z_i]$ in $\Cay^1(G,X)$; each admits a van~Kampen diagram with sublinear radius; deduce that each side  of $[\mathbf{x}, \mathbf{y}, \mathbf{z}]$ is in a $0$-neighbourhood of the other two sides.  
\newqed

\bs

\ni Remarkably, weaker conditions than \ref{subquadratic} and \ref{sublinear rad} imply hyperbolicity.  If $\PP$ is a finite presentation for which there exist $L,N, \epsilon$ such that every null-homotopic word $w$ with $$N \ \leq \ \Area(w) \ \leq \ LN$$ has $\Area(w) \leq \epsilon \, \ell(w)^2$, then $\PP$ presents a hyperbolic group \cite{Drutu2, Gromov4, Ol2, Papasoglu3}.  And if there exists $L>0$ such that $\Rad(n) \leq n/73$ for all $n \leq L$ then $\PP$ presents a hyperbolic group \cite{Gromov4, Papasoglu4}; it is anticipated that $8$ could replace $73$ in this inequality and the result still hold \cite{Papasoglu4}.

\bs

\begin{exer} \label{L-delta}
For $x,y,z,t$ in a metric space $(\mathcal{X},d)$, define $\delta(x,y,z;t)$ to be
$$ \max  \{ \ d(x,t)+d(t,y) - d(x,y), \  d(y,t)+d(t,z) - d(y,z), \  d(z,t)+d(t,x) - d(z,x) \ \}.$$
We say that $(\mathcal{X},d)$ enjoys the $L_{\delta}$ property of Chatterji~\cite{Chatterji} if for all $x,y,z \in \mathcal{X}$ there exists $t \in \mathcal{X}$ such that 
$\delta(x,y,z;t) \leq \delta$.
\begin{enumerate}
\item Which of the following spaces are $L_{\delta}$ for some $\delta \geq 0$?
\begin{enumerate}
\item $\Z^n$ with the word metric associated to the presentation $\langle x_1, \ldots, x_n \mid [x_i,x_j], \forall 1\leq i<j \leq n \rangle$. 
\item $\mathbb{R}^n$ with the Euclidean metric. 
\item A hyperbolic group with a word metric.
\end{enumerate}
\item Prove that if $(\mathcal{X},d)$ is an $L_{\delta}$ metric space for some $\delta \geq 0$ then all its asymptotic cones are $L_0$.  (In contrast to Theorem~\ref{R-trees}, the converse is false for general metric spaces;  however it is an open problem for spaces admitting cocompact group actions.) 
\end{enumerate}
 \end{exer}

\section{Groups with simply connected asymptotic cones} \label{simply conn cones}

We now give a characterisation of finitely generated groups with simply connected asymptotic cones. The implication $(\textit{1}) \Rightarrow (\textit{2})$  was
proved by R.~Handel \cite{Handel} and by Gromov~\cite[5.F]{Gromov}; an account is given by Dru\c{t}u~\cite{Drutu}. The reverse implication appears in \cite[page~793]{Papasoglu}. More general arguments in \cite{Riley3} develop those in \cite{Drutu,  Gromov,  Handel, Papasoglu}.

\begin{thm}\label{Partitioning words} 
Let $G$ be a group with finite generating set $X$. Fix any non-principal ultrafilter $\omega$. The following are equivalent.
\begin{enumerate}
    \item The asymptotic cones $\textup{Cone}_{\omega}(G,\mathbf{1},\textbf{s})$ are
   simply connected for all $\textbf{s}=(s_n)$ with $s_n \to \infty$.
    \item Let $\lambda \in (0,1)$.  There exist $K,L \in \Naturals$  such that for all
    null-homotopic words $w$ there is an equality
    \begin{eqnarray}
    w \ = \ \prod^{K}_{i=1}{u_i w_i u^{-1}_i}
    \end{eqnarray}
    in the free group $F(X)$
    for some words $u_i$ and $w_i$
    such that the $w_i$ are null-homotopic
    and $\ell(w_i) \leq \lambda \ell(w) + L$ for all $i$.
     \item Let $\lambda \in (0,1)$.  There exist $K,L \in \Naturals$ such that for all  null-homotopic words $w$ there is a diagram around the  boundary of which reads $w$, and that possesses at most $K$ 2-cells, the boundary circuits of which are labelled by null-homotopic words of length at most $\lambda \ell(w)+L$. 
\end{enumerate}
\end{thm}

\ni \emph{Sketch proof.}
The equivalence of (\emph{2}) and  (\emph{3}) when the $\lambda$'s are each the same is proved in the same way as Lemma~\ref{vK Lemma}.  We leave the task of extending this to the case where the $\lambda$'s differ as an exercise. 

We sketch a proof by contradiction, essentially from \cite{Bridson}, that shows that (\emph{1}) implies  (\emph{3}).     Fix $\lambda \in (0,1)$.  Suppose there are null-homotopic words $w_n$ with $\ell(w_n) \to \infty$ and  such that if $\Delta_n$ is a diagram with $\partial \Delta_n$ labelled by $w_n$ and whose 2-cells have boundaries labelled by null-homotopic words of length at most $\lambda \ell(w)$, then $\Area(\Delta_n) \geq n$.  (Note that, despite the lack of mention of $L$ in this last statement, its negation implies (\emph{3}).)   
Let $\gamma_n: \partial ( [0,1]^2) \to \Cay^1(G,X)$ be an edge-circuit based at the identity that follows a path labelled by $w_n$.  Define $s_n:= \ell(w_n)$ and define $\gamma: \partial ([0,1]^2) \to \textup{Cone}_{\omega}(G,\mathbf{1},\textbf{s})$ by $\gamma(r)= (\gamma_n(r))$.  As $\textup{Cone}_{\omega}(G,\mathbf{1},\textbf{s})$ is simply connected,  $\gamma$ can be extended to a continuous map $\overline{\gamma} : [0,1]^2 \to \textup{Cone}_{\omega}(G,\mathbf{1},\textbf{s})$.  By uniform continuity there exists $\epsilon >0$ such that for all $a,b \in [0,1]^2$, if $d(a,b) < \epsilon$ (in the Euclidean metric) then $d(\overline{\gamma}(a), \overline{\gamma}(b) ) \leq \lambda/16$.  For convenience, we can assume $1/\epsilon$ is an integer.  Consider the images $\mathbf{x}^{i,j}$ under $\overline{\gamma}$ of the  $(1+ 1/\epsilon)^2$ lattice points of the subdivision of $[0,1]^2$ into squares of side $1/\epsilon$.  The images $\mathbf{x}^{i,j}, \mathbf{x}^{i+1,j},\mathbf{x}^{i,j+1}, \mathbf{x}^{i+1,j+1}$ of four corners of a square of side $1/\epsilon$ are represented by sequences $\left(x^{i,j}_n\right), \left(x^{i+1,j}_n\right), \left(x^{i,j+1}_n\right), \left(x^{i+1,j+1}_n\right)$ with $$\lim_{\omega} \frac{d\left(x^{i,j}_n, x^{i+1,j}_n\right) + d\left(x^{i+1,j}_n, x^{i+1,j+1}_n\right) + d\left(x^{i+1,j+1}_n, x^{i,j+1}_n\right) + d\left(x^{i,j+1}_n, x^{i,j}_n\right)}{s_n}   \ \leq \ \frac{\lambda}{4}$$ and hence $$d\left(x^{i,j}_n, x^{i+1,j}_n\right) + d\left(x^{i+1,j}_n, x^{i+1,j+1}_n\right) + d\left(x^{i+1,j+1}_n, x^{i,j+1}_n\right) + d\left(x^{i,j+1}_n, x^{i,j}_n\right) \  < \ \frac{ \lambda \ell(w_n)}{ 2}$$ for all $n$ in a set $\mathcal{S}$ of $\omega$-measure $1$ (and hence infinite).  Moreover $\mathcal{S}$ can be taken to be independent of the choice of square, as there are only finitely many in the subdivision of $[0,1]^2$.  All adjacent $x_n^{i,j}$ can be joined  by geodesics to  give fillings of infinitely many $w_n$ with diagrams with $1/ \epsilon^2$ 2-cells, each with boundary length at most $\lambda \ell(w_n)$. This is the contradiction we seek.   

\ms

Finally, we assume (\emph{3}) and sketch a proof of (\emph{1}).  Fix $\lambda \in (0,1)$ and let $K,L$ be as per (\emph{3}).  Let $s_n$ be a sequence of real numbers with $n \to \infty$.  We will show that $\mathcal{C} := \textup{Cone}_{\omega}(G,\mathbf{1},\textbf{s})$ is simply connected.  We begin by proving that all rectifiable (that is, finite length) loops $\gamma: \partial \mathbb{D}^2 \to \mathcal{C}$  are null-homotopic.  

Assume $\gamma$ is parametrised proportional to arc length.  Inscribe a regular $m$-gon in $\mathbb{D}^2$ with $m>2/\lambda$.  Let $\mathbf{a}^0, \ldots, \mathbf{a}^{m-1}$ be the images under $\gamma$ of the  vertices  of the $m$-gon.  For $i=0, \ldots, m-1$, let $(a^i_n)_{n\in \N}$ be sequences of vertices in $\Cay^1(G, X)$ representing $\mathbf{a}^i$.  For fixed $n$ and for all $i$, join $a^i_n$ to $a^{i+1}_n$ (upper indices modulo $m$) by a geodesic and fill the resulting  edge-circuit as in (\emph{3}).  Roughly speaking, as there are no more than $K$ 2-cells in each filling, one topological configuration occurs for all $n$ in some set of $\omega$-measure $1$, and the fillings converge in $\mathcal{C}$ to a filling of the geodesic $m$-gon by a diagram of that type in which the lengths of the boundaries of the $\leq K$ 2-cells are all at most $\lambda \ell(\gamma)$ (the additive $L$ term disappears in the limit).  Adding in the regions between the $m$-gon and $\gamma$, we have a filling with at most $K + m$ regions each with boundary length at most $\lambda \ell(\gamma)$.        Now iterate the process, refining the filling further, each successive time decreasing the length of the boundaries of the regions by a factor of $\lambda$.  Asymptotic cones are complete metric spaces \cite{Wilkie}, and this fact is used in defining a limit that is a continuous extension of $\gamma$ across $\mathbb{D}^2$.  

\begin{figure}[ht]
\psfrag{p}{$p$}
\psfrag{q}{$q$}
\centerline{\epsfig{file=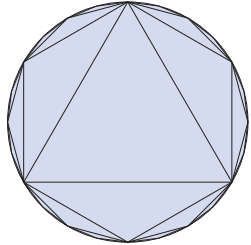}} 
\caption{A tessellation of $\mathbb{D}^2$ by ideal triangles.} \label{tessellation} 
\end{figure}

To show that arbitrary loops $\gamma: \partial \mathbb{D}^2 \to \mathcal{C}$ are null-homotopic, first tessellate $\mathbb{D}^2$ by ideal triangles as illustrated in Figure~\ref{tessellation}, and then extend $\gamma$ to the 1-skeleton of the tessellation by mapping each edge to a geodesic connecting its endpoints.  Each geodesic triangle is a rectifiable loop of length at most three times the diameter of the image of  $\gamma$ (which is finite, by compactness).    So each triangle can be filled in the way already explained, and the result  is an extension of $\gamma$ to a continuous map $\overline{\gamma} : \mathbb{D}^2 \to \mathcal{C}$.  This extension can be proved to be continuous by an argument that uses continuity of $\gamma$ to estimate the perimeter of geodesic triangles in the tessellation whose domains are close to $\partial \mathbb{D}^2$.   Details are in \cite{Riley3}.  
\newqed

\bs

\begin{exer}  
Recall that if $X$ and $X'$ are two finite generating sets for $G$ then $(G,d_X)$ and $(G,d_X')$ are quasi-isometric and so their asymptotic cones are Lipschitz-equivalent.  It follows that conditions (\emph{2}) and (\emph{3}) of Theorem~\ref{Partitioning words}  do not depend on the particular finite generating set $X$.  Prove this result directly, i.e.\ without using asymptotic cones. 
\end{exer}

\medskip
\ni The polynomial bound (\ref{Dehn bound}) of the following theorem is due to Gromov \cite[$5F''_1$ ]{Gromov}. Dru\c{t}u \cite[Theorem~5.1]{Drutu} also provides a proof. The upper bound on $\IDiam$ appears as a remark of Papasoglu at the end of \cite{Papasoglu} and is proved in detail in \cite{Bridson}.  The idea is to repeatedly refine the 2-cells of a diagram as per  Theorem~\ref{Partitioning words}~(3) until the lengths of the boundary loops of the 2-cells are all at most a constant, but there is a technical issue of maintaining planarity during each refinement.  Below, we use the Dehn proof system to avoid this issue (or, to be more honest, push it elsewhere: to the proof of Proposition~\ref{dehn proof prop}, essentially).  

\begin{thm}\label{various bounds} Fix any non-principal ultrafilter $\omega$.  Suppose, for
all sequences of scalars $\textbf{s}$ with $s_n \to \infty$, the asymptotic cones $\textup{Cone}_{\omega}\!\left(G,\mathbf{1},\textbf{s}\right)$
of a finitely generated group $G$ are simply connected.  Then
there exists a finite presentation $\langle X \mid
R \rangle$ for $G$ with respect to which
\begin{eqnarray}
\Area(n) &\preceq& n^{\alpha}, \label{Dehn bound} \\ \IDiam(n)
&\preceq&  n,\label{isodiametric bound} \\
  \FL(n) &\preceq& n,\label{filling length bound}
\end{eqnarray}
where $\alpha = \log_{1/\lambda} K$ and $K$ is as in Theorem~\ref{Partitioning words}.  Further, these bounds are realisable simultaneously.
\end{thm}

\begin{proof}
Fix $\lambda \in (0,1)$ and let $K,L \in \N$ be as in Theorem~\ref{Partitioning words}~(3).  Suppose $w$ is null-homotopic.  Then, for similar reasons to those used in Section~\ref{Dehn proof system}, there is a sequence $(w_i)_{i=0}^m$ with each $w_{i+1}$ obtained from $w_i$ by either a free reduction/expansion or by 
\begin{eqnarray} \label{big r}
w_i \ = \ \alpha \beta \ \mapsto \ \alpha u \beta \ = \ w_{i+1}
\end{eqnarray}
where $u$ is a null-homotopic word with $\ell(u) \leq \lambda \ell(w) + L$.  Moreover, the number of moves of type (\ref{big r}) is at most $K$ and $$\max_i \ell(w_i)  \ \leq \ (K+1) (\lambda n +L) +n.$$  This  bound on $\max_i \ell(w_i)$ holds because the total number of edges in the diagram is at most $K(\lambda n +L) +n$, and a further $\lambda n +L$ is added because the moves of type  (\ref{big r}) insert the whole of $u$ rather than exchanging one part of a relator for another like an application-of-a-relator move.  

Let $R$ be the set of null-homotopic words on $X$ of length at most $1+ (L / (1- \lambda))$.  We will show by induction on $\ell(w)$ that $w$ admits a null-$\PP$-sequence where $\PP = \langle X \mid R \rangle$.  This is clearly the case if  $\ell(w)  \leq  1+ (L / (1- \lambda))$.  Assume $\ell(w) > 1+ (L / (1- \lambda))$.
Then for $u$ as in (\ref{big r})  $$\ell(u) \  \leq \  \lambda \ell(w) + L \ < \ \ell(w),$$ and so by the induction hypothesis and Proposition~\ref{dehn proof prop}, each $u$  admits a null-sequence involving at most $\Area(u)$ application-of-a-relator moves and words of length at most $FL(u)$.   In place of each move  (\ref{big r}) insert such a null-sequence, run backwards, to create the $u$.  The result is a null-$\PP$-sequence for $w$.  Moreover, this null-$\PP$-sequence shows that  
\begin{eqnarray}
\Area(n) & \leq & K \Area(\lambda n +L) \ \ \textup{ and} \label{area ind} \\
\FL(n)  & \leq & \FL(\lambda n +L) +  (K+1) (\lambda n +L) + n. \label{FL ind}
\end{eqnarray}
Let $k = 1 + \lfloor \log_{\lambda} (1/n) \rfloor$.  Repeatedly apply (\ref{area ind}) until, after $k$ iterations, we have $\Area(n)$ in terms of 
$$\Area \left( \lambda^k n + \frac{L}{1- \lambda} \right).$$ But this is $1$ (or $0$) because the argument is at most  $1+ (L / (1- \lambda))$.  So $$\Area(n) \ \leq \ K^k \ \leq \ K^{1- \log_{\lambda} n} \ = \ K \, n^{\log_{1/\lambda} K},$$ which proves (\ref{Dehn bound}).  
We get (\ref{filling length bound}) by summing a geometric series arising from iteratively applying (\ref{FL ind}), and (\ref{isodiametric bound}) then follows.  
\end{proof}

\ms

\begin{example}
Finitely generated groups $(G,X)$ satisfying the $L_{\delta}$ condition of Exercise~\ref{L-delta} enjoy conditions (\emph{2}) and (\emph{3}) of Theorem~\ref{Partitioning words} with $K=3$, with $L =  \delta$ and $\lambda =2/3$.  So all the asymptotic cones of $G$ are simply connected, $G$ is finitely presentable, and 
$$\Area(n) \ \preceq \  n^{\log_{3/2}{3}}, \ \ \ \FL(n) \ \preceq \ n, \ \ \ \textup{and} \ \ \ \IDiam(n) \ \preceq \ n.$$   
This upper bound on Dehn function is subcubic ($\log_{3/2}{3} \simeq 2.71$) and is due to Elder~\cite{Elder2}.  It is an open problem to find a group satisfying the $L_{\delta}$ condition for some $\delta$ but not having $\Area(n) \preceq n^2$.    
\end{example}
	
\ms

\ni Having Dehn function growing at most polynomially fast does not guarantee that the asymptotic cones of a group $G$ will all be simply connected: there are groups \cite{Bridson, SBR} with such Dehn functions but with $\IDiam$ growing faster than linearly.  Indeed, Ol'shanskii and Sapir have recently constructed a group with with $\IDiam(n) \simeq n$ and $\Area(n) \simeq n^3$ but no asymptotic cones simply connected  \cite{OS2}:
$$\langle \ a, b, c, k \ \mid \  [a,b], \, [a,c], \,
k^b = ka, \, k^c=ka \ \rangle.$$
They claim that $S$-machines from \cite{OS3} could be be used to obtain an example with $\IDiam(n) \simeq n$ and $\Area(n) \simeq n^2 \log n$.  So  the following theorem of Papasoglu \cite{Papasoglu} is near sharp. 

\begin{thm} If a group has quadratic Dehn function then its asymptotic cones are all simply connected.  
\end{thm}

\ni Papasoglu's proof shows that groups with quadratic Dehn function have linear $\IDiam$ and it then proceeds to the criterion in Theorem~\ref{Partitioning words}.  In particular, by Theorem~\ref{various bounds},  one can deduce (without using asymptotic cones) that groups with quadratic Dehn functions have linear $\FL$, and this seems a non-trivial result.  More generally, it is shown in \cite{GR1} that groups with $ \Area(n) \preceq n^{\alpha}$ for some $\alpha \geq 2$ have $\IDiam(n) \preceq n^{\alpha -1}$,  and in \cite[Theorem~8.2]{GR4} the additional conclusion $\GL(n) \preceq n^{\alpha -1}$ is drawn.   

\ms

\begin{exer}
Give a direct proof that finite presentations with $\Area(n) \preceq n^2$ have $\IDiam(n) \preceq n$, and (\emph{harder}) $\FL(n) \preceq n$.
\end{exer}

\section{Higher dimensions} 

The notions of finite generability and finite presentability are $\mathcal{F}_1$ and $\mathcal{F}_2$ in a family of \emph{finiteness conditions} $\mathcal{F}_n$ for a group $G$.  We say $G$ is of \emph{type $\mathcal{F}_n$} if there is a $K(G,1)$ which is a complex with finite $n$-skeleton.  (A $K(G,1)$ is a space with fundamental group $G$ and all other homotopy groups trivial.)
For example, if $G$ has finitely presentation $\PP$ then  one can attach 3-cells, then 4-cells, and so on, to the presentation 2-complex of $\PP$ in such a way as to kill off $\pi_2$, then $\pi_3$, and so on, of its universal cover, producing a $K(G,1)$ with finite 1- and 2-skeleta.      

Suppose a group $G$ is of type $\mathcal{F}_n$  for some $n \geq 2$.  Then there are reasonable notions of higher dimensional filling inequalities concerning filling $(n-1)$-spheres with $n$-discs in some appropriate sense.  The following approach of Bridson~\cite{Bridson2} is a geometric/combinatorial interpretation of the algebraic definitions of Alonso, Pride \& Wang \cite{AWP}.   By way of warning we remark that there is no clear agreement between this definition and that of Epstein et al.\ \cite{Epstein}, that concerns spheres in a Riemannian manifold on which $G$ acts properly, discontinuously and cocompactly.  

Starting with a finite presentation  for $G$ build a \emph{presentation $(n+1)$-complex} by  attaching finitely many $3$-cells to kill $\pi_2$ of the presentation 2-complex, and then  finitely many $4$-cells to kill $\pi_3$, and so on up to  $(n+1)$-cells to kill $\pi_n$.  The data involved in this construction is termed a \emph{finite $(n+1)$-presentation} $\PP_{n+1}$.  Let $\Cay^{n+1}(\PP_{n+1})$ denote the universal cover of the presentation $(n+1)$-complex of $\PP_{n+1}$.  An important technicality here is that the $(j+1)$-cells may be attached not by combinatorial maps from their boundary combinatorial $n$-spheres, but by \emph{singular combinatorial maps}.  These can collapse $i$-cells, mapping them into the $(i-1)$-skeleton of their range -- see \cite{Bridson2} for more details.

Consider a singular combinatorial map $\gamma: S^{n} \to \Cay^{n+1}(\PP_{n+1})$, where $S^{n}$ is some combinatorial $n$-sphere. We can \emph{fill} $\gamma$ by giving a singular combinatorial extension $\bar{\gamma}: D^{n+1} \to \Cay^{n+1}(\PP_{n+1})$ with respect to some combinatorial  $(n+1)$-disc $D^{n+1}$ such that $S^n=\partial D^{n+1}$. We define $\Vol_n({\gamma})$ and $\Vol_{n+1}(\bar{\gamma})$ to be the number of $n$-cells $e$ in $S^n$ and $D^{n+1}$, respectively, such that $\gamma \restricted{e}$ is a homeomorphism, and we define the \emph{filling volume} $\FVol(\gamma)$ of $\gamma$ to be the minimum amongst all $\Vol_{n+1}(\bar{\gamma})$ such that $\bar{\gamma}$ fills $\gamma$.  This leads to the definition of the $n$-th order Dehn function, which, up to  $\simeq$ equivalence, is independent of the choice of finite $(n+1)$-presentation \cite{AWP}. 

We will be concerned with inequalities concerning not just filling volume but also  (intrinsic) diameter.  So define $\Diam_{n}(\gamma)$ and $\Diam_{n+1}(\bar{\gamma})$ to be the maxima of the  distances between two vertices on $S^n$ or $D^{n+1}$, respectively, in the combinatorial metric on their 1-skeleta.  And define $\FDiam(\gamma)$ to be the minimum of $\Diam_{n+1}(\bar{\gamma})$ amongst all $\bar{\gamma}$ filling $\gamma$.   

Recall from Section~\ref{combable section} how we coned off a loop (a 1-sphere) as in Figure~\ref{cockleshell} to get an upper bound on the Dehn function.  In the same way it is possible to cone off a singular combinatorial $n$-sphere and then fill the \emph{rods} (the cones over each of the $\Vol_n(\gamma)$ non-collapsing $n$-cells) to prove the following generalisation of Theorem~\ref{Combable groups}.

\begin{thm}  \textup{\cite{Epstein, Gersten5}}  \label{higher diml combing}
Suppose $G$ is a finitely generated, asynchronously combable group with length function $k \mapsto \textup{L}(k)$. Then $G$ is of type $\mathcal{F}_{n}$ for all $n$. Further, given a finite $(n+1)$-presentation $\PP_{n+1}$ for $G$, every
singular combinatorial $n$-sphere $\gamma: S^n \to \Cay^{n+1}(\PP_{n+1})$ can be filled by some singular combinatorial $(n+1)$-disc $\bar{\gamma}: D^{n+1} \to \Cay^{n+1}(\PP_{n+1})$ with 
\begin{eqnarray*}
\FVol_{n+1}(\gamma) & \preceq &  \Vol_n(\gamma) \, \textup{L}(\Diam_n(\gamma)), \\ 
\FDiam_{n+1}(\gamma) &\preceq&   \textup{L}(\Diam_n(\gamma)).
\end{eqnarray*}
\ni Moreover, these bounds are realisable simultaneously.
\end{thm}

\ni Similarly, by coning off and then filling the \emph{rods} it is possible to obtain filling inequalities for groups with $n$-connected asymptotic cones.   (A space is \emph{$n$-connected} if its homotopy groups $\pi_0, \pi_1, \ldots, \pi_n$ are all trivial.)
Each rod has filling volume that is at most polynomial in its diameter  for reasons similar to those that explain the appearance of the polynomial area bound in  Theorem~\ref{various bounds}.

\begin{thm} \textup{\cite{Riley3}} \label{Mainthm about groups with highly conn cones}
Suppose $G$ is a finitely generated group whose asymptotic cones are all
$n$-connected ($n \geq 1$). Then $G$ is of type
$\mathcal{F}_{n+1}$ and, given any finite $(n+1)$-presentation $\PP_{n+1}$ for $G$, every
singular combinatorial $n$-sphere $\gamma: S^n \to \Cay^{n+1}(\PP_{n+1})$ can be filled by some singular combinatorial $(n+1)$-disc $\bar{\gamma}: D^{n+1} \to \Cay^{n+1}(\PP_{n+1})$ with
\begin{eqnarray*}
\FVol_{n+1}(\gamma) & \preceq &  \Vol_n(\gamma) \, (\Diam_n(\gamma))^{\alpha_n}, \\ 
\FDiam_{n+1}(\gamma) &\preceq&   \Diam_n(\gamma),
\end{eqnarray*}
for some $\alpha_n$ depending only on $\PP_{n+1}$.  Moreover, these bounds are realisable simultaneously.
\end{thm}

\begin{Open Question}
Do the higher order Dehn functions of hyperbolic groups admit linear upper bounds?
\end{Open Question}

\ni Perhaps the characterisation of hyperbolic groups in terms of $\mathbb{R}$-trees in Theorem~\ref{hyperbolicity via Dehn fns} can be used to resolve this question.  Mineyev~\cite{Mineyev} gets linear upper bounds on the volumes of \emph{homological} fillings.

\begin{Open Question} 
Is there a sequence $G_n$ of groups (discrete if possible) such that the asymptotic cones of $G_n$ are all $n$-connected but not all $(n+1)$-connected?
\end{Open Question}

\ni Gromov \cite[$\S2.B_1$]{Gromov} makes the tantalising suggestion that $\SL_n(\mathbb{Z})$ might provide a family of examples.   Epstein \& Thurston  \cite[Chapter~10]{Epstein} show, roughly speaking, that any $(n-2)$-st order isoperimetric function for  $\SL_n(\mathbb{Z})$  is at least exponential and hence, by Theorem~\ref{Mainthm about groups with highly conn cones}, the asymptotic cones of $\textup{SL}_n(\mathbb{Z})$ are not all $(n-2)$-connected.  It remains to show that all the asymptotic cones of $\textup{SL}_n(\mathbb{Z})$  are $(n-3)$-connected. 

In \cite[$\S2.\textup{B}.$(f)]{Gromov} Gromov outlines a strategy for showing that  a certain sequence of solvable Lie groups $S_n$ will have every asymptotic cone $(n-3)$-connected but of non-trivial  (uncountably generated, even) $\pi_{n-2}$.  The reason one expects the non-triviality of $\pi_{n-2}$ is that  similar arguments to those in \cite{Epstein} should show the $(n-2)$-nd order isoperimetric function of $S_n$ again to be at least exponential.

Stallings gave the first example of a finitely presented group that is not of type $F_3$ \cite{Stallings}. Bieri generalised this to a family of groups $\textup{SB}_n$ of type $F_{n-1}$ but not $F_n$ \cite{Bieri}.  These groups are discussed in \cite{BradyNotes}. Their asymptotic cones are ripe for investigation.

\bibliographystyle{plain}
\bibliography{bibli}
\addcontentsline{toc}{chapter}{Bibliography}



\end{document}